\documentclass[leqno, notitlepage]{amsart}

\usepackage{amssymb}
\usepackage[all, poly]{xy}

\DeclareMathOperator{\im}{im}
\DeclareMathOperator{\re}{re}
\DeclareMathOperator{\Hom}{Hom}
\DeclareMathOperator{\Ext}{Ext}
\DeclareMathOperator{\Ends}{Ends}
\DeclareMathOperator{\In}{In}
\DeclareMathOperator{\Out}{Out}
\DeclareMathOperator{\Ker}{Ker}
\DeclareMathOperator{\Coker}{Coker}
\DeclareMathOperator{\Deg}{deg}
\DeclareMathOperator{\Ob}{Ob}

\DeclareMathOperator{\Reg}{reg}
\DeclareMathOperator{\Irr}{irr}
\DeclareMathOperator{\Id}{Id}
\DeclareMathOperator{\Supp}{Supp}

\DeclareMathOperator{\Ind}{Ind}

\DeclareMathOperator{\ad}{ad}

\DeclareMathOperator{\Spec}{Spec}

\DeclareMathOperator{\rank}{rank}

\newcommand{\modN}{\:(\text{mod}\: N)}

\swapnumbers
\theoremstyle{plain} 
\newtheorem*{theorem}{Theorem}
\newtheorem*{corollary}{Corollary}

\newtheorem*{proposition}{Proposition}

\newtheorem{numproposition}[subsubsection]{Proposition}

\numberwithin{equation}{subsubsection}

\numberwithin{enumi}{subsubsection}

\begin{document}

\title{ Affine Lie algebras and tame quivers }

\author{ Igor Frenkel }
\address{ Department of Mathematics, Yale University,
P.O. Box 208283, New Haven, CT 06520-8283 }
\email{frenkel@math.yale.edu}

\author{ Anton Malkin }
\address{ Department of Mathematics, Yale University,
P.O. Box 208283, New Haven, CT 06520-8283 }
\email{malkin@math.yale.edu}

\author{ Maxim Vybornov }
\address{ Department of Mathematics, Duke University,
P.O. Box 90320, Durham, NC 27708-0320 }
\email{mv@math.duke.edu}

\date{May 11, 2000}

\maketitle

\setcounter{section}{-1}

\tableofcontents

\section{Introduction}

\subsubsection{}
There is a remarkable connection between the theory
of representations of quivers and the structure theory
of Lie algebras. The first manifestation of
this connection was discovered by P. Gabriel 
\cite{Gabriel1972}.
Let $Q$ be a quiver obtained by orienting edges
of the Dynkin graph corresponding to 
a simple simply laced Lie algebra $\mathfrak{g}$.
Gabriel proved that the set $\mathcal{T}$ of
isomorphism classes of indecomposable 
complex representations of
$Q$ is in one-to-one correspondence with the set 
$R_+$ of positive roots of $\mathfrak{g}$.

Gabriel's result was soon extended by J. Bernstein, 
I. Gelfand, and V. Ponomarev 
\cite{BernsteinGelfandPonomarev}, who
introduced reflection functors $\mathcal{S}_i$
corresponding to Coxeter generators of the Weyl group 
of $\mathfrak{g}$. The reflection functor
$\mathcal{S}_i$ acts from the category
$\mathcal{M} (Q)$ of complex 
representations of $Q$ to the category
$\mathcal{M} (Q')$, where $Q'$ differs from $Q$
only by orientation. Using the reflection functors
Bernstein, Gelfand and Ponomarev were able
to give another proof of the Gabriel theorem. 

\subsubsection{}\label{IntroFinite}
The set $R_+$ of positive roots corresponds to a
basis of a maximal nilpotent subalgebra 
$\mathfrak{n}$ of $\mathfrak{g}$, and one might guess 
that there exists an intrinsic Lie bracket
on the $\mathbb{C}$-linear span of the set 
$\mathcal{T}$,
such that the resultant Lie algebra is isomorphic
to $\mathfrak{n}$. However this Lie
bracket was introduced only 18 years later by 
C. M. Ringel \cite{Ringel1988, Ringel1990}.
Ringel actually considers
representations of $Q$ over finite fields rather
than over complex numbers, so we use a variant of 
his definition
due to A. Schofield \cite{Schofield} and 
G. Lusztig \cite{Lusztig1991b}.

Let $[\mathbf{P}_{\alpha}]$, $[\mathbf{P}_{\beta}] 
\in \mathcal{T}$
be isomorphism classes of
indecomposable representations of $Q$. Then their
Lie bracket is defined as follows
\begin{equation}\label{RingelLieBracket}
[[\mathbf{P}_{\alpha}] , 
[\mathbf{P}_{\beta}]] =
\sum_{[\mathbf{P}_{\gamma}]\in\mathcal{T}}
(\chi (N_{\mathbf{P}_{\alpha} , 
\mathbf{P}_{\beta} ; \mathbf{P}_{\gamma}})-
\chi (N_{\mathbf{P}_{\beta} , 
\mathbf{P}_{\alpha} ; \mathbf{P}_{\gamma}}))
[\mathbf{P}_{\gamma}] ,
\end{equation}
where 
$\chi (N_{\mathbf{P}_{\alpha} , 
\mathbf{P}_{\beta} ; \mathbf{P}_{\gamma}})$
is the Euler characteristic with compact support
of the algebraically constructible set
$N_{\mathbf{P}_{\alpha} , \mathbf{P}_{\beta} ; 
\mathbf{P}_{\gamma}}$
of all subrepresentations 
$\mathbf{V}\subset\mathbf{P}_{\gamma}$
such that $\mathbf{V}$ is isomorphic to
$\mathbf{P}_{\alpha}$, and 
$\mathbf{P}_{\gamma}/\mathbf{V}$
is isomorphic to $\mathbf{P}_{\beta}$.

Thus we obtain a complex Lie algebra denoted by
$\mathfrak{n}^{\ast}$ with a distinguished 
basis $\{ [\mathbf{P}] \}_{[ \mathbf{P}] \in 
\mathcal{T}}$.
In \cite{Ringel1990} Ringel proved that 
$\mathfrak{n}^{\ast}$ and $\mathfrak{n}$ are
isomorphic as $R_+$-graded complex Lie algebras, 
and,
moreover, he was able to find the structure 
constants
in $\mathfrak{n}^{\ast}$. Namely, let 
$\mathbf{P}_{\alpha}$
be the unique up to an isomorphism indecomposable 
representation
of $Q$ corresponding to a root $\alpha \in R_+$. 
Then given $\alpha$, $\beta \in R_+$ one has 
\begin{equation}\label{IntroAstBracket}
[[\mathbf{P}_{\alpha}] ,
[\mathbf{P}_{\beta}]] =
\begin{cases}
\epsilon (\alpha, \beta)
[\mathbf{P}_{\alpha + \beta}] &
\text{if $\alpha + \beta \in R_+$ ,}\\
0 &
\text{if $\alpha + \beta \notin R_+$ },
\end{cases}
\end{equation}
where $\epsilon$ is a bimultiplicative 
two-cocycle on the root lattice 
of $\mathfrak{g}$, uniquely defined
by its values on pairs of simple roots:
$\epsilon (i,i)=-1$ for any simple root
$i$, $\epsilon (i,j)=-1$ if
there is an arrow $i\rightarrow j$
in $Q$, $\epsilon (i,j)=1$
otherwise. The choice of bimultiplicative
cocycle $\epsilon$ corresponds
to a choice of orientation of edges which
converts a Dynkin graph into quiver.

It is interesting to note that the choice of  
Chevalley basis in $\mathfrak{n}$ and,
more generally, in $\mathfrak{g}$, such that
structure constants are given by the cocycle
$\epsilon$, was first introduced in relation
to a vertex operator construction of the affine
Lie algebra $\Hat{\mathfrak{g}}$ associated 
to $\mathfrak{g}$ (see \cite{FrenkelKac1980, 
Segal1981}).

Let us also remark that there is an equivalent
definition of the cocycle $\epsilon$ in 
terms of the category $\mathcal{M} (Q)$. Namely,
\begin{equation}\label{IntroHomExt}
\epsilon (\alpha, \beta) =
(-1)^{(\dim_{\mathbb{C}}\Hom_{\mathcal{M} (Q)} 
(\mathbf{P}_{\alpha}, \mathbf{P}_{\beta})
-
\dim_{\mathbb{C}}\Ext^1_{\mathcal{M} (Q)} 
(\mathbf{P}_{\alpha}, \mathbf{P}_{\beta}))}.
\end{equation}

Ringel's proof of \eqref{IntroAstBracket}
is based on a case-by-case study of all 
possible varieties 
$N_{\mathbf{P}_{\alpha} , 
\mathbf{P}_{\beta} ; \mathbf{P}_{\gamma}}$
and is rather long.
In this paper we give a new and short proof of
the Ringel theorem using reflection functors
of Bernstein, Gelfand, and Ponomarev. Instead of
studying varieties 
$N_{\mathbf{P}_{\alpha} , 
\mathbf{P}_{\beta} ; \mathbf{P}_{\gamma}}$
and calculating their Euler characteristics
for a particular quiver $Q$ we
consider all quivers with the same
underlying Dynkin graph and use 
functorial properties of the Lie algebra
$\mathfrak{n}^*$.

More explicitly, given a quiver $Q$ we define 
a complex Lie algebra $\mathfrak{n}^{\epsilon}$. 
As a linear
space $\mathfrak{n}^{\epsilon}$ has a basis 
$\{ \Tilde{e}_{\alpha} \}_{\alpha \in R_+}$
corresponding to isomorphism classes of 
indecomposable 
representations of $Q$. The Lie bracket 
is defined as
follows
\begin{equation}\nonumber
[\Tilde{e}_{\alpha} ,
\Tilde{e}_{\beta}] =
\begin{cases}
\epsilon (\alpha, \beta)
\Tilde{e}_{\alpha + \beta} &
\text{if $\alpha + \beta \in R_+$ ,}\\
0 &
\text{if $\alpha + \beta \notin R_+$ }.
\end{cases}
\end{equation}
Because of \eqref{IntroHomExt}
$\mathfrak{n}^{\epsilon}$ is functorial with
respect to  reflection functors.
The same is true for 
$\mathfrak{n}^{\ast}$ by definition.

It turns out that there are enough 
reflection functors
to ensure that a Lie algebra 
with a basis parameterized by 
isomorphism classes of indecomposable
representations of $Q$ and with
functorial structure constants
is unique. 
More precisely, we introduce a homomorphism
$\Xi : \mathfrak{n}^{\epsilon}
\rightarrow \mathfrak{n}^{\ast}$ 
given on generators by
$\Xi ( \Tilde{e}_i ) = [ \mathbf{P}_i ]$ 
for any simple root $i$. Then we use reflection 
functors to prove
that $\Xi ( \Tilde{e}_{\alpha} ) = [ 
\mathbf{P}_{\alpha} ]$
for any $\alpha \in R_+$,
which is equivalent to the Ringel theorem.

\subsubsection{}
Shortly after appearance of Gabriel's paper
L. A. Nazarova \cite{Nazarova1973} and independently 
P. Donovan and M. R. Freislich 
\cite{DonovanFreislich1973} classified  
indecomposable
representations of quivers
associated to simply laced
extended Dynkin graphs, which we call 
quivers of affine type (they
are also called tame quivers).

For quivers of affine type the bijection 
between isomorphism classes of indecomposable 
representations
and positive roots of the corresponding affine Lie
algebra $\mathfrak{g}$ holds only for positive
real roots $\alpha \in R^{\re}_+$. In the case of a
positive imaginary
root $\alpha \in R^{\im}_+$ there exists
an uncountable family $\mathcal{T}_{\alpha}$
of non-isomorphic indecomposable representations
corresponding to $\alpha$. According to Dlab
and Ringel \cite{DlabRingel} the set 
$\mathcal{T}_{\alpha}$
does not depend on $\alpha \in R^{\im}_+$ and 
admits a surjection
\begin{equation}\nonumber
\mu : \mathcal{T}_{\alpha}
\rightarrow
\mathbb{CP}^1
\end{equation}
such that $\mu$ is injective for all points of
$\mathbb{CP}^1$ except for $L\leq 3$ exceptional points
for which fibers are finite sets.

Since $\mathcal{T}_{\alpha}$ is an infinite set for
$\alpha \in R^{\im}_+$ one has to adjust the
definition of $\mathfrak{n}^{\ast}$. One way
to do it, proposed by Ringel \cite{Ringel1990p, 
Ringel1993},
consists of replacing representations of $Q$
by formal composition series. We adopt another 
approach
due to A. Schofield \cite{Schofield} and
G. Lusztig \cite{Lusztig1991a, Lusztig1991b}. 
Instead of formal linear
combinations of elements of the set $\mathcal{T}$ of
indecomposable representations we consider complex 
valued constructible 
(in some sense) functions on $\mathcal{T}$. 
One can generalize 
the definition of Lie bracket \eqref{RingelLieBracket}.
Finally we define 
$\mathfrak{n}^*$ to be the Lie subalgebra
of the Lie algebra of constructible functions 
on $\mathcal{T}$ generated by characteristic
functions of simple representations.

The main result of the present paper is a 
generalization of the Ringel theorem
described in \ref{IntroFinite} to affine case.
We show, in particular, that the Lie algebra
$\mathfrak{n}^{\ast}$ is isomorphic to the
positive part $\mathfrak{n}$ of the affine
Lie algebra $\mathfrak{g}$ corresponding to 
the quiver $Q$. The Lie algebra $\mathfrak{n}$
is graded by the set $R_+$ of positive roots
and is the precise analogue of a maximal nilpotent 
subalgebra of a finite dimensional simple Lie 
algebra.

Our results imply that the Lie algebra
$\mathfrak{n}^{\ast}$ contains the characteristic 
function
of indecomposable representation
$\mathbf{P}_{\alpha}$ corresponding to any 
positive real root $\alpha$, and we can find
the structure constants for the Lie bracket of
two such characteristic functions.
Namely, if $\alpha$, $\beta$, $\alpha + \beta
\in R^{\re}_+$ then
\begin{equation}\nonumber
[[\mathbf{P}_{\alpha}] ,
[\mathbf{P}_{\beta}]] =
\epsilon' (\alpha, \beta)
[\mathbf{P}_{\alpha + \beta}]. 
\end{equation}
Here abusing notation we denote 
by $[\mathbf{P}_{\alpha}]$ the
characteristic function of 
$[\mathbf{P}_{\alpha}] \in \mathcal{T}$.
The structure constants 
$\epsilon'$ are related to
the bimultiplicative two-cocycle
$\epsilon$ defined in \ref{IntroFinite}
via a twist by coboundary
\begin{equation}\nonumber
\epsilon' (\alpha, \beta)=
\epsilon (\alpha, \beta) 
\xi (\alpha +\beta)
\xi^{-1} (\alpha) 
\xi^{-1} (\beta) ,
\end{equation}
where $\xi ( \alpha )=
(-1)^{1+\dim_{\mathbb{C}} \Hom_{\mathcal{M} (Q)} 
(\mathbf{P}_{\alpha}, \mathbf{P}_{\alpha})}$. Note
that in the case of a quiver with underlying
Dynkin graph  
$\Hom_{\mathcal{M} (Q)} 
(\mathbf{P}_{\alpha}, \mathbf{P}_{\alpha})
=\mathbb{C}$ for any indecomposable
$\mathbf{P}_{\alpha}$ and, therefore, the twist
does not affect Ringel's structure constants.

The case of an imaginary root $\alpha \in R^{\im}_+$
is more involved. We give a description of
the root space $\mathfrak{n}^{\ast}_{\alpha}$ in terms
of the surjection
$\mu : \mathcal{T}_{\alpha}
\rightarrow \mathbb{CP}^1$. Namely,
$\mathfrak{n}^{\ast}_{\alpha}$
consists of all functions $f$ on $\mathcal{T}_{\alpha}$
such that
\begin{equation}\label{IntroConst}
\sum_{[\mathbf{P}] \in \mu^{-1} (z)}
f ([\mathbf{P}]) \quad 
\text{ does not depend on $z \in \mathbb{CP}^1$ }.
\end{equation}
The following remarkable identity guarantees that
$\mathfrak{n}^{\ast}_{\alpha}$ has the 
correct dimension
\begin{equation}\nonumber
\sum_{z \in \mathbb{CP}^1} (N_z - 1) =
\rank \mathfrak{g} - 2 ,
\end{equation}
where $N_z$ is the order of the fiber
$\mu^{-1} (z)$.

Considering constructible functions on
$\mathcal{T}$ with values in $\mathbb{Z}$
satisfying condition
\eqref{IntroConst}
on $\mathcal{T}_{\alpha}$ for any $\alpha \in R^{\im}_+$
one obtains a lattice 
$\aleph^{\ast} \in \mathfrak{n}^{\ast}$. Our
results imply
that $\aleph^{\ast}$ is closed with respect to
the Lie bracket in $\mathfrak{n}^{\ast}$
(i.e. it is an order). We also prove
that $\aleph^{\ast}$ is characterized by the property
that it is the minimal order in $\mathfrak{n}^{\ast}$
containing simple root generators. In particular,
the lattice $\aleph^{\ast}_{\alpha}$ for any imaginary
root $\alpha$ is isomorphic to the root lattice of
a simple Lie algebra $\mathfrak{g}_0$, such that 
$\mathfrak{g} = \Hat{\mathfrak{g}}_0$.

Our strategy in affine case is similar to
the one we use in the proof
of the Ringel theorem. We define some ad hoc Lie algebra 
$\mathfrak{n}^{\epsilon}$ using the cocycle $\epsilon$,
and a homomorphism
$\Xi : \mathfrak{n}^{\epsilon}
\rightarrow
\mathfrak{n}^{\ast}$. Then we use functorial
properties of $\mathfrak{n}^{\epsilon}$ and
$\mathfrak{n}^{\ast}$ to describe $\Xi$. In the affine
case, however, the set of reflection functors is 
not enough to fix $\Xi$, and we employ additional 
functors $\mathcal{C}_z$ introduced by Dlab and Ringel. 
A functor $\mathcal{C}_z$ corresponds to a 
point $z \in \mathbb{CP}^1$ and is a full, faithful, 
exact
functor from $\mathcal{M} (C_{N_z})$ to
$\mathcal{M} (Q)$, where $C_{N_z}$ is the cyclic
quiver with $N_z$ vertices. Both
$\mathfrak{n}^{\ast}$ and $\mathfrak{n}^{\epsilon}$
are functorial with respect to functors
$\mathcal{C}_z$ by construction. Using the reflection 
functors and
the functors $\mathcal{C}_z$ we are able to describe 
the map
$\Xi$ for any real root space and a codimension one
subspace in any imaginary root space. The proof is 
completed
with analysis of the positive part of an affine Lie 
algebra
$\Hat{sl}_2$ embedded into $\mathfrak{g}$.

\subsubsection{}
Our description of the isomorphism of the Lie algebra
$\mathfrak{n}^{\ast}$ constructed via representation 
theory of
a quiver of affine type with the positive part 
$\mathfrak{n}$ of the
corresponding affine Lie algebra reveals a certain 
fine structure
of affine Lie algebras and their root systems. 
It turns out that
for any orientation of edges of the extended 
Dynkin graph
one can canonically associate $L \leq 3$ affine 
Lie algebras
of type $A^{(1)}_n$ embedded into
$\mathfrak{g}$. In $D^{(1)}_n$ and $E^{(1)}_n$ cases
these subalgebras correspond precisely to
the connected components of the non-extended 
Dynkin graph 
with the branching vertex removed. In particular, 
it follows that
$L=3$ in $DE$ case.

If orientation of $Q$ is such that each vertex 
is either
a sink or a source, then one can interpret 
representations of
$Q$ in terms of a finite subgroup $\Gamma$  
of $SL (2, \mathbb{C} )$ associated to $Q$
via McKay correspondence \cite{Lusztig1992}.
In this case the type $A^{(1)}_n$ subalgebras 
arise from maximal cyclic
subgroups of $\Gamma$ and their positive roots 
correspond to representations of $\Gamma$ induced from
the cyclic subgroups. Thus our construction
of the positive part $\mathfrak{n}$ 
of an affine Lie algebra $\mathfrak{g}$
can be viewed as a far-reaching expansion
of the McKay correspondence between group
$\Gamma \subset SL (2, \mathbb{C} )$ and
extended Dynkin graph of $\mathfrak{g}$.

\subsubsection{}\label{IntroConclusion}
We conclude Introduction with remarks
concerning some generalizations and
interpretations of our results.

First let us note that using species instead of
quivers one can extend results of the paper
to all simple and affine types of Lie algebras.
Though we consider only simply laced case
in order to preserve clarity of the exposition,
the statements and the proofs
can be generalized almost word-by-word
to species.

Second we remark that Ringel's construction
of the nilpotent Lie algebra $\mathfrak{n}$ was 
extended by 
L. Peng and J. Xiao \cite{PengXiao}
to the whole simple Lie algebra $\mathfrak{g}$
via the root category 
$\mathcal{R} (Q) = D^b \mathcal{M} (Q) / T^2$. 
Our proof of the Ringel theorem 
can be extended to this setting. Moreover 
it would become
more transparent as the reflection functors
act more naturally in the root category
than in abelian category $\mathcal{M} (Q)$.
It would be very interesting to extend the 
construction of
$\mathfrak{n}^{\ast}$ to the root category 
of a quiver
of affine type.

Let us finally mention that results and 
techniques used in this paper could be
stated in the language of algebraic stacks.
For example, description \eqref{IntroConst} 
of imaginary root
subspaces of $\mathfrak{n}^{\ast}$ indicates that
the right setting for $\mathfrak{n}^{\ast}$ is
that of cohomology of the set $\mathcal{T}$ 
considered
as a stack. The details will be provided elsewhere.

\subsubsection*{Acknowledgements}\label{ack} 
The authors are grateful 
to Mikhail Khovanov for his careful reading of 
the manuscript and many valuable comments.
This research was supported in part by NSF grant 
DMS-9700765.

\section{Cartan datum, Lie algebras, 
and quivers}\label{Prelim}

Throughout the paper $\mathbb{Z}$ 
and  $\mathbb{C}$
denote, respectively, integer 
numbers
and complex numbers;
$\mathbb{Z}_{+} = 
\{ k \in \mathbb{Z}, k\geq 0 \}$.

\subsection{Cartan datum and Dynkin graph}

\subsubsection{}
A \emph{Cartan datum} is a pair 
$(I,<,>)$,
consisting of a
finite set $I$ and a bilinear form $<,>$ on
the free abelian group $\mathbb{Z} [I]$, 
with values 
in $\mathbb{Z}$. The bilinear form should 
satisfy the
following conditions:
\begin{align}\nonumber
\label{a} <i,i>& =2 \text{ for all } i \in I ,  
\\\nonumber
<i,j> & \leq 0 \text{ for all } i \neq j ,
\\\nonumber
<i,j> & =  0 \text{ if } <j,i> = 0 .
\end{align}  

The matrix $a_{ij}=<i,j>$ is called the 
Cartan matrix of the 
Cartan datum $(I,<,>)$.

\subsubsection{}

A Cartan datum 
is said to be \emph{irreducible} 
if the corresponding Cartan 
matrix cannot be made block-diagonal 
by simultaneous 
permutations of rows and columns.
 
A Cartan datum 
is said to be of \emph{finite type} if 
the corresponding Cartan matrix is 
positive definite.

A Cartan datum 
is said to be of \emph{affine type} if the
corresponding Cartan matrix is irreducible,
positive semi-definite, but not positive 
definite.

A Cartan datum 
is said to be \emph{symmetric} if 
$a_{ij}=a_{ji}$
for all $i$, $j$.

A Cartan datum 
is said to be \emph{simply laced} if 
$a_{ij}\in \{ 0,-1 \}$
for all $i \neq j$.

If a Cartan datum is simply laced then it is 
symmetric.
 
\subsubsection{}

{F}rom now on all the Cartan data is assumed to be 
symmetric.
See, however, \ref{IntroConclusion}.

\subsubsection{}

For a given finite set $S$ we denote by 
$\mathcal{P}_2(S)$
the set of all two-element subsets of $S$.

By definition, a finite graph is a triple 
$(I, E, \Ends)$, 
consisting of two finite sets $I$ (vertices) 
and $E$ (edges), 
and a map $\Ends: E \rightarrow 
\mathcal{P}_2 (I)$.

To a symmetric Cartan datum $(I,<,>)$ 
we associate 
a \emph{Dynkin graph} $(I, E, \Ends)$ as follows: 
the set of vertices
coincides with $I$, two vertices $i$ and $j$ 
are joined by
$-a_{ij}$ edges, there are no edges joining 
a vertex
with itself.

A Dynkin graph is called irreducible 
(resp. of finite type, of affine type) if 
the corresponding Cartan datum is irreducible 
(resp. of finite type, of affine type).

\subsubsection{}\label{Graphs}

The lists of all 
irreducible Dynkin graphs 
of finite and affine types 
are contained in Figures 1 and 2
respectively (cf. \cite{Bourbaki}).

\begin{figure}[ht]
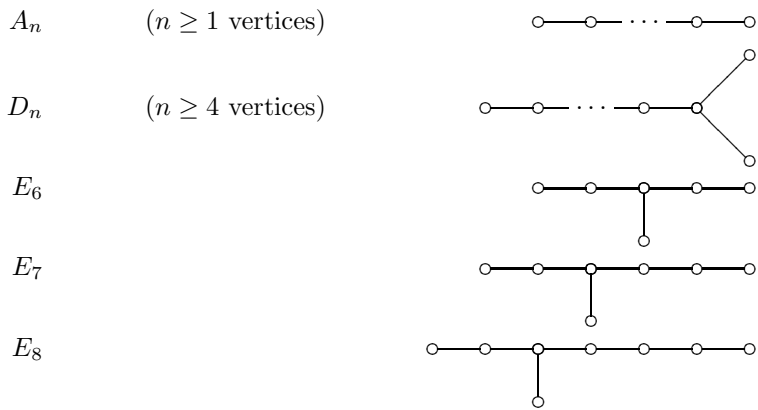

\begin{align}\nonumber
\text{$A_n$} && \text{($n \geq 1$ vertices)} &&
\xygraph{
[]
!{<0pt,0pt>;<20pt,0pt>:} 
*\cir<2pt>{}  - [r] *\cir<2pt>{} 
- [r] *{ \;  \dots \; }
- [r] *\cir<2pt>{} - [r] *\cir<2pt>{} 
}\\\nonumber
\text{$D_n$} && \text{($n \geq 4$ vertices)} &&
\xygraph{
[] 
!{<0pt,0pt>;<20pt,0pt>:} 
*\cir<2pt>{}  - [r] *\cir<2pt>{} 
- [r] *{ \;  \dots \; }
- [r] *\cir<2pt>{}
- [r] *\cir<2pt>{} - [ru] *\cir<2pt>{} 
[ld] *\cir<2pt>{} - [rd] *\cir<2pt>{}
}\\\nonumber
\text{$E_6$} && &&
\xygraph{
[] 
!{<0pt,0pt>;<20pt,0pt>:} 
*\cir<2pt>{} - [r] *\cir<2pt>{} 
- [r]  *\cir<2pt>{} - [d] *\cir<2pt>{} 
[u] *\cir<2pt>{}
- [r] *\cir<2pt>{} - [r] *\cir<2pt>{}
}\\\nonumber 
\text{$E_7$} && &&
\xygraph{
[] 
!{<0pt,0pt>;<20pt,0pt>:} 
*\cir<2pt>{} - [r] *\cir<2pt>{} 
- [r]  *\cir<2pt>{} 
- [d] *\cir<2pt>{} [u] *\cir<2pt>{}
- [r] *\cir<2pt>{}
- [r] *\cir<2pt>{} - [r] *\cir<2pt>{}
}\\\nonumber 
\text{$E_8$} && &&
\xygraph{
[] 
!{<0pt,0pt>;<20pt,0pt>:} 
*\cir<2pt>{} - [r] *\cir<2pt>{} 
- [r]  *\cir<2pt>{} 
- [d] *\cir<2pt>{} 
[u] *\cir<2pt>{}
- [r] *\cir<2pt>{} - [r] *\cir<2pt>{}
- [r] *\cir<2pt>{} - [r] *\cir<2pt>{}
} 
\end{align}
\caption{Finite Irreducible Dynkin Graphs}
\end{figure}

\begin{figure}[ht]
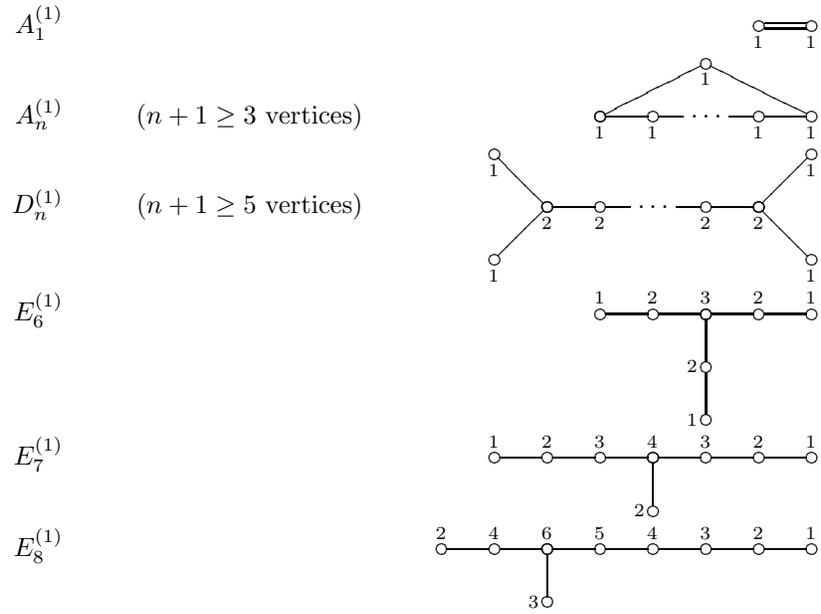

\begin{align}\nonumber
\text{$A_{1}^{(1)}$} &&  &&
\xygraph{
[]
!{<0pt,0pt>;<20pt,0pt>:} 
*\cir<2pt>{} 
!{\save -<0pt,6pt>*\txt{$_1$}  \restore}
-@2 [r] 
*\cir<2pt>{}
!{\save -<0pt,6pt>*\txt{$_1$}  \restore}
[d]
}\\\nonumber
\text{$A^{(1)}_n$} && \text{($n+1 \geq 3$ 
vertices)} &&
\xygraph{
[]
!{<0pt,0pt>;<20pt,0pt>:} 
*\cir<2pt>{}  
!{\save -<0pt,6pt>*\txt{$_1$}  \restore}
- [r] 
*\cir<2pt>{}
!{\save -<0pt,6pt>*\txt{$_1$}  \restore}
- [r] 
*{ \;  \dots \; }
- [r] 
*\cir<2pt>{} 
!{\save -<0pt,6pt>*\txt{$_1$}  \restore}
- [r] 
*\cir<2pt>{}
!{\save -<0pt,6pt>*\txt{$_1$}  \restore}
- [llu] 
*\cir<2pt>{}
!{\save -<0pt,6pt>*\txt{$_1$}  \restore}
- [lld] 
*\cir<2pt>{}
}\\\nonumber
\text{$D^{(1)}_n$} && \text{($n+1 \geq 5$ 
vertices)} &&
\xygraph{
[] 
!{<0pt,0pt>;<20pt,0pt>:} 
[u] 
*\cir<2pt>{}  
!{\save -<0pt,6pt>*\txt{$_1$}  \restore}
- [rd] 
*\cir<2pt>{}
!{\save -<0pt,6pt>*\txt{$_2$}  \restore}
- [ld] 
*\cir<2pt>{} 
!{\save -<0pt,6pt>*\txt{$_1$}  \restore}
[ru] 
*\cir<2pt>{}
- [r] 
*\cir<2pt>{}
!{\save -<0pt,6pt>*\txt{$_2$}  \restore}
- [r] 
*{ \;  \dots \; }
- [r] 
*\cir<2pt>{}
!{\save -<0pt,6pt>*\txt{$_2$}  \restore}
- [r] 
*\cir<2pt>{} 
!{\save -<0pt,6pt>*\txt{$_2$}  \restore}
- [ru] 
*\cir<2pt>{}
!{\save -<0pt,6pt>*\txt{$_1$}  \restore}
[ld] 
*\cir<2pt>{} 
- [rd] 
*\cir<2pt>{}
!{\save -<0pt,6pt>*\txt{$_1$}  \restore}
}\\\nonumber
\text{$E^{(1)}_6$} && &&
\xygraph{
[] 
!{<0pt,0pt>;<20pt,0pt>:} 
*\cir<2pt>{} 
!{\save -<0pt,-6pt>*\txt{$_1$}  \restore}
- [r] 
*\cir<2pt>{} 
!{\save -<0pt,-6pt>*\txt{$_2$}  \restore}
- [r]  
*\cir<2pt>{} 
!{\save -<0pt,-6pt>*\txt{$_3$}  \restore}
- [d] 
*\cir<2pt>{} 
!{\save -<5pt,0pt>*\txt{$_2$}  \restore}
- [d] 
*\cir<2pt>{}
!{\save -<5pt,0pt>*\txt{$_1$}  \restore}
[uu] 
*\cir<2pt>{}
- [r] 
*\cir<2pt>{} 
!{\save -<0pt,-6pt>*\txt{$_2$}  \restore}
- [r] 
*\cir<2pt>{}
!{\save -<0pt,-6pt>*\txt{$_1$}  \restore}
}\\\nonumber 
\text{$E^{(1)}_7$} && &&
\xygraph{
[] 
!{<0pt,0pt>;<20pt,0pt>:} 
*\cir<2pt>{} 
!{\save -<0pt,-6pt>*\txt{$_1$}  \restore}
- [r] *\cir<2pt>{} 
!{\save -<0pt,-6pt>*\txt{$_2$}  \restore}
- [r]  *\cir<2pt>{} 
!{\save -<0pt,-6pt>*\txt{$_3$}  \restore}
- [r]  *\cir<2pt>{}
!{\save -<0pt,-6pt>*\txt{$_4$}  \restore}
- [d] *\cir<2pt>{} 
!{\save -<5pt,0pt>*\txt{$_2$}  \restore}
[u] *\cir<2pt>{}
- [r] *\cir<2pt>{}
!{\save -<0pt,-6pt>*\txt{$_3$}  \restore}
- [r] *\cir<2pt>{} 
!{\save -<0pt,-6pt>*\txt{$_2$}  \restore}
- [r] *\cir<2pt>{}
!{\save -<0pt,-6pt>*\txt{$_1$}  \restore}
}\\\nonumber 
\text{$E^{(1)}_8$} && &&
\xygraph{
[] 
!{<0pt,0pt>;<20pt,0pt>:} 
*\cir<2pt>{} 
!{\save -<0pt,-6pt>*\txt{$_2$}  \restore}
- [r] *\cir<2pt>{} 
!{\save -<0pt,-6pt>*\txt{$_4$}  \restore}
- [r]  *\cir<2pt>{} 
!{\save -<0pt,-6pt>*\txt{$_6$}  \restore}
- [d] *\cir<2pt>{} 
!{\save -<5pt,0pt>*\txt{$_3$}  \restore}
[u] *\cir<2pt>{}
- [r] *\cir<2pt>{} 
!{\save -<0pt,-6pt>*\txt{$_5$}  \restore}
- [r] *\cir<2pt>{}
!{\save -<0pt,-6pt>*\txt{$_4$}  \restore}
- [r] *\cir<2pt>{} 
!{\save -<0pt,-6pt>*\txt{$_3$}  \restore}
- [r] *\cir<2pt>{}
!{\save -<0pt,-6pt>*\txt{$_2$}  \restore}
- [r] *\cir<2pt>{}
!{\save -<0pt,-6pt>*\txt{$_1$}  \restore}
}
\end{align} 
\caption{Affine Dynkin Graphs}
\end{figure}

The integers on vertices
of affine Dynkin graphs on Figure 2
are components of the first imaginary root
$\delta$ (see \ref{AffineRoots}).

Note that all symmetric Cartan data of
finite or affine type are simply laced 
except for $A^{(1)}_1$.

\subsubsection{}\label{DynkinFiniteAffine}

If one removes a vertex together with
adjacent edges from a Dynkin graph of affine
type one obtains a Dynkin graph of finite type (not
necessarily irreducible).

\subsection{A subalgebra 
of a Kac-Moody algebra}

\subsubsection{}\label{DefinitionOfKMPreliminary}

Let $\mathfrak{n}$ be the quotient of the 
free Lie algebra over $\mathbb{C}$ 
with generators 
$\{ e_i \}_{i \in I}$ by the ideal generated by the
following relations

\begin{equation}\label{KMSerre}
(\ad ( e_i ))^{1-a_{ij}} e_j = 0
\text{ if $i \neq j$ }  ,
\end{equation}
where $\ad (x) y = [x,y]$.

Relations \ref{KMSerre} are called the Serre 
relations.

The Lie algebra $\mathfrak{n}$ is 
$\mathbb{Z}_+ [I]$-graded with $\Deg e_i =i$.

\subsection{Quiver and the Euler cocycle}

\subsubsection{} 
A \emph{quiver} $Q$ is an oriented graph, that is
a quadruple $(I, \Omega, \In, \Out)$, consisting of  
two finite sets $I$ (vertices) and $\Omega$ 
(oriented edges), and
two maps $\In$ and $\Out$ from $\Omega$ to $I$. 
The underlying
non-oriented graph is given by $(I, \Omega, 
\{\In, \Out\})$.

A quiver is said to be of finite (resp. affine) 
type 
if the underlying non-oriented
graph is a Dynkin graph of finite (resp. affine) 
type.

\subsubsection{}\label{DefinitionOfKM} 

For a quiver $Q$ of finite or affine type 
we denote by $\mathfrak{n} (Q)$ the 
Lie algebra over $\mathbb{C}$ associated with
the underlying Dynkin graph of $Q$
as in \ref{DefinitionOfKMPreliminary}.

\subsubsection{}

The abelian group $\mathbb{Z}[I]$ is called
the \emph{root lattice}.

\subsubsection{}\label{DefinitionOfEpsilon}

Given a quiver $Q$ we 
introduce a biadditive form 
$e : \mathbb{Z} [I] \times \mathbb{Z} [I]
\rightarrow \mathbb{Z}$ given by the following 
values on generators:
\begin{equation}\nonumber
e(i,j) = \delta_{ij}
- \text{(number of $h\in \Omega$ 
such that  $i=\Out (h)$ and $j=\In (h)$),} 
\end{equation}
where $\delta_{ij}$ is the Kronecker symbol.
The form $e$ is called Euler form (see \ref{HomExt} 
for a justification of this name).

If $Q$ is of finite or affine type then 
\begin{equation}\label{EpsilonPlusEpsilon}
e ( \alpha , \beta ) + e ( \beta , \alpha ) =
< \alpha , \beta >.
\end{equation}

Having the Euler form $e$ we
define \emph{Euler cocycle} $\epsilon$ to be
a bimultiplicative function from 
$\mathbb{Z}[I] \times \mathbb{Z}[I]$
to $\{ \pm 1 \}$ given by 
\begin{equation}\nonumber
\epsilon (\alpha , \beta ) = 
(-1)^{e(\alpha , \beta )}.
\end{equation}

If we want to specify the quiver $Q$ used in
the definition of $\epsilon$ we write
$\epsilon_{Q}$.

Euler cocycle $\epsilon$ is the main building block of
the structure constants of $\mathfrak{n}$ in the
natural basis associated to the quiver $Q$. The precise
construction of the natural basis will be carried out
in Sections \ref{FiniteSection} and \ref{AffineSection}
in the finite and the affine cases respectively. 

\subsection{Category of representations of a quiver}
\label{DefinitionOfQmod}

\subsubsection{} 

Let $Q=(I, \Omega, \In, \Out)$ be a quiver.

A (finite-dimensional) \emph{representation} of  
$Q$ 
(over $\mathbb{C}$)
is the following data:

\begin{itemize}
\item a finite dimensional $I$-graded vector space 
$V = \oplus_{i \in I} V_i$
over $\mathbb{C}$,
\item a collection of linear maps 
$x = \oplus_{h \in \Omega} x_h 
\in \oplus_{h \in \Omega} 
\Hom_{\mathbb{C}} (V_{\Out(h)} , V_{\In (h)})$.
\end{itemize}

A morphism from a representation 
$( V , x) $
to another representation
$( V' ,x') $
is an $I$-graded $\mathbb{C}$-linear map 
$\phi : V \rightarrow V'$, 
such that 
$x'_h \phi_{\Out (h)} = \phi_{\In (h)} x_h$ 
for each $h \in \Omega$.
The composition of morphisms is the
composition of linear maps.

We denote the category of representations 
of the quiver by 
$\overline{\mathcal{M}} (Q)$. 
It is an abelian category with respect to the
natural additive structure on morphisms.

\subsubsection{}

Let $\mathcal{M} (Q)$ be a full subcategory of
$\overline{\mathcal{M}} (Q)$ with objects being 
$\mathbf{M}=(V,x)\in\Ob (\overline{\mathcal{M}} (Q))$
for which there exists $N(\mathbf{M}) \in \mathbb{N}$
such that $ x_{h_N} \dots x_{h_1} = 0$
for any sequence 
$h_1, \dots ,h_N \in \Omega$ with $\In (h_i) = 
\Out (h_{i+1})$.
We call objects of $\mathcal{M} (Q)$ \emph{nilpotent}
representations of $Q$. 

The subcategory $\mathcal{M} (Q)$ 
is closed with respect to extensions.
If $Q$ does not contain an oriented cycle then
$\mathcal{M} (Q) = 
\overline{\mathcal{M}} (Q)$.

\subsubsection{}

\emph{Dimension} of a representation 
$\mathbf{M} = ( V , x )$
of $Q$ is an element of $\mathbb{Z} [I]$ given by
the graded dimension of $V$:
\begin{equation}\nonumber
\dim_{\mathcal{M} (Q)} \mathbf{M} =  
\sum_{i\in I} \dim_{\mathbb{C}} (V_i) i .
\end{equation}

\subsubsection{}\label{DefinitionOfSimple}

It is easy to see that isomorphism
classes of simple objects in
$\mathcal{M} (Q)$ 
are in one-to-one correspondence
with vertices of $Q$ (we recall that
$\Ob (\mathcal{M} (Q))$ 
consists of nilpotent representations). 
The simple object
corresponding to a vertex $i \in I$ is 
given (up to an isomorphism) by 

\begin{equation}\nonumber
\mathbf{P}_{i} =
(V, x),
\text{ where } 
\begin{cases}
V_j = \{ 0 \}, j \neq i, \\
V_i = \mathbb{C} , \\
x = 0 .
\end{cases} 
\end{equation}

Note that 
$\dim_{\mathcal{M} (Q)} 
\mathbf{P}_i = i \in I 
\subset \mathbb{Z} [I]$. 

\subsubsection{}\label{DimSeries}

Let $\mathbf{M} \in \Ob (\mathcal{M} (Q))$.
It follows from \ref{DefinitionOfSimple} that
$(\dim_{\mathcal{M} (Q)} \mathbf{M})_i$ is equal to
the number of factors isomorphic to $\mathbf{P}_i$
in Jordan-H\"older series of $\mathbf{M}$.

\subsubsection{}\label{HomExt}

The following proposition explains 
the homological meaning
of the Euler form $e$. 

\begin{proposition}
\begin{enumerate}
\item\label{Hereditary}
The category $\mathcal{M} (Q)$ is hereditary,
i.e. $\Ext^n_{\mathcal{M} (Q)} 
(\mathbf{M}, \mathbf{N}) = 0$
for any $\mathbf{M},\mathbf{N} \in 
\Ob (\mathcal{M} (Q))$ and any $n \geq 2$.
\item\label{eExtExt}
The following equality holds
for any $\mathbf{M},\mathbf{N} \in 
\Ob (\mathcal{M} (Q))$.
\begin{multline}\nonumber
e(\dim_{\mathcal{M} (Q)} \mathbf{M},
\dim_{\mathcal{M} (Q)} \mathbf{N})= \\
= \dim_{\mathbb{C}} 
\Hom_{\mathcal{M} (Q)} (\mathbf{M},\mathbf{N}) - 
\dim_{\mathbb{C}} 
\Ext^1_{\mathcal{M} (Q)} (\mathbf{M},\mathbf{N}) 
= \\ =
\sum_{k=0}^{\infty} (-1)^k 
\dim_{\mathbb{C}} 
\Ext^k_{\mathcal{M} (Q)} (\mathbf{M},\mathbf{N})
\end{multline}
\end{enumerate}
\end{proposition}

\begin{proof}
The proof is standard. Let us consider the following
exact sequence of functors from 
$(\mathcal{M} (Q))^{\circ} \times \mathcal{M} (Q)$ 
to the category of 
$\mathbb{C}$-linear spaces.

\setcounter{equation}{2}

\begin{multline}\label{HomHomHomExt}
0 \xrightarrow{} 
\Hom_{\mathcal{M}(Q)} (\mathbf{M},\mathbf{N})
\xrightarrow{\pi}
\oplus_{i \in I} 
\Hom_{\mathbb{C}} ( V_i , W_i )
\xrightarrow{\rho}
\\ 
\xrightarrow{\rho}
\oplus_{h \in \Omega}
\Hom_{\mathbb{C}} (V_{\Out (h)}, W_{\In (h)})
\xrightarrow{\sigma}
\Ext^1_{\mathcal{M} (Q)} (\mathbf{M} , \mathbf{N})
\xrightarrow{} 0 .
\end{multline}
Here $\mathbf{M}=(V,x)$, $\mathbf{N}=(W,y) \in 
\Ob (\mathcal{M} (Q))$, and $\pi$, $\rho$, $\sigma$
are natural transformations given as follows
\begin{equation}\nonumber
\begin{split}
\pi ( \phi ) &= \oplus \phi_i ,
\\
\rho (\oplus \phi_i ) &=
\oplus (\phi_{\In (h)}  x_h - y_h \phi_{\Out (h)}) ,
\\
\sigma (\oplus \psi_h ) &= 
[ 0 \rightarrow \mathbf{N} \rightarrow 
\mathbf{E} \rightarrow \mathbf{M}
\rightarrow 0 ],
\end{split}
\end{equation}
where $\mathbf{E}=(\oplus_{i \in I}
(V_i \oplus W_i), \oplus_{h \in \Omega}
\bigl(\begin{smallmatrix}
x_h & 0 \\ \psi_h & y_h  
\end{smallmatrix}\bigr) )$.
The two middle terms in the exact sequence 
\ref{HomHomHomExt} are exact functors
(both with respect to the first and
to the second arguments). Thus 
$\Ext^1_{\mathcal{M} (Q)}$ is right exact,
which proves \ref{Hereditary}. The equality
\ref{eExtExt} follows from evaluation of dimensions
in \ref{HomHomHomExt}.
\end{proof}

\subsubsection{}

Let $\mathcal{X} (Q)$ be the set of isomorphism 
classes
of objects of $\mathcal{M} (Q)$.
We denote the isomorphism class of 
$\mathbf{M} \in \Ob (\mathcal{M}(Q))$ by 
$[\mathbf{M}]$. 

We also use the following notation:

\begin{align}\nonumber
\mathcal{X}_{\alpha} (Q)=& 
\{ [\mathbf{M}] \in \mathcal{X} (Q)  | 
\dim_{\mathcal{M} (Q)} \mathbf{M} = \alpha \} ,
\\\nonumber
\mathcal{T} (Q) = &
\{ [\mathbf{M}] \in \mathcal{X} (Q) | 
\text{ $\mathbf{M}$ is indecomposable }  \} , 
\\\nonumber
\mathcal{T}_{\alpha} (Q) = &
\{ [\mathbf{M}] \in \mathcal{T} (Q)  | 
\dim_{\mathcal{M} (Q)} \mathbf{M} = \alpha \} .
\end{align}

So $\mathcal{T} (Q)$ is the set of isomorphism 
classes
of indecomposable objects.
For example, $\mathcal{X}_{i} (Q) =
\mathcal{T}_i (Q) = \{ [\mathbf{P}_i] \}$.

\subsubsection{}\label{DefinitionOfE}

The set $\mathcal{X}_{\alpha} (Q)$ has
a natural structure of the 
orbit space for an action of
an algebraic group on an 
affine algebraic variety over
$\mathbb{C}$. 

Let $\alpha = \sum_{i \in I} 
\alpha_i i \in \mathbb{Z}_+ [I]$.
Consider an affine space over $\mathbb{C}$ 
\begin{equation}\nonumber
\mathbf{E}_{\alpha} = \bigoplus_{h \in \Omega}
\Hom_{\mathbb{C}} 
( \mathbb{C}^{\alpha_{\Out (h)}} , 
\mathbb{C}^{\alpha_{\In (h)}} ).  
\end{equation}

Let $\mathbf{G}_{\alpha} = \prod_{i \in I}
GL ( \alpha_i , \mathbb{C} )$. The group
$\mathbf{G}_{\alpha}$ acts on 
$\mathbf{E}_{\alpha}$
by $x^g_h = g_{\In (h)} x_h g_{\Out (h)}^{-1}$,
$g \in \mathbf{G}_{\alpha}$.

When we want to specify the quiver $Q$ used
in the definition of $\mathbf{E}_{\alpha}$
and $\mathbf{G}_{\alpha}$ we write
$\mathbf{E}_{\alpha} (Q)$
and $\mathbf{G}_{\alpha} (Q)$.

There is a natural bijection between the
set of isomorphism classes of objects in
$\overline{\mathcal{M}} (Q)$ 
of dimension $\alpha$ and the set
of orbits of $\mathbf{G}_{\alpha}$ in 
$\mathbf{E}_{\alpha}$. 
Namely, given a point
$x \in \mathbf{E}_{\alpha}$, 
$( \mathbb{C}^{\alpha} ,
x)$ is a representative
of the isomorphism class of 
objects in $\overline{\mathcal{M}} (Q)$, 
corresponding to this point. This class depends
only on the orbit to which the point belongs.
And vice versa, in each isomorphism class
there are objects of the form 
$(\mathbb{C}^{\alpha} , x)$. 

\subsubsection{}\label{change1--mv}
Let $c=(h_1, \dots , h_n)$ be an oriented cycle, 
that is
$n$-tuple of elements of $\Omega$, such that
$\In (h_i) = \Out (h_{i+1})$ for $1 \leq i 
\leq n-1$ and
$\In (h_n) = \Out (h_1)$. Given $\mathbf{M} = (V,x)
\in \Ob (\overline{\mathcal{M}} (Q))$, we call 
$x_{h_n} \dots x_{h_1} 
\in \Hom_{\mathbb{C}} 
(V_{\Out (h_1)},V_{\Out (h_1)})$
the holonomy of $x$ around the cycle $c$.

An object $\mathbf{M}=(V,x)\in \Ob (\overline
{\mathcal{M}} (Q))$ is nilpotent if and
only if the trace of the holonomy of $x$ 
around every cycle $c$ is equal to zero. 
Thus the set $\mathbf{E}^{nil}_{\alpha} = 
\{ x\in \mathbf{E}_{\alpha} | 
(\mathbb{C}^{\alpha}, x) 
\in \mathcal{M}(Q) \}$ is an affine 
subvariety in
$\mathbf{E}_{\alpha}$.
One can consider functions on 
$\mathcal{X}_{\alpha} (Q)$ as
$\mathbf{G}_{\alpha}$-invariant functions
on $\mathbf{E}^{nil}_{\alpha}$.

\section{Convolution algebra}\label{Convolution}

\subsection{Algebraically constructible functions}

\subsubsection{}

Let $X$ be an algebraic variety over $\mathbb{C}$. 

A \emph{constructible} set in $X$ is a set
obtained from subvarieties in $X$ by finitely
many standard set theoretic operations.

A function on $X$ is called \emph{constructible}
if $X$ has a finite partition into constructible
sets such that the function is constant on each
of them.

We denote by $M (X)$ the set of all constructible
functions on $X$ with values in $\mathbb{C}$.
The set $M (X)$ is naturally a 
$\mathbb{C}$-linear space.

Let $G$ be an algebraic group acting on $X$.
We denote by $M_{G} (X)$ a subspace of $M (X)$,
consisting of all $G$-invariant functions.

\subsubsection{}

Let $f: X \rightarrow Y$ be a morphism of
algebraic varieties.
Then $f^\ast$ denotes a $\mathbb{C}$-linear 
map from $M (Y)$ to $M (X)$ 
defined as follows:

\begin{equation}\nonumber
(f^\ast (\phi)) (x) =  \phi( f (x)) .
\end{equation}

Let $f: X \rightarrow Y$ be a proper morphism of
algebraic varieties.
Then $f_\ast$ denotes a $\mathbb{C}$-linear 
map from $M (X)$ to $M (Y)$ 
defined as follows \cite{MacPherson1974}:

\begin{equation}\nonumber
(f_\ast (\phi)) (y) = 
\sum_{a\in\mathbb{C}}
a \chi(f^{-1}(y) \cap \phi^{-1} (a)),
\end{equation}
where $\chi$ denotes the Euler characteristic 
with compact support. 

Note that if $f$ is an equivariant
morphism of $G$-varieties, then
the restriction of $f^\ast$ (resp. $f_\ast$) 
gives a $\mathbb{C}$-linear map from 
$M_G (Y)$ to $M_G (X)$ 
(resp. from $M_G (X)$ to $M_G (Y)$).

Let $G$ and $H$ be algebraic groups,
$X$ and $Y$ be algebraic varieties
with $H$ actions,
and $f: X\rightarrow Y$ be a locally
trivial $H$-equivariant principal $G$-bundle. 
Then $f_\flat$ is a $\mathbb{C}$-linear 
map from $M_{G\times H} (X)$ to $M_H (Y)$ 
defined as follows:

\begin{equation}\nonumber
(f_\flat (\phi)) (y) =  \phi (f^{-1} (y)).
\end{equation}

The map $f_\flat (\phi)$ is well defined because
$\phi \in M_G (X)$.

\begin{proposition}\label{AstFlatFunctorial}
The maps $f^\ast$, $f_\ast$, and $f_\flat$ have the 
following
properties:
\begin{enumerate}
\item $f^\ast$, $f_\ast$, and $f_\flat$ 
are functorial:
\begin{gather}\nonumber
(f\circ g)^\ast = g^\ast \circ f^\ast,
\\\nonumber
(f\circ g)_\ast = f_\ast \circ g_\ast,
\\\nonumber
(f\circ g)_\flat = f_\flat \circ g_\flat,
\end{gather}
when the right-hand side if defined,
\item
if $f: X\rightarrow Y$ 
is a locally trivial 
$H$-equivariant principal $G$-bundle
then $f_\flat$ and $f^*$ are inverse to each
other and give an isomorphism
between $M_{G \times H} (X)$ and $M_{H} (Y)$.
\end{enumerate}
\end{proposition}
\begin{proof}
All the statements except the functoriality of 
$f_{\ast}$ are obvious. The functoriality of
$f_{\ast}$ follows from properties of
the Euler characteristic (namely, additivity
with respect to algebraic stratifications,
and multiplicativity for fiber bundles). 
\end{proof}

\subsection{Convolution algebra and Lie subalgebras}

\subsubsection{}\label{DefinitionOfCCal}

We introduce the following 
$\mathbb{C}$-linear spaces:
\begin{gather}\nonumber
\mathcal{L}_{\alpha} (Q) = 
M_{\mathbf{G_{\alpha}}} (\mathbf{E}_{\alpha}^{nil}), 
\\\nonumber
\mathcal{L} (Q) = 
\bigoplus_{\alpha \in \mathbb{Z}_+ [I]}
\mathcal{L}_{\alpha} (Q)
\end{gather}

Due to \ref{DefinitionOfE} we can consider
elements of $\mathcal{L} (Q)$ as functions
on $\mathcal{X} (Q)$ with values in 
$\mathbb{C}$. Note however, that not every
function on $\mathcal{X} (Q)$ belongs to
$\mathcal{L} (Q)$.

In particular, given a representation $\mathbf{M}$
of $Q$ we sometimes say "the characteristic function
of $[M]$" instead of "the characteristic function
of $\mathbf{G}_{\alpha}$-orbit 
in $\mathbf{E}_{\alpha}$, corresponding to
$[M]$". In either case, one should make sure that
the function in question is constructible. It is
always true, for example, if the number
of $\mathbf{G}_{\alpha}$-orbits is finite,
or if the function is in the image of 
a Hall map (see \ref{HallMaps}).

\subsubsection{}\label{DefinitionOfEi}

For example
\begin{equation}\nonumber
\mathcal{L}_{i} (Q) =  
\mathbb{C} E_i
\text{ for } i\in I \subset \mathbb{Z} [I] ,
\end{equation}
where $E_i$ is the function equal to
$1$ on $\mathbf{E}_{i}=\mathbf{E}_i^{nil}$
(which is a point).

\subsubsection{}\label{EDiagrams}

We endow $\mathcal{L}(Q)$ with a bilinear
product $\ast$, graded by $\mathbb{Z}_+ [I]$,

\begin{equation}\label{GradingForAst}
\mathcal{L}_{\alpha} (Q) \ast 
\mathcal{L}_{\beta} (Q)
\subset \mathcal{L}_{\alpha + \beta} (Q).
\end{equation}

By linearity it is enough to define $\ast$-product for
a pair of functions 
$f_{\alpha} \in \mathcal{L}_{\alpha} (Q)$
and 
$f_{\beta} \in \mathcal{L}_{\beta} (Q)$. Moreover
$f_{\alpha} \ast f_{\beta}$ should belong to
$\mathcal{L}_{\alpha + \beta} (Q)$  for the grading 
property \eqref{GradingForAst} to hold.
 
The product $f_{\alpha} \ast f_{\beta}$ 
is defined as follows \cite{Lusztig1991a, Lusztig1991b}.
Consider a diagram of varieties:
\begin{equation}
\mathbf{E}_{\alpha}^{nil} \times 
\mathbf{E}_{\beta}^{nil}
\xleftarrow{p_1}
\mathbf{E'} 
\xrightarrow{p_2}
\mathbf{E''}
\xrightarrow{p_3}
\mathbf{E}_{\alpha + \beta}^{nil},
\end{equation}
where the notation is as follows:

$\mathbf{E}_{\alpha}^{nil}$, 
$\mathbf{E}_{\beta}^{nil}$,
$\mathbf{E}_{\alpha + \beta}^{nil}$ 
are defined in \ref{DefinitionOfE},

$\mathbf{E''}$ is the variety of all pairs $(x, W)$, 
consisting
of $x \in \mathbf{E}^{nil}_{\alpha + \beta}$ 
and an $x$-stable
$I$-graded subspace of $\mathbb{C}^{\alpha + \beta}$
such that $\dim W = \alpha$,

$\mathbf{E'}$ is the variety of all quadruples 
$(x, W, R', R'')$, where $(x, W) \in \mathbf{E''}$,
$R'$ is an isomorphism 
$\mathbb{C}^{\alpha} \Tilde{\rightarrow} W$,
$R''$ is an isomorphism 
$\mathbb{C}^{\beta} \Tilde{\rightarrow} 
\mathbb{C}^{\alpha+\beta} / W$.

$p_2(x,W,R',R'')=(x,W)$, 

$p_3(x,W)=x$,

$p_1(x,W,R',R'')=(x',x'')$, where 
$x_h R'_{\Out (h)}=R'_{\In (h)} x_h'$,
and $x_h R''_{\Out (h)}=R''_{\In (h)} x_h''$ 
for all $h\in H$. 

Note that 
$p_2$ is a principal 
$\mathbf{G}_{\alpha} \times \mathbf{G}_{\beta}$ 
fibration, and $p_3$ is proper.  

Given $f_{\alpha} \in \mathcal{L}_{\alpha} (Q)$ and
$f_{\beta} \in \mathcal{L}_{\beta} (Q)$, let 
$g( x_1, x_2 )= f_{\alpha} (x_1) f_{\beta} (x_2)$ be
an algebraically constructible function on
$\mathbf{E}^{nil}_{\alpha} \times 
\mathbf{E}^{nil}_{\beta}$.
By definition 
\begin{equation}\nonumber
f_{\alpha}*f_{\beta}=
(p_3)_\ast (p_2)_\flat (p_1)^\ast (g).
\end{equation}

\begin{theorem}The space 
$\mathcal{L} (Q)$ equipped with the $*$-product
is a $\mathbb{Z}_+ [I]$-graded associative
algebra over $\mathbb{C}$. 
\end{theorem}

\begin{proof}
Associativity follows from functorial properties
\ref{AstFlatFunctorial} 
of the maps $p^{\ast}, p_{\ast}, p_{\flat}$.
\end{proof}

\subsubsection{} The algebra $\mathcal{L} (Q)$ 
with the $\ast$-product
defined above 
is called the \emph{Hall algebra}. 
In our definition the $\ast$-product is the 
opposite of 
the one given in \cite{Lusztig1991b} (which
goes back to Ringel \cite{Ringel1988}).

\subsubsection{}

Let $Q$ be a quiver of finite or affine type.
Consider $\mathcal{L} (Q)$ as
a $\mathbb{Z}_+ [I]$-graded
Lie algebra over $\mathbb{C}$, 
with the following Lie bracket:
\begin{equation}\label{ConvolutionBracket}
[f,g]=f \ast g - g \ast f .
\end{equation}

We denote by $\mathfrak{n}^{\ast} (Q)$ the Lie 
subalgebra
of $\mathcal{L} (Q)$ generated by 
$\{ E_i \}_{i \in I}$ (see \ref{DefinitionOfEi}). 

\begin{numproposition}
$\mathfrak{n}^{\ast} (Q)$ has the 
following properties:
\begin{enumerate}
\item
$\mathfrak{n}^{\ast} (Q)$ is a 
$\mathbb{Z}_+ [I]$-graded Lie algebra: 
$\mathfrak{n}^{\ast} (Q) = 
\bigoplus_{\alpha \in \mathbb{Z}_+ [I]}
\mathfrak{n}^{\ast}_\alpha (Q)$, 
\item
$\dim_\mathbb{C} \mathfrak{n}^{\ast}_\alpha (Q) < 
\infty$
for any $\alpha \in \mathbb{Z}_+ [I]$,
\item
$\mathfrak{n}^{\ast}_{i} (Q) = \mathbb{C} E_i$  for 
$i \in I \subset \mathbb{Z}_+ [I]$,
where $E_i$ is defined in \ref{DefinitionOfEi},
\item\label{LRSSerre}
the generators $E_i$ satisfy the following
relations:
\begin{equation}\nonumber
(\ad ( E_i ))^{1-a_{ij}} E_j = 0
\text{ for $i \neq j$ },
\end{equation}
where $\ad (X) Y = [X,Y]$, and $a_{ij}=<i,j>$.
\end{enumerate}
\end{numproposition}
\begin{proof}
\ref{LRSSerre} is a simple calculation
using the definition
of the $\ast$-product,
all the rest is obvious. 
\end{proof}

\subsubsection{}\label{DefinitionOfNMMM}

Given $\mathbf{M}_1$, $\mathbf{M}_2$,
$\mathbf{M}_3 \in \Ob (\mathcal{M} (Q))$
we denote by 
$ n^Q_{\mathbf{M}_1, \mathbf{M}_2 ;
\mathbf{M}_3 }$ the Euler characteristic 
with compact support of the variety
$N^Q_{\mathbf{M}_1, \mathbf{M}_2;\mathbf{M_3}}$
of all
subobjects $\mathbf{V}$ of $\mathbf{M}_3$ such that
$[\mathbf{V}]=[\mathbf{M}_1]$ and
$[\mathbf{M}_3 / \mathbf{V}] =[\mathbf{M}_2]$
(we consider $N^Q_{\mathbf{M}_1, \mathbf{M}_2;
\mathbf{M_3}}$
as a constructible subvariety in a product of 
Grassmannians):
\begin{equation}\nonumber
n^Q_{\mathbf{M}_1, \mathbf{M}_2 ;
\mathbf{M}_3 } =
\chi ( N^Q_{\mathbf{M}_1, \mathbf{M}_2 ;
\mathbf{M}_3 }).
\end{equation}

We denote by $\mathcal{L}^{ind} (Q)$ the subspace
of $\mathcal{L} (Q)$ consisting of all
$f \in \mathcal{L} (Q)$ such that
$f ([\mathbf{M}]) = 0$ if
$\mathbf{M} \in \Ob (\mathcal{M} (Q))$
is decomposable.

\begin{numproposition}\label{Reidtmann}
\begin{enumerate}
\item\label{Reidtmann1}
Let $\mathbf{M'}$ and $\mathbf{M''}$ be two 
indecomposable 
objects of $\mathcal{M} (Q)$. 
If $n^Q_{\mathbf{M'} , \mathbf{M''} ; \mathbf{M}} 
\neq 0$
then either $\mathbf{M}$ is indecomposable or 
$\mathbf{M} = \mathbf{M'} \oplus \mathbf{M''}$.
\item\label{Reidtmann2}
Let $\mathbf{M'}$ and $\mathbf{M''}$ be two 
indecomposable 
objects of $\mathcal{M} (Q)$. Then
\begin{equation}\nonumber
n^Q_{\mathbf{M'} , \mathbf{M''} ; 
\mathbf{M'} \oplus \mathbf{M''}} = 
n^Q_{\mathbf{M''} , \mathbf{M'} ; 
\mathbf{M'} \oplus \mathbf{M''}} =
\begin{cases}
1 & \text{ if  
$\mathbf{M'}$ is not isomorphic to $\mathbf{M''}$ 
}, \\
2 & \text{ if 
$\mathbf{M'}$ is isomorphic to $\mathbf{M''}$ 
} .
\end{cases} 
\end{equation}
\item\label{Reidtmann3}
The subspace $\mathcal{L}^{ind} (Q) \subset 
\mathcal{L} (Q)$
is closed with respect to Lie bracket 
\ref{ConvolutionBracket}.
\item\label{ReidtmannFinal}
$\mathfrak{n}^{\ast} (Q) \subset 
\mathcal{L}^{ind} (Q)$,
that is if $f\in  \mathfrak{n}^{\ast} (Q)$ then
$f ([\mathbf{M}])=0$ for any decomposable $\mathbf{M}
\in \mathcal{M} (Q)$.
\end{enumerate}
\end{numproposition}
\begin{proof}
\ref{Reidtmann1} and \ref{Reidtmann2} 
have been proven by
Ch. Riedtmann \cite{Riedtmann1994}.

\ref{Reidtmann3} follows from \ref{Reidtmann1} 
and \ref{Reidtmann2}.

\ref{ReidtmannFinal} follows from \ref{Reidtmann3}
and from the fact that $E_i \in 
\mathcal{L}^{ind} (Q)$
for any $i\in I$.
\end{proof}

\subsubsection{}

Because of \ref{LRSSerre} we have a surjective
$\mathbb{Z} [I]$-graded 
homomorphism of Lie algebras $\Xi^{\ast}$ from 
the Lie algebra $\mathfrak{n} (Q)$ 
introduced in \ref{DefinitionOfKM}
to $\mathfrak{n}^{\ast} (Q)$, defined as follows:
\begin{gather}\nonumber
\Xi^{\ast}:
\mathfrak{n} (Q) \rightarrow 
\mathfrak{n}^{\ast}(Q) , 
\\\nonumber
\Xi^{\ast} (e_i) = E_i , \quad i \in I. 
\end{gather}

Our goal in this paper is to investigate 
the homomorphism $\Xi^{\ast}$.
For this study we utilize functorial
properties of the $\ast$-product.

\section{Functors}\label{Functors}

\subsection{Hall functors and Hall maps}
\label{HallMaps}

\subsubsection{}

Let $Q = (I, \Omega, \In ,\Out)$, 
$Q' = (I', \Omega', \In ' ,\Out ')$ be quivers,
$\mathcal{F} : 
\mathcal{M}(Q) \rightarrow \mathcal{M} (Q')$ be
a full, faithful, exact functor, such that
$\im \mathcal{F}$ is 
\emph{\'epaisse} (closed
with respect to extensions).

We denote by $\dim \mathcal{F}$ an additive 
map from $\mathbb{Z} [I]$ to $\mathbb{Z} [I']$
such that 

\begin{equation}\nonumber
\dim\mathcal{F} 
(\dim_{\mathcal{M} (Q)} \mathbf{M}) = 
\dim_{\mathcal{M} (Q')}
(\mathcal{F} (\mathbf{M}))
\end{equation}
for any $\mathbf{M}\in \mathcal{M} (Q)$. The map
$\dim \mathcal{F}$ is well-defined due to
\ref{DimSeries}.

\subsubsection{}\label{PhiPsi}

Let $Q = (I, \Omega, \In ,\Out)$, 
$Q' = (I', \Omega', \In ' ,\Out ')$ be quivers, 
$\mathcal{F} : 
\mathcal{M}(Q) \rightarrow \mathcal{M} (Q')$ be
a full, faithful, exact functor with 
\emph{\'epaisse} image.
We call $\mathcal{F}$ a
\emph{Hall functor} if there exist two 
sets of linear maps
\begin{gather}\nonumber
\{ \phi_i \in \oplus_{h' \in \Omega'} 
\Hom_{\mathbb{C}} (
\mathbb{C}^{(\dim\mathcal{F} (i))_{\Out (h')}}, 
\mathbb{C}^{(\dim\mathcal{F} (i))_{\In (h')}} 
) \}_{i \in I} ,
\\\nonumber
\{ \psi_h \in \oplus_{h' \in \Omega'} 
\Hom_{\mathbb{C}} (
\mathbb{C}^{(\dim\mathcal{F} (\Out (h)))_{\Out (h')}}, 
\mathbb{C}^{(\dim\mathcal{F} (\In (h)))_{\In (h')}} 
) \}_{h \in \Omega},
\end{gather} 
such that $\mathcal{F}$ is given by the following 
action on 
objects of $\mathcal{M} (Q)$
\begin{gather}\nonumber
\mathcal{F} ((V,x)) = (W,z),
\text{ where }
\\\label{blocks}
W = \oplus_{i \in I} ( V_i \otimes_{\mathbb{C}}
\mathbb{C}^{\dim\mathcal{F} (i)}) ,
\\\nonumber
z = (\oplus_{i \in I} (\Id_{V_i} \otimes
\phi_i))
\oplus
(\oplus_{h \in \Omega} (x_h \otimes
\psi_h )) ,  
\end{gather}
and the natural action on morphisms.

We borrow all functors we use from the theory 
of representations of quivers, however we want them 
to be given
by specific formulas \eqref{blocks}. The following 
proposition
asserts that this condition is not too restrictive.

\begin{proposition}
Let $Q$ and  $Q'$ be quivers, 
$\mathcal{H} : 
\mathcal{M}(Q) \rightarrow \mathcal{M} (Q')$ be
a full, faithful, exact functor with 
\emph{\'epaisse} image.
Then there exist a Hall functor $\mathcal{F}:
\mathcal{M} (Q) \rightarrow \mathcal{M} (Q')$
and an automorphism $\mathcal{G}$ of the category
$\mathcal{M} (Q)$ such that
$\mathcal{H}$ is naturally equivalent
to $\mathcal{F} \circ \mathcal{G}$.
\end{proposition}
\begin{proof}
Given $h \in \Omega$ let $\mathbf{P}_h = (U,t) \in 
\Ob (\mathcal{M} (Q))$, where
\begin{align}\nonumber
U_i &= \mathbb{C} 
\text{ if $i = \Out (h)$ 
or $i = \In (h)$} ,
\\\nonumber
U_i &= \{ 0 \} 
\text{ otherwise },
\\\nonumber
t_h &= \Id_{\mathbb{C}} ,
\\\nonumber
t_s &= 0 
\text{ if } s \neq h .
\end{align}
We fix an object of the form 
$(\mathbb{C}^{\dim\mathcal{F} (i)}, \phi_i )$
(resp. 
$(\mathbb{C}^{\dim\mathcal{F} (\Out (h))} \oplus
\mathbb{C}^{\dim\mathcal{F} (\In (h))} $, $
\bigl( \begin{smallmatrix}
\phi_{\Out (h)} & 0 \\
\psi_h & \phi_{\In (h)}
\end{smallmatrix}
\bigr) )$) in the isomorphism class
of $\mathcal{H} ( \mathbf{P}_i )$
for each $i \in I$
(resp. $\mathcal{H} (\mathbf{P}_h)$
for each $h \in \Omega$).
Using $\{ \phi_i \}_{i \in I}$ and
$\{ \psi_h \}_{h \in \Omega }$ we construct
a functor $\mathcal{F} :
\mathcal{M} (Q) \rightarrow \mathcal{M} (Q')$. 
Namely the action
of $\mathcal{F}$ on objects is given by 
\eqref{blocks}, and
the action on morphisms is the natural one.
It is easy to see that $\mathcal{F}$ is exact,
$\mathcal{F} (\mathbf{P}_i)$ is isomorphic to
$\mathcal{H} (\mathbf{P}_i)$ for any $i \in I$,
and the natural map 
$\Ext^1_{\mathcal{M} (Q)} (\mathbf{P}_i , 
\mathbf{P}_j ) 
\rightarrow
\Ext^1_{\mathcal{M} (Q')} ( \mathcal{F} ( 
\mathbf{P}_i ), 
\mathcal{F} ( \mathbf{P}_j ))$ is a bijection
for any $i$, $j \in I$. Now one
can deduce the statement of the proposition by
induction on length of
Jordan-H\"{o}lder series of an object 
of $\mathcal{M} (Q)$.
\end{proof}

\subsubsection{}\label{HallMap}

Let $\mathcal{F} :
\mathcal{M} (Q) \rightarrow 
\mathcal{M} (Q')$ be a Hall functor.
As $\mathcal{F}$ is full and faithful
it induces an injective map from the set
$\mathcal{X} (Q)$
of isomorphism classes of representations of $Q$
to the set 
$\mathcal{X} (Q')$ of isomorphism classes of 
representations
of $Q'$ (or from the set of 
$\mathbf{G}_{\alpha} (Q)$-orbits in
$\mathbf{E}^{nil}_{\alpha} (Q)$ to the set of
$ \mathbf{G}_{\dim\mathcal{F} ( 
\alpha )} (Q')$-orbits in
$\mathbf{E}^{nil}_{\dim\mathcal{F} ( \alpha )} (Q')$). 
Therefore
one can define a push-forward map 
$\mathfrak{f}$ from $\mathcal{L} (Q)$
to the set of $\mathbb{C}$-valued functions on  
$\mathcal{X} (Q')$ as follows

\begin{gather}\nonumber
(\mathfrak{f} (g)) 
([\mathcal{F} (\mathbf{M})]) = g ([ \mathbf{M} ]) ,
\\\nonumber
(\mathfrak{f} (g)) ([\mathbf{X}]) = 0
\text{ if $\nexists \mathbf{N} \in 
\Ob (\mathcal{M} (Q))$ 
such that $\mathbf{X}
= \mathcal{F} (\mathbf{N})$ } .
\end{gather}

We call $\mathfrak{f}$ the \emph{Hall map} 
associated with the
Hall functor $\mathcal{F}$.

\begin{proposition}
Let $\mathcal{F} : 
\mathcal{M} (Q) \rightarrow 
\mathcal{M} (Q')$ be a Hall
functor, and
$\mathfrak{f}$ be the corresponding
Hall map. Then
\begin{enumerate}
\item
$\im \mathfrak{f} \subset \mathcal{L} (Q')$, that is 
$\mathfrak{f}$ maps constructible 
$\mathbf{G}_{\alpha} (Q)$-invariant functions on 
the variety
$\mathbf{E}^{nil}_{\alpha} (Q)$ to
constructible 
$\mathbf{G}_{\dim\mathcal{F} ( 
\alpha )} (Q')$-invariant 
functions on the variety 
$\mathbf{E}^{nil}_{\dim\mathcal{F} ( \alpha )} (Q')$,
\item
$\mathfrak{f} : \mathcal{L} (Q) \rightarrow
\mathcal{L} (Q')$ is a homomorphism of algebras:
\begin{equation}\nonumber
\mathfrak{f} (g) \ast \mathfrak{f} (h) =
\mathfrak{f} ( g \ast h )
\text{ for any $g$, $h \in \mathcal{L} (Q)$}.
\end{equation}
\end{enumerate}
\end{proposition}
\begin{proof}
Follows from definitions. Note that being a Hall 
functor
$\mathcal{F}$ induces algebraic
maps between varieties 
$\mathbf{E}^{nil}_{\alpha}$, 
$\mathbf{E}'$, 
$\mathbf{E}''$
(used in \eqref{EDiagrams}) for the quiver $Q$ and
the corresponding varieties for the quiver $Q'$.
\end{proof}

\subsubsection{}\label{HallInd}

Let $\mathcal{F}:
\mathcal{M} (Q)
\rightarrow
\mathcal{M} (Q')$ be a Hall functor, in particular,
$\mathcal{F}$ is full and faithful. Then
$\mathcal{F} (\mathbf{M})$ is indecomposable
if and only if $\mathbf{M}$ is indecomposable.
It follows that $\mathcal{F}$ induces an injective
map from the set 
$\mathcal{T} (Q)$ of isomorphism classes of
indecomposable objects of 
$\mathcal{M} (Q)$ to the set 
$\mathcal{T} (Q')$ of isomorphism classes of 
indecomposable objects of
$\mathcal{M} (Q')$.

\subsubsection{}\label{HallAndEuler}

The following proposition describes
functorial properties of the Euler cocycle 
$\epsilon$
and the bilinear form $<,>$ with respect 
to a Hall functor.

\begin{proposition}
Let $\mathcal{F} : \mathcal{M} (Q) \rightarrow
\mathcal{M} (Q')$ be a Hall functor. Then
\begin{equation}\nonumber
\begin{split}
\epsilon_Q 
( \dim_{\mathcal{M} (Q)} \mathbf{M}_1, 
\dim_{\mathcal{M} (Q)} \mathbf{M}_2 ) &=
\epsilon_{Q'} 
( \dim_{\mathcal{M} (Q')} \mathcal{F}(\mathbf{M}_1), 
\dim_{\mathcal{M} (Q')} \mathcal{F}(\mathbf{M}_2) ) , 
\\ < \dim_{\mathcal{M} (Q)} \mathbf{M}_1, 
\dim_{\mathcal{M} (Q)} \mathbf{M}_2 >_{Q} &=
< \dim_{\mathcal{M} (Q')} \mathcal{F}(\mathbf{M}_1), 
\dim_{\mathcal{M} (Q')} 
\mathcal{F}(\mathbf{M}_2) >_{Q'}.
\end{split}
\end{equation}
for any $\mathbf{M}_1$, $\mathbf{M}_2
\in \Ob ( \mathcal{M} (Q))$.
\end{proposition}

\begin{proof}
The proposition follows from \ref{HomExt} and from
the fact that $\mathcal{F}$ is full, faithful,
exact, and has \emph{\'epaisse} image. 
\end{proof}

\subsubsection{}\label{EmbeddingFunctor}

The simplest example of a Hall functor is an embedding
functor. 

We call a subquiver $Q' = 
(I', \Omega', \In', \Out')$
of a quiver $Q = (I, \Omega, \In , \Out )$ 
\emph{full} if 
any $h \in \Omega$ such that $\In (h)\in I'$ and
$\Out (h) \in I'$ belongs to $\Omega'$.

Let $Q'$ be a full subquiver of $Q$.
Then the natural embedding 
$\mathcal{I}_{Q' \subset Q} : \mathcal{M} (Q')
\rightarrow \mathcal{M} (Q)$ clearly satisfies
all the axioms of a Hall functor, and therefore
induces an injective homomorphism 
$\mathfrak{i}_{Q'\subset Q} :
\mathcal{L} (Q') \rightarrow
\mathcal{L} (Q)$. Moreover,  
\begin{equation}\nonumber
\mathfrak{i}_{Q'\subset Q} (E_k) =
E_k ,
\end{equation}
where $k$ is a vertex of $Q'$. It follows that
$\mathfrak{i}_{Q'\subset Q}$
induces an injective homomorphism from
$\mathfrak{n}^{\ast} (Q')$ 
to $\mathfrak{n}^{\ast} (Q)$.

In the next section we consider a more intricate
example of an (almost) Hall functor -- 
the reflection functor, and
in Section \ref{AffineSection} we have still
more examples.

\subsection{Reflection functors}

\subsubsection{}

Let $(I,<,>)$ be a Cartan datum. We denote
by $\sigma_i$ an endomorphism of $\mathbb{Z}[I]$
given by the reflection with respect to
a vertex $i\in I$:
\begin{equation}\nonumber
\sigma_i ( \alpha ) = \alpha - <i,\alpha> i.
\end{equation}
It is easy to see that $\sigma_i$ preserves the inner
product $<,>$, and that $\sigma_i \sigma_i = id$. 
The subgroup of the group of endomorphisms of
$\mathbb{Z} [I]$
generated by $\{ \sigma_i \}_{i \in I}$
is called the \emph{Weyl group}.
We denote the Weyl group by $W$.

\subsubsection{}

Let $Q=(I, \Omega, \In, \Out)$ be a quiver with
the underlying Cartan datum \mbox{$(I,<,>)$}. 

A vertex $i \in I$ is called a \emph{source} if
there is no $h \in \Omega$ such that $\In (h)=i$.

A vertex $i \in I$ is called a \emph{sink} if
there is no $h \in \Omega$ such that $\Out (h)=i$.

A vertex is called \emph{admissible} if it
is either a source or a sink.
 
A reflection $\sigma_i$ is called 
\emph{admissible} if $i$ is an admissible
vertex.
 
We denote by $\sigma_i Q = (I, \Omega, \In ', \Out ' )$
a quiver with the same sets of vertices and edges
and with the maps $\In '$ and $\Out '$ defined 
as follows:
\begin{gather}\nonumber
\In ' (h) = \In (h) , \Out ' (h) = \Out (h)
\text{ if } 
\In (h) \neq i \text{ and } \Out (h) \neq i ,
\\\nonumber
\In ' (h) = \Out (h) , \Out ' (h) = \In (h)
\text{ if } 
\In (h) = i \text{ or } \Out (h) = i .
\end{gather} 
In other words, $\sigma_iQ$ is obtained from $Q$
by reversing all the arrows coming to $i$
and going out of $i$.

\subsubsection{}\label{Coxeter}

We call a sequence of vertices 
$i_1, \dots , i_K \in I$ \emph{admissible} if
for any $l \in \{1, \dots , K \}$ the vertex $i_l$
is a sink for the quiver 
$\sigma_{i_{l+1}} \dots  \sigma_{i_K} Q$.

We call $c=\sigma_{i_1} \dots \sigma_{i_{|I|}} \in W$
a \emph{Coxeter element} if $i_1, \dots , i_{|I|}$ is
an admissible sequence of vertices in which
each vertex is present exactly once (in other
words it is an ordering of vertices). It is
easy to see that if the quiver has no oriented
cycles then there exists a unique Coxeter element.

It is clear that $cQ=Q$.

\subsubsection{}\label{ReflectionsOfEpsilon}

The following proposition asserts that
the Euler cocycle $\epsilon$ 
introduced in \ref{DefinitionOfEpsilon} is 
functorial with respect to admissible reflections.

\begin{proposition}
Let $i$ be an admissible vertex for $Q$. Then
\begin{equation}\nonumber
\epsilon_{\sigma_i Q} (\sigma_i \alpha, 
\sigma_i \beta ) = 
\epsilon_Q ( \alpha, \beta ) ,
\end{equation}
for any $\alpha$, $\beta\in \mathbb{Z}[I]$.
\end{proposition}

\begin{proof}
An easy calculation using the definition of 
$\epsilon_Q$. See also Remark \ref{FactorCategories}.
\end{proof}

\subsubsection{}\label{DefinitionIN}
For an admissible $i \in I$ we denote by ${^i 
\mathcal{M}} (Q)$ 
the full subcategory of $\mathcal{M} (Q)$ 
defined as follows. Let
$i$ be a source (resp. sink) and 
$(V, x)$ be an object of 
$\mathcal{M} (Q)$. Then it is
an object of ${^i \mathcal{M}} (Q)$ if 
$\bigoplus_{h | \Out (h) =i} x_h :
V_i \rightarrow \bigoplus_{h | \Out (h) =i} V_{\In (h)}$
is injective (resp.
$\bigoplus_{h | \In (h) =i} x_h :
\bigoplus_{h | \In (h) =i} V_{\Out (h)} \rightarrow V_i$
is surjective). Let us note that the
subcategory ${^i \mathcal{M}} (Q)$ is not
abelian in general.

\begin{proposition}
\begin{enumerate}
\item\label{ExtClosed}
The subcategory ${^i \mathcal{M}} (Q)$ is \emph{\'epaisse}.
\item\label{IMInd}
$\Ob ({^i \mathcal{M}} (Q))$ contains all
indecomposable objects of $\mathcal{M} (Q)$
except for the simple object $\mathbf{P}_i$.
\end{enumerate}
\end{proposition}
\begin{proof}
Follows from the definition.
\end{proof}

\subsubsection{}

Let $i$ be an admissible vertex.
Following \cite{BernsteinGelfandPonomarev}
we define \emph{reflection functor}
$\mathcal{S}_i (Q): \mathcal{M} (Q) \rightarrow 
\mathcal{M} (\sigma_i Q)$.

If $i$ is a sink then the action of the
functor $\mathcal{S}_i (Q)$ on objects
is defined
by $\mathcal{S}_i (Q)((V, x))=
(V', x')$, where
\begin{gather}\nonumber
V'_k = V_k  \text{ if } k \neq i ,
\\\nonumber
V'_i = \Ker (
\bigoplus_{h | \In (h) =i} x_h :
\bigoplus_{h | \In (h) =i} V_{\Out (h)} 
\rightarrow V_i ), 
\\\nonumber
x'_h = x_h \text{ if } \In (h) \neq i , 
\\\nonumber
x'_h: V'_i \rightarrow V_{\Out (h)} 
\text{ is the inclusion composed 
with the projection if } 
\In (h) = i .
\end{gather}

If $i$ is a source then the action of the
functor $\mathcal{S}_i (Q)$ on objects
is defined
by $\mathcal{S}_i (Q) ( (V, x))=
(V', x')$, where
\begin{gather}\nonumber
V'_k = V_k  \text{ if } k \neq i ,
\\\nonumber
V'_i = \Coker (
\bigoplus_{h | \Out (h) =i} x_h :
V_i \rightarrow
\bigoplus_{h | \Out (h) =i} V_{\In (h)} ), 
\\\nonumber
x'_h = x_h \text{ if } \Out (h) \neq i , 
\\\nonumber
x'_h: V_{\In (h)} \rightarrow V'_i 
\text{ is the inclusion composed 
with the projection if } 
\Out (h) = i . 
\end{gather}

The action of $\mathcal{S}_i (Q)$ on morphisms is
the natural one.

\begin{proposition}\label{PropertiesOfSCal}
The functor $\mathcal{S}_i (Q)$ has the 
following properties:
\begin{enumerate}
\item\label{ImageOfSCal}
The image of $\mathcal{S}_i (Q)$ coincides
with the subcategory ${^i \mathcal{M}} (\sigma_iQ)$.
\item\label{InverseOfSCal}
The restriction of the 
functor $\mathcal{S}_i (Q)$ defines an
equivalence of categories:
$\mathcal{S}_i (Q): 
{^i \mathcal{M}} (Q)
\Tilde{\rightarrow}
{^i \mathcal{M}} (\sigma_i Q)$,
the inverse functor is 
$\mathcal{S}_i (\sigma_i Q)$.
\item\label{SCalInd}
An object $\mathbf{M} \in \Ob ({^i \mathcal{M}} (Q))$
is indecomposable if and only if
$\mathcal{S}_i (Q) \mathbf{M}$ is
indecomposable.
\item\label{SCalAndDim}
Let $\mathbf{M} \in 
\Ob({^i \mathcal{M}} (Q))$,  
$\dim_{\mathcal{M} (Q)} \mathbf{M}=\alpha$. 
Then $\dim_{\mathcal{M} (\sigma_i Q)} 
\mathcal{S}_i (\mathbf{M}) =
\sigma_i (\alpha) $.
\item\label{ExactnessOfSCal} 
If 
\begin{equation}\nonumber
0 \rightarrow \mathbf{M}_1
\rightarrow \mathbf{M}_2
\rightarrow \mathbf{M}_3
\rightarrow 0
\end{equation}
is an exact sequence of objects of
${^i \mathcal{M}} (Q)$
then
\begin{equation}\nonumber
0 \rightarrow \mathcal{S}_i (Q) (\mathbf{M}_1)
\rightarrow   \mathcal{S}_i (Q) (\mathbf{M}_2)
\rightarrow   \mathcal{S}_i (Q) (\mathbf{M}_3)
\rightarrow 0
\end{equation}
is exact in $^i \mathcal{M} (\sigma_i Q)$.
\end{enumerate}  
\end{proposition}

\begin{proof}
\ref{ImageOfSCal} and \ref{SCalAndDim}
follow from the definition of
$\mathcal{S}_i (Q)$,

\ref{InverseOfSCal} 
Let $i$ be a sink.
Then there is the following split
exact sequence
\begin{equation}\nonumber
0 \rightarrow 
\mathcal{S}_i (\sigma_i Q)\mathcal{S}_i (Q) (\mathbf{M})
\rightarrow   \mathbf{M}
\rightarrow   \mathbf{L}_i
\rightarrow 0 ,
\end{equation}
for any $\mathbf{M} \in \mathcal{M} (Q)$.
Here $\mathbf{L}_i = (V, x)$, where
\begin{gather}\nonumber
V_k = 0  \text{ if } k \neq i ,
\\\nonumber
V_i=\Coker (
\bigoplus_{h | \In (h) =i} x_h :
\bigoplus_{h | \In (h) =i} V_{\Out (h)} 
\rightarrow V_i), 
\\\nonumber
x=0 .
\end{gather}

If $\mathbf{M} \in {^i \mathcal{M}} (Q)$ then
$\mathbf{L}_i = 0$.

The proof for $i$ being a source
is analogous.

\ref{SCalInd} follows from \ref{InverseOfSCal}.

\ref{ExactnessOfSCal} follows from the Snake Lemma.

\end{proof}

\subsubsection{}\label{FactorCategories}\emph{Remark.}
Roughly speaking Proposition \ref{PropertiesOfSCal}
means that the
restriction of the functor
$\mathcal{S}_i$ to the
subcategory ${^i\mathcal{M}} (Q)$ has
associated dimension mapping 
$\dim\mathcal{S}_i=\sigma_i$ and
has properties of a Hall functor. 
In particular together with
\ref{HallAndEuler} it explains
\ref{ReflectionsOfEpsilon}.
However one should be careful, as
the subcategory ${^i \mathcal{M}} (Q)$ is
not abelian. It would be more accurate
to define the reflection functor
$\mathcal{S}_i$ not on a subcategory,
but on a quotient category (in the sense of
Serre) -- namely, on the quotient of 
$\mathcal{M}(Q)$ over the embedded category
of representations of the quiver with one vertex
$i$ and no edges.

Authors decided not to go
into the treatment of Hall mappings for
quotient categories in order to save space.
Interested reader can easily invent all
the necessary notions.
 
\subsubsection{}\label{DefinitionOfiL}

Though $\mathcal{S}_i$ is not strictly 
speaking a Hall functor, one can still define
the corresponding Hall map $\mathfrak{s}_i$
if one restricts its domain. 

We denote by ${^i \mathbf{E}}^{nil}_{\alpha}$ 
the subvariety
of $\mathbf{E}^{nil}_{\alpha}$ given by
${^i \mathbf{E}}^{nil}_{\alpha}
=\{ x \in  \mathbf{E}^{nil}_{\alpha} |
(\mathbb{C}^{\alpha}, x) \in 
{^i \mathcal{M}} (Q) \}$.

Let ${^i \mathcal{L}}_{\alpha} (Q)$ be a subspace of
$\mathcal{L}_{\alpha} (Q)$ consisting 
of all functions with
support in ${^i \mathbf{E}}^{nil}_{\alpha}$, and
${^i \mathcal{L}} (Q) =
\oplus_{\alpha\in \mathbb{Z}_+ [I]}
{^i \mathcal{L}_{\alpha} (Q)}$.
It follows from
\ref{ExtClosed} that ${^i \mathcal{L}} (Q)$
is a subalgebra of $\mathcal{L} (Q)$
with respect to $\ast$-product, and a Lie subalgebra
with respect to the Lie bracket 
\eqref{ConvolutionBracket}.

Now one can repeat the construction of the Hall
map in \ref{HallMap} for the functor 
$\mathcal{S}_i$ using Proposition
\ref{PropertiesOfSCal}, and get an
isomorphism of algebras
\begin{equation}\nonumber
\mathfrak{s}_i :
{^i \mathcal{L}} (Q) \rightarrow
{^i \mathcal{L}} (\sigma_i Q),
\end{equation}
such that
\begin{equation}\nonumber
(\mathfrak{s}_i (g)) 
([\mathcal{S}_i (\mathbf{M})]) = g ([\mathbf{M}])
\end{equation}
for $g \in {^i \mathcal{L}} (Q)$.

\section{Quiver of finite type}\label{FiniteSection}

\setcounter{subsubsection}{-1}
\subsubsection{}

Throughout this chapter $Q=(I, \Omega, \In, \Out )$
denotes a quiver of finite 
type with underlying Cartan datum $(I, <,>)$.

\subsection{Root system for a Cartan datum of finite
type}

\subsubsection{}

An element $\alpha$ of $\mathbb{Z} [I]$ is called a 
\emph{root} if $<\alpha, \alpha > =2$.
We denote the set of all roots by
$R$. Let $R_{+}$ denote the set of 
\emph{positive roots} :
$R_{+} = R \cap \mathbb{Z}_+ [I]$.  

For example elements of $I \subset \mathbb{Z} [I]$ 
are (positive) roots. 
We call $i \in I \subset R_{+}$ a \emph{simple root}.

Since the symmetric form $<,>$ is positive definite,
the set of roots is finite. Therefore the Weyl
group $W$ (which preserves the root system $R$) is
a finite group. 

\subsection{Lie algebra based on the Euler cocycle}

\subsubsection{}

We denote by $\mathfrak{n}^{\epsilon} (Q)$ 
a $\mathbb{C}$-linear space with
basis 
$\{ \Tilde{e}_{\alpha} \}_{\alpha \in R_+}$
equipped with the bilinear
bracket  

\begin{equation}\label{EpsilonCommutatorFinite}
[\Tilde{e}_\alpha, \Tilde{e}_\beta] = 
\begin{cases}
\epsilon (\alpha , \beta ) 
\Tilde{e}_{\alpha + \beta} &
\text{ if } \alpha + \beta \in R_+ , \\ 
0 & 
\text{ if } \alpha + \beta \notin R_+ , 
\end{cases}
\end{equation}
where the Euler cocycle $\epsilon$ 
is defined in \ref{DefinitionOfEpsilon}. 

\begin{theorem}[Frenkel, Kac \cite{FrenkelKac1980},
Segal \cite{Segal1981}]
The space
$\mathfrak{n}^{\epsilon} (Q)$ equipped with
the bracket \eqref{EpsilonCommutatorFinite} 
is a Lie algebra
over $\mathbb{C}$.  
\end{theorem}
\begin{proof}One has to check that the bracket is
skew-symmetric and satisfies the Jacobi
identity. Both statements follow from
the definition of $\epsilon$.
\end{proof}

\begin{numproposition}\label{eSerreFinite}
The Lie algebra
$\mathfrak{n}^{\epsilon} (Q)$ is generated by
$\{ \Tilde{e}_i \}_{i\in I}$,
the following relations hold in 
$\mathfrak{n}^{\epsilon} (Q)$ :
\begin{align}\nonumber
[\Tilde{e}_i, \Tilde{e}_j] = 0 
&\text{ if $<i,j>=0$ }, 
\\\nonumber
[\Tilde{e}_i, [\Tilde{e}_i, \Tilde{e}_j]] = 0 
&\text{ if $<i,j>=-1$ }.
\end{align}
\end{numproposition}

\begin{proof}
Induction using the fact that for any 
$\alpha \in R_+$ there exists $i\in I$
such that $(\alpha - i) \in (R_+ \cup 0)$
\cite{Bourbaki}.
\end{proof}

\subsubsection{}
Because of Proposition \ref{eSerreFinite}
we have a surjective homomorphism of Lie algebras
\begin{equation}\nonumber
\Xi^{\epsilon} : \mathfrak{n} (Q) 
\rightarrow 
\mathfrak{n}^{\epsilon} (Q)  
\end{equation}
induced by
\begin{equation}\nonumber
\Xi^{\epsilon} (e_{i}) =  \Tilde{e}_i, 
\end{equation}
where $\mathfrak{n} (Q)$ is 
introduced in \ref{DefinitionOfKM}.

It is known that $\mathfrak{n} (Q)$
is isomorphic to the nilpotent radical of
a Borel subalgebra of the semisimple Lie algebra
associated to the Dynkin graph underlying the 
quiver $Q$. In particular 
$\dim_\mathbb{C} 
\mathfrak{n} (Q)
= |R_+| =  
\dim_\mathbb{C} \mathfrak{n}^{\epsilon}(Q)$.

Comparison of dimensions implies 
that the surjective homomorphism 
$\Xi^{\epsilon}$ 
is actually an isomorphism and
$\Xi_Q =
\Xi^{\ast} \circ (\Xi^{\epsilon})^{-1}$
is a well-defined, surjective, 
$\mathbb{Z} [I]$-graded homomorphism from
$\mathfrak{n}^{\epsilon} (Q)$ to 
$\mathfrak{n}^{\ast} (Q)$.

The following commutative diagram shows all the Lie
algebras involved 
and the corresponding $\mathbb{Z} [I]$-graded 
homomorphisms

\begin{equation}\nonumber
\xymatrix{
& \mathfrak{n} (Q) 
\ar[dl]_{\Xi^{\epsilon}}
\ar[dr]^{\Xi^{\ast}} \\
\mathfrak{n}^{\epsilon} (Q)
\ar[rr]^{\Xi_Q}
&&
\mathfrak{n}^{\ast} (Q)
}
\end{equation}

Our goal is to study the map $\Xi_Q$.

\subsection{Lie algebra $\mathfrak{n}^{\ast}$
for a quiver of finite type}

\subsubsection{}\label{BGPTheorem}

The following theorem is a combination of results of
Gabriel \cite{Gabriel1972}, and Bernstein, 
Gelfand, and Ponomarev \cite{BernsteinGelfandPonomarev}.

\begin{theorem}
Let $Q$ be a quiver of finite type. 
\begin{enumerate}
\item
If $\alpha \notin R_+$ then there are no
indecomposable objects in $\mathcal{M} (Q)$
with dimension $\alpha$
(in other words, $\mathcal{T}_{\alpha}$ is empty).
\item\label{FiniteEAlpha}
There is a unique (up to an isomorphism) 
indecomposable object
$\mathbf{P}_\alpha$
in $\mathcal{M} (Q)$ with dimension
$\alpha$ for any $\alpha \in R_+$ 
In other words, 
$\mathcal{T}_{\alpha}$ consists of
a single point. We denote the characteristic
function of $[\mathbf{P}_\alpha]$ by 
$E_\alpha \in \mathcal{L}_{\alpha} (Q)$.
\item\label{GoodSequence}
For any given $\alpha \in R_+$ there
exists an admissible sequence 
$\{ i_t \}_{t=0}^{N}$
of vertices of $Q$ (depending on
$\alpha$),
such that
\begin{itemize} 
\item
$\alpha = \sigma_{i_N} 
\sigma_{i_{N - 1}} 
\dots \sigma_{i_{1}} i_{0}$,
\item
$\mathbf{P}_\alpha = \mathcal{S}_{i_N} 
\mathcal{S}_{i_{N - 1}} 
\dots \mathcal{S}_{i_{1}} 
\mathbf{P}_{i_{0}}$,
\item
$E_\alpha = \mathfrak{s}_{i_N} 
\mathfrak{s}_{i_{N - 1}} 
\dots 
\mathfrak{s}_{i_{1}} E_{i_{0}}$,
\item 
for any $t \in \{ 1, \dots , N \}$ the element
$\sigma_{i_{t}} 
\sigma_{i_{t - 1}} 
\dots \sigma_{i_{1}} 
i_{0}$ of $R$
belongs to $R_+$.
\end{itemize}
\end{enumerate} 
\end{theorem}

Note that the notation 
$\mathbf{P}_{\alpha}$ and
$E_{\alpha}$
is consistent with the notation
$\mathbf{P}_i$ and 
$E_i$ which we use above.

\subsubsection{Remark}\label{DefinitionOfn}

It follows from Theorem \ref{BGPTheorem}
that in the case of a quiver of
finite type there are finitely many non-isomorphic 
objects of $\mathcal{M} (Q)$ in each dimension.
In other words, there are finitely many
$\mathbf{G}_\alpha$-orbits in $\mathbf{E}_{\alpha}$
for each $\alpha$. Therefore $\mathcal{L}_\alpha (Q)$ 
is finite-dimensional with a basis consisting
of characteristic functions
of the orbits. Thus $\mathcal{L} (Q)$
can be thought of as a 
$\mathbb{C}$-linear space 
with a basis given by the set 
$\mathcal{X} (Q)$ of isomorphism
classes of objects of $\mathcal{M} (Q)$.

The $\ast$-product can be rewritten as
\begin{equation}\nonumber
[\mathbf{M}_1] \ast [\mathbf{M}_2] =
\sum_{[\mathbf{M}_3] \in \mathcal{X} (Q)}
n^Q_{\mathbf{M}_1, \mathbf{M}_2 ;
\mathbf{M}_3 } 
[\mathbf{M}_3],
\end{equation}
where $\mathbf{M}_1$, $\mathbf{M}_2$, $\mathbf{M}_3$
are representatives of isomorphism classes and
$ n^Q_{\mathbf{M}_1, \mathbf{M}_2 ;
\mathbf{M}_3 }$ is 
defined in \ref{DefinitionOfNMMM}.

\subsubsection{}

Our main result for quivers of finite type is
the following

\begin{theorem}\label{FiniteTheorem}
Let $Q$ be a quiver of finite type. Then
$\Xi_Q ( \Tilde{e}_\alpha ) = E_\alpha $
for any $\alpha \in R_+$.
\end{theorem}

\subsubsection{}

Before proving Theorem \ref{FiniteTheorem}
we give the following immediate Corollary:

\begin{corollary}[Ringel \cite{Ringel1990}]
\label{FiniteCorollary}
\begin{enumerate}
\item
The $\mathbb{C}$-linear space
$\mathfrak{n}^{\ast} (Q)$ coincides with the 
space of all 
$\mathbb{C}$-valued functions
on the set  $\mathcal{T}$ of isomorphism classes
of indecomposable representations.
\item
The Lie bracket in $\mathfrak{n}^{\ast} (Q)$ is given
by the Euler cocycle $\epsilon$:
\begin{equation}\nonumber
[E_\alpha , E_\beta ]=
\begin{cases} 
\epsilon (\alpha , \beta)
E_{\alpha + \beta} , 
&\text{ if } \alpha + \beta \in R_+ , \\
0 &\text{ if } \alpha + \beta \notin R_+ .
\end{cases}
\end{equation}
\end{enumerate}
\end{corollary}

\subsubsection{Remark}
Statement \ref{FiniteCorollary}
was first established by 
C. M. Ringel \cite{Ringel1990},
who used slightly different definition of
the Lie algebra $\mathfrak{n}^{\ast} (Q)$.
He found all the varieties 
$ N^Q_{\mathbf{M}_1, \mathbf{M}_2 ;
\mathbf{M}_3 }$ used in the definition
of the structure constants 
$ n^Q_{\mathbf{M}_1, \mathbf{M}_2 ;
\mathbf{M}_3 }$ (see \ref{DefinitionOfn}).

Our proof does not rely on the explicit form
of the varieties $ N^Q_{\mathbf{M}_1, \mathbf{M}_2 ;
\mathbf{M}_3 }$.

\subsection{Reflections revisited or
a proof of the Ringel theorem}

\subsubsection{}\label{DefinitionOfEij}

Let $i$ be an admissible vertex.
We denote by ${^i \mathfrak{n}}^{\ast} (Q)$
the intersection
of $\mathfrak{n}^{\ast}(Q)$
with ${^i \mathcal{L}} (Q)$ 
(see \ref{DefinitionOfiL}).
Note that $E_i \notin {^i \mathfrak{n}^{\ast}} (Q)$.

Let $j$ be a vertex connected with $i$ by an edge
(we recall that $Q$ being of finite type
is simply laced).
Let $E_{i+j} \in {^i \mathcal{L}} (Q)$ be a
function on $\mathbf{E}_{i+j}$ equal
to the characteristic function of the
$\mathbf{G}_{i+j}$-orbit, corresponding to the 
unique, up to an isomorphism, indecomposable 
representation $\mathbf{P}_{i+j}$ of
dimension $i+j \in \mathbb{Z} [I]$.

It follows from the definition of the
$\ast$-product that $E_{i+j} = [E_j, E_i]$ if
$i$ is a source, and 
$E_{i+j} = [E_i, E_j]$ if $i$ is a sink.
In particular, $E_{i+j} \in {^i \mathfrak{n}^{\ast}} (Q)$.

\begin{proposition}
\begin{enumerate}
\item\label{FiniteSAction}
Let $\alpha \in R_+$, $\alpha  \neq i$. We denote by
$E'_{\alpha} \in {^i \mathcal{L}}(Q)$  and
$E''_{\alpha} \in {^i \mathcal{L}}(\sigma_i Q)$
the characteristic functions defined in
\ref{FiniteEAlpha} for quivers $Q$ and $\sigma_i Q$
respectively. Then
\begin{equation}\nonumber
\mathfrak{s}_i (E'_{\alpha}) = E''_{\sigma_i \alpha} .
\end{equation}
\item
The subspace ${^i \mathfrak{n}}^{\ast} (Q)$ is a 
$\mathbb{Z}_+ [I]$-graded Lie subalgebra 
of $\mathfrak{n}^{\ast} (Q)$ generated by elements
$E_k$ for $k \neq i$ and $E_{i+j}$ for all
$j$ connected to $i$ by an edge.
\item\label{ReflectionIsomorphismLRS}
The restriction of $\mathfrak{s}_i$
to ${^i \mathfrak{n}}^{\ast} (Q)$ 
yields an isomorphism of Lie algebras
${^i \mathfrak{n}}^{\ast} (Q)$ and
${^i \mathfrak{n}}^{\ast} ( \sigma_iQ)$.
\end{enumerate}
\end{proposition}
\begin{proof}
Note that $\sigma_i (j) = i+j$ because
$Q$ is simply laced. Now the proposition 
follows from \ref{DefinitionOfiL}, 
\ref{SCalInd}, and \ref{SCalAndDim}.
The action of $\mathfrak{s}_i$ on the generators
of ${^i \mathfrak{n}}^{\ast} (Q)$ is given by
\begin{gather}\nonumber
\mathfrak{s}_i (E'_j) = E''_j 
\quad\text{ if $i$ is not connected with $j$ by an edge },
\\\nonumber
\mathfrak{s}_i (E'_{i+j}) = E''_j , \qquad
\mathfrak{s}_i (E'_j) = E''_{i+j}  
\qquad\text{ if $i$ is connected with $j$ by an edge },
\end{gather}
where the notation is as in \ref{FiniteSAction}.
\end{proof}

\subsubsection{}\label{DefinitionOfeij}

For an admissible vertex $i$ we denote by 
$^i \mathfrak{n}^{\epsilon}$  the
Lie subalgebra of $\mathfrak{n}^{\epsilon}$,
generated by elements 
$\{ \Tilde{e}_k \}_
{\substack{k \in I \\ k \neq i}}$ and 
$\{ \Tilde{e}_{i + j} \}_
{\substack{j \in I \\ i+j \in R_+}}$ 
(cf. \ref{DefinitionOfEij}).
Note that $i+j \in R_+$ if and only if $j$ is connected 
to $i$ by an edge. One has
$\mathfrak{n}^{\epsilon}  =
{^i \mathfrak{n}^{\epsilon}}  
\oplus \mathbb{C} \Tilde{e}_i$ as a
$\mathbb{C}$-linear space,
and 
$\Tilde{e}_{i+j} = [ \Tilde{e}_j , \Tilde{e}_i ]$
(resp. 
$\Tilde{e}_{i+j} = [ \Tilde{e}_i , \Tilde{e}_j ]$)
if $i$ is a source (resp. a sink).

\subsubsection{}\label{FiniteSEAction}

We denote by $\mathfrak{s}^{\epsilon}_i$ 
a $\mathbb{C}$-linear map 
$^i \mathfrak{n}^{\epsilon} (Q) \rightarrow 
{^i \mathfrak{n}^{\epsilon} (\sigma_i Q)} $
given by 
$\mathfrak{s}^{\epsilon}_i (\Tilde{e}'_{\alpha} ) = 
\Tilde{e}''_{\sigma_i (\alpha)}$,
where  
$\Tilde{e}'_{\alpha}$, $\Tilde{e}''_{\alpha}$ are
basic vectors of 
${^i \mathfrak{n}}^{\epsilon} (Q)$ and 
${^i \mathfrak{n}^{\epsilon} (\sigma_i Q)}$
respectively. In particular, for the generators
\begin{gather}\nonumber
\mathfrak{s}^{\epsilon}_i 
(\Tilde{e}'_j) = \Tilde{e}''_j 
\qquad \text{ if $i$ is not 
connected with $j$ by an edge },
\\\nonumber
\mathfrak{s}^{\epsilon}_i 
(\Tilde{e}'_{i+j}) = \Tilde{e}''_j , \quad
\mathfrak{s}^{\epsilon}_i 
(\Tilde{e}'_j) = \Tilde{e}''_{i+j}  
\qquad \text{ if $i$ is connected with $j$ by an edge }.
\end{gather}

\begin{proposition}
\label{ReflectionIsomorphismEpsilon}
The map $\mathfrak{s}^{\epsilon}_i$ 
is well-defined and is an isomorphism of
Lie algebras.
\end{proposition}
\begin{proof}
Follows from the fact that 
$\sigma_i (R_+ \backslash \{ i \}) = 
R_+ \backslash \{ i \}$
and from \ref{ReflectionsOfEpsilon}.
\end{proof}

\subsubsection{}

It follows from \ref{DefinitionOfEij},
\ref{DefinitionOfeij} that
$\Xi_Q 
(^i \mathfrak{n}^{\epsilon} (Q) ) \subset
{^i \mathfrak{n}^{\ast} (Q)}$,
and that $\Xi_Q (\Tilde{e}_{i+j}) =
E_{i+j}$.

Using Propositions \ref{ReflectionIsomorphismLRS}
and \ref{ReflectionIsomorphismEpsilon} we
get the following diagram
of Lie algebra homomorphisms:

\begin{equation}\label{FiniteSDiagram}
\xymatrix{
^i \mathfrak{n}^{\epsilon} (Q)
\ar[d]^{\Xi_Q}
\ar[r]^{\mathfrak{s}^{\epsilon}_i} & 
^i \mathfrak{n}^{\epsilon} (\sigma_i Q)
\ar[d]^{\Xi_{\sigma_i Q}} \\
^i \mathfrak{n}^{\ast} (Q) 
\ar[r]^{\mathfrak{s}_i} &
^i \mathfrak{n}^{\ast} (\sigma_i Q) ,
}
\end{equation}
where the horizontal maps are
isomorphisms. The diagram \eqref{FiniteSDiagram}
is commutative on generators
$\Tilde{e}_k$, $\Tilde{e}_{i+j}$
of $^i\mathfrak{n}^{\epsilon} (Q)$ and,
therefore, is commutative.

\subsubsection{Proof of Theorem \ref{FiniteTheorem}}

Let $G (Q) \subset R_+$ be the set of all
$\alpha \in R_+$ such that 
$\Xi_Q  (\Tilde{e}_{\alpha} ) = E_{\alpha}$.
Then $I \subset G (Q)$ for any $Q$, and
it follows from \ref{FiniteSAction}, 
\ref{FiniteSEAction}, and the commutative diagram
\eqref{FiniteSDiagram} that 
$\sigma_i (G (Q) \backslash \{ i \} ) = 
G (\sigma_i Q) \backslash \{ i \}$
for any admissible vertex $i$.
Using the sequence of reflections 
\eqref{GoodSequence}
one obtains that $\alpha \in G (Q)$ for
any $\alpha \in R_+$. Therefore 
$\Xi_Q (\Tilde{e}_{\alpha} ) = E_{\alpha}$
for any $\alpha \in R_+$ ,
which proves the theorem.
\qed

\section{Quiver of affine type}\label{AffineSection}

\setcounter{subsubsection}{-1}
\subsubsection{}

In this chapter we consider 
the Lie algebra $\mathfrak{n}^{\ast} (Q)$ in
the case of $Q$ being of affine type.

Throughout the chapter $Q=(I, \Omega, \In, \Out )$
denotes a quiver of affine 
type with underlying Cartan datum $(I, <,>)$.
In particular the bilinear form $<,>$ is
positive semi-definite.

\subsection{Root system for a Cartan datum of affine
type}

\subsubsection{}\label{AffineRoots}

An element $\alpha \in \mathbb{Z} [I]$ is called a 
\emph{root} if $<\alpha, \alpha > \leq 2$ and
$\alpha \neq 0$.
We denote the set of all roots by
$R$. Let $R_{+}$ denote the set of 
\emph{positive roots} :
$R_{+} = R \cap \mathbb{Z}_{+} [I]$.  

For example elements of $I \subset \mathbb{Z} [I]$ 
are (positive) roots. 
We call $i \in I \subset R_{+}$ a \emph{simple}
root.

It is easy to see that if $\alpha \in R$ then
$<\alpha , \alpha > \in \{ 0, 2 \}$. A root
$\alpha$ such that  
$<\alpha , \alpha > =  2$ (resp. 
$<\alpha , \alpha > =  0$ )
is called \emph{real} (resp. \emph{imaginary})
root. We denote the set of all positive
real (resp. positive imaginary) roots by $R^{\re}_+$ 
(resp. $R^{\im}_+$). 
So $R_+ = R^{\re}_+ \cup R^{\im}_+$.

It is known that 
$R^{\im}_+ = (\mathbb{Z}_+ \backslash 
\{ 0 \}) \delta$, where
$\delta \in R_+$. We call $\delta$ the
\emph{first imaginary root}.

The components of $\delta$ are given in the 
table of affine Dynkin graphs (Figure 2).

\subsubsection{}\label{ExtendingVertex}

A vertex $p \in I$ is called an \emph{extending} vertex
if $\delta_p =1$, where $\delta$ is the first imaginary
root. Such a vertex always exists, but it is not unique
in general. If one removes from $Q$ 
an extending vertex together with adjacent edges one
gets an irreducible quiver of finite type.

\subsection{Lie algebra based on the Euler cocycle}

\subsubsection{}\label{AffineNEpsilon}

We consider the following set of
$\mathbb{C}$-linear spaces:
\begin{itemize}
\item
for each positive real root $\alpha$ 
a one-dimensional $\mathbb{C}$-linear space 
$\mathfrak{n}_{\alpha}^{\epsilon} (Q) =
\mathbb{C} \Tilde{e}_{\alpha} $ 
with generator $\Tilde{e}_{\alpha}$,
\item 
for each positive imaginary root $n \delta$
a $\mathbb{C}$-linear space 
$\mathfrak{n}_{n \delta}^{\epsilon} (Q) =
\mathbb{C} [I] / \mathbb{C} \delta$,
where we consider $\delta$ as an
element of $R_+ \subset \mathbb{Z} [I]
\subset \mathbb{C} [I]$.
For $h \in \mathbb{C} [I]$ we denote by 
$h (n)$ the image of $h$ 
under the natural projection map
$\mathbb{C} [I] \rightarrow
\mathfrak{n}_{n \delta}^{\epsilon} (Q)$.
\end{itemize}

We denote by 
$\mathfrak{n}^{\epsilon} (Q)$ the following
$R_+$-graded $\mathbb{C}$-linear space:
\begin{equation}\nonumber
\mathfrak{n}^{\epsilon} (Q) = 
\bigoplus_{\alpha \in R_+}
\mathfrak{n}_{\alpha}^{\epsilon} (Q) ,
\end{equation}
equipped with the bilinear bracket

\begin{equation}\label{EpsilonCommutatorAffine}
\begin{split}
[\Tilde{e}_\alpha, \Tilde{e}_\beta] &= 
\begin{cases}
\epsilon (\alpha , \beta ) 
\Tilde{e}_{\alpha + \beta} &
\text{ if } \alpha + \beta \in R^{\re}_+ , \\ 
\epsilon (\alpha , \beta ) 
\alpha (k) &
\text{ if } \alpha + \beta = k \delta , \\ 
0 & 
\text{ if } \alpha + \beta \notin R_+ , 
\end{cases}
\\
[h (n), \Tilde{e}_{\alpha} ] =
-[\Tilde{e}_{\alpha} , h(n)] &=
\epsilon (n \delta , \alpha ) <h, \alpha>  
\Tilde{e}_{\alpha + n \delta} , 
\\
[h (n), h (m) ] &= 0 ,
\end{split}\end{equation}
where $h \in \mathbb{C} [I]$, 
$\epsilon$ is the Euler cocycle 
(see \ref{DefinitionOfEpsilon}).
Note that $<h' , \alpha > = <h'' , \alpha>$ if
$h' (n) = h'' (n)$   because 
$\delta \in \Ker<,>$. Therefore the bracket
\eqref{EpsilonCommutatorAffine} is well-defined. 

\begin{theorem}
The $\mathbb{C}$-linear space
$\mathfrak{n}^{\epsilon} (Q)$ equipped with
the bracket \eqref{EpsilonCommutatorAffine} 
is a Lie algebra.  
\end{theorem}
\begin{proof}One has to check that the bracket is
skew-symmetric and satisfies the Jacobi
identity. Both statements follow (after lengthy,
but straightforward  calculations) from
the definition of $\epsilon$.
A more conceptual proof is to consider 
$\mathfrak{n}^{\epsilon}$ as embedded
into an algebra of vertex operators, associated
with the quiver $Q$ (see \cite{Frenkel1985,
MoodyRaoYokonuma1990}).
\end{proof}

\begin{numproposition}\label{eSerreAffine}
The Lie algebra $\mathfrak{n}^{\epsilon} (Q)$ 
is generated by $\{ \Tilde{e}_i \}_{i\in I}$,
the following relations hold in 
$\mathfrak{n}^{\epsilon} (Q)$ : 
\begin{equation}\nonumber
(\ad ( \Tilde{e}_i ))^{1-a_{ij}} 
\Tilde{e}_j = 0 
\text{ if $i \neq j$ },  
\end{equation}
where $\ad (x) y = [x,y]$, and $a_{ij}=<i,j>$.
\end{numproposition}

\begin{proof}
Induction on $n ( \alpha )$, where
$n ( \alpha ) = \sum_{i\in I} \alpha_i$
for a root $\alpha$.
\end{proof}

\subsubsection{}\label{DefinitionOfXiAffine}
Because of Proposition \ref{eSerreAffine}
we have a surjective 
$R_+$-graded homomorphism of Lie algebras:
\begin{equation}\nonumber
\Xi^{\epsilon} : \mathfrak{n} (Q) 
\rightarrow 
\mathfrak{n}^{\epsilon} (Q)  
\end{equation}
induced by
\begin{equation}\nonumber
\Xi^{\epsilon} (e_{i}) =  \Tilde{e}_i, 
\end{equation}
where $\mathfrak{n} (Q)$ is the Lie algebra
introduced in \ref{DefinitionOfKM}.

It is known that 
\begin{equation}\nonumber
\dim_\mathbb{C} 
\mathfrak{n}_{\alpha} (Q)
= 
\begin{cases}
1 & \text{ if $\alpha \in R^{\re} $  } \\
|I| - 1 & \text{ if $\alpha \in R^{\im}$ } 
\end{cases}
\Biggr\}=  
\dim_\mathbb{C} \mathfrak{n}^{\epsilon}_{\alpha} (Q).
\end{equation}

Comparison of dimensions implies 
that the surjective homomorphism 
$\Xi^{\epsilon}$ 
is actually an isomorphism and
$\Xi_Q =
\Xi^{\ast} \circ (\Xi^{\epsilon})^{-1}$
is a well-defined, surjective, 
$\mathbb{Z} [I]$-graded homomorphism from
$\mathfrak{n}^{\epsilon} (Q)$ to 
$\mathfrak{n}^{\ast} (Q)$.

As in the case of a quiver of finite type
we want to study 
the map $\Xi_Q$ and,
in particular, its image. There are two (related)
differences with the finite case.
First, the description of indecomposable objects
in the category $\mathcal{M} (Q)$ is more complicated.
In particular, there are infinite families of
indecomposable objects in some graded dimensions. 
Second, some of the graded 
components of $\mathfrak{n}^{\epsilon} (Q)$ are not
one-dimensional.

Because of these complications we consider
special cases first, and then build the general 
affine case on them.

\subsection{Jordan quiver}

\subsubsection{}

In this subsection we consider quiver $C_1$ 
that consists
of one vertex connected by an edge with itself:
$C_1 = (\{ \delta \} , \{ e \}, \In , \Out )$, where
$\In (e) = \Out (e) =\delta$. One can draw $C_1$ 
as follows:
 
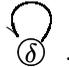
\begin{equation}\nonumber
\xygraph{[]*+[o][F-]{\delta} 
:@(ur,ul) [] *+[o][F-]{\delta}} .
\end{equation}

Though $C_1$ is not a Dynkin quiver, the definition
of the category of nilpotent 
representations of $C_1$ makes perfect sense.
Objects of $\mathcal{M} (C_1)$
are pairs $(V, x)$ consisting of a
$\mathbb{C}$-linear space $V$ and
a nilpotent $x\in \Hom_{\mathbb{C}} (V,V)$.

\subsubsection{}

An indecomposable object of $\mathcal{M} (C_1)$
with dimension $n$ is 
isomorphic to the pair 
$(\mathbb{C}^n, J_n ) = \mathbf{J}_n$,
where $J_n$ is the Jordan block
\begin{equation}\label{Jordan}
J_n  =
\begin{pmatrix}
0 & 1 & 0 & \dots & 0 & 0 \\
0 & 0 & 1 & \dots & 0 & 0 \\
\hdotsfor{6} \\
0 & 0 & 0 & \dots & 1 & 0 \\
0 & 0 & 0 & \dots & 0 & 1 \\
0 & 0 & 0 & \dots & 0 & 0 
\end{pmatrix} .
\end{equation}

\begin{numproposition}\label{HomJJ}
$\dim_{\mathbb{C}} \Hom_{\mathcal{M} (C_1)}
(\mathbf{J}_n , \mathbf{J}_n) = n$.
\end{numproposition}

\begin{proof}
An element $f$ of $\Hom_{\mathcal{M} (C_1)}
(\mathbf{J}_n , \mathbf{J}_n)$ is, by definition, 
an $n \times n$ matrix such that $f J_n = J_n f$. 
It follows that
\begin{equation}
f  =
\begin{pmatrix}
a_1 & a_2 & a_3 & \dots & a_{n-1} & a_n \\
0 & a_1 & a_2 & \dots & a_{n-2} & a_{n-1} \\
\hdotsfor{6} \\
0 & 0 & 0 & \dots & a_2 & a_3 \\
0 & 0 & 0 & \dots & a_1 & a_2 \\
0 & 0 & 0 & \dots & 0 & a_1 
\end{pmatrix} 
\end{equation}
for some set $\{ a_i \in \mathbb{C} \}_{i=1}^{n}$.
\end{proof}

\subsubsection{}

We
call the vertex $\delta$ of $C_1$ the first imaginary 
root
when $\delta$ is considered as an element of 
$\mathbb{Z}[\delta]$.

One has $e_{C_1} (\alpha , \beta) = 0$
for any $\alpha, \beta \in \mathbb{Z}[ \delta ]$.
We put $<\alpha, \beta >_{C_1} = 
e_{C_1} (\alpha , \beta) + e_{C_1} 
(\beta , \alpha) = 0$ and
$\epsilon_{C_1} (\alpha , \beta )=
(-1)^{e_{C_1} (\alpha , \beta )} = 1$.

We use the notation $r$ (rotation transformation,
cf. \ref{Rotation}) for the identity transformation of
$\mathbb{Z}[ \delta ]$.

\subsection{Kronecker quiver}

\subsubsection{}

In this section we consider a quiver $K$ that consists
of two vertices connected by two edges, which 
are directed in 
the same way. So $K=(\{ 0,1 \} , 
\{ a,b \}, \In , \Out )$, where
$\In (a) = \In (b) =1$, $\Out (a) = \Out (b) =0$:

\begin{equation}\nonumber
\xygraph{[]
*+[o][F-]{0} :@/^/^a [r]*+[o][F-]{1}
                  [l] 
*+[o][F-]{0} :@/_/_b [r]*+[o][F-]{1}
}
\end{equation}

The set of indecomposable  representations 
of $K$ was described by  
Kronecker \cite{Kronecker1890}.

\subsubsection{}

It is easy to see that 
$\alpha \in \mathbb{Z} \oplus\mathbb{Z}$  
is a positive real root 
(resp. positive imaginary root) 
if $\alpha =(n+1, n)$ or $\alpha =(n, n+1)$ 
(resp. $\alpha =(n+1, n+1)$) for some 
$n \in \mathbb{Z}_+$.
In particular $\alpha_0 = (1,0)$ and 
$\alpha_1 =(0,1)$ are
the simple roots, and $\delta = (1,1)$ 
is the first imaginary root.

The Lie bracket \eqref{EpsilonCommutatorAffine}
in $\mathfrak{n}^{\epsilon} (K)$ can be 
written as follows:

\begin{equation}\label{KBracket}
\begin{split}
[\Tilde{e}_{(n,n+1)}, \Tilde{e}_{(m,m+1)}] &= 0 ,\\
[\Tilde{e}_{(n+1,n)}, \Tilde{e}_{(m+1,m)}] &= 0 ,\\
[\alpha_1 (n) , \alpha_1 (m)] &= 0 ,\\
[\alpha_1 (n) , \Tilde{e}_{(m,m+1)}] &= 
2 (-1)^{n} \Tilde{e}_{(m+n,m+n+1)} ,\\
[\alpha_1 (n) , \Tilde{e}_{(m+1,m)}] &= 
2 (-1)^{n+1} \Tilde{e}_{(m+n+1,m+n)} ,\\
[\Tilde{e}_{(n,n+1)}, \Tilde{e}_{(m+1,m)}] &= 
(-1)^{n+m} \alpha_1 (m+n+1).
\end{split}
\end{equation}

Note that $\alpha_0 (n) = - \alpha_1 (n)$.

\subsubsection{}

According to Section \ref{DefinitionOfQmod} a 
representation 
$\mathbf{M}$ of the quiver $K$ is a 
quadruple $((V_0 , V_1 ), (x_a , x_b ))$
consisting  of two 
$\mathbb{C}$-linear spaces 
$V_0$ and $V_1$ and two linear maps 
$x_a , x_b \in \Hom_{\mathbb{C}} (V_0 , V_1)$.
The dimension of the representation $\mathbf{M}$ is 
given by
$\dim_{\mathcal{M} (K)} 
\mathbf{M} = (\dim_\mathbb{C} V_0 , 
\dim_\mathbb{C} V_1)$.

\subsubsection{}
Let us describe the
structure of the category $\mathcal{M} (K)$.
There are two simple objects
$\mathbf{U}^0_0 = ((\mathbb{C}, 0),
(0, 0))$ and
$\mathbf{U}^1_0 = ((0, \mathbb{C}),
(0, 0))$. If one applies repeatedly
the reflection functors $\mathcal{S}_0$
and $\mathcal{S}_1$ to the simple objects
one obtains two families of indecomposable objects
$\mathbf{U}^0_n$ and $\mathbf{U}^1_n$, given by 
\begin{gather}\nonumber
\mathbf{U}^0_n = ((\mathbb{C}^{n+1}, \mathbb{C}^n),
(A, B)), 
\quad \dim_{\mathcal{M} (K)} 
\mathbf{U}^0_n = (n+1,n), 
\\\nonumber
\mathbf{U}^1_n = ((\mathbb{C}^n, \mathbb{C}^{n+1}),
(A^t, B^t)), 
\quad \dim_{\mathcal{M} (K)} 
\mathbf{U}^1_n = (n,n+1), 
\end{gather}
where
\begin{equation}\label{AB}
A=
\begin{pmatrix}
1 & 0 & \dots & 0 & 0 \\
0 & 1 & \dots & 0 & 0 \\
\hdotsfor{5} \\
0 & 0 & \dots & 0 & 0 \\
0 & 0 & \dots & 1 & 0 
\end{pmatrix} , 
\quad
B=
\begin{pmatrix}
0 & 1 & \dots & 0 & 0 \\
0 & 0 & \dots & 0 & 0 \\
\hdotsfor{5} \\
0 & 0 & \dots & 1 & 0 \\
0 & 0 & \dots & 0 & 1 
\end{pmatrix} .
\end{equation}

Note that the set of graded dimensions of the 
indecomposables
obtained in this way coincides with the set of 
positive real roots.

\subsubsection{}\label{CategoryC}

In this section we use the argument and notation 
borrowed
(with slight adjustments) from 
\cite[Section 2]{Lusztig1992}.

The category $\mathcal{M} (K)$ is isomorphic 
to the following
category $\mathcal{C}$: an object 
of $\mathcal{C}$ is
a triple $((V_0,V_1),\Delta)$ 
consisting of a pair $(V_0, V_1)$
of $\mathbb{C}$-linear spaces, and 
a $\mathbb{C}$-linear
map $\Delta : V_0 
\otimes_{\mathbb{C}} \mathbb{C}^2
\rightarrow V_1$.
A morphism from $((V_0,V_1),\Delta)$
to $((V'_0,V'_1),\Delta ')$
is a pair $(\phi_0,\phi_1)$ of
$\mathbb{C}$-linear maps
$\phi_i : V_i \rightarrow V'_i$ such that
$\phi_1 \circ \Delta =
\Delta' \circ (\phi_0 \otimes 
\Id_{\mathbb{C}^2})$. 
Let us introduce the following notation: 
$\Delta_e (v)=\Delta (v \otimes e)$, where 
$e \in \mathbb{C}^2$, $v \in V_0$.
The isomorphism of categories
$\mathcal{C}$ and $\mathcal{M} (K)$ is given by a 
functor $\mathcal{J}$ that takes 
$((V_0,V_1), \Delta)$ to
$((V_0,V_1), (
\Delta_{\bigl( 
\begin{smallmatrix}1\\0
\end{smallmatrix}\bigr)},
\Delta_{\bigl(
\begin{smallmatrix}0\\1
\end{smallmatrix}\bigr)}))$,
and acts naturally on morphisms.

Let $\mathbf{M}=((V_0,V_1),\Delta)
\in \Ob (\mathcal{C})$. We denote
by $\Spec \mathbf{M}$ (\emph{spectrum} of
$\mathbf{M}$) the set of all lines
$z\in \mathbb{CP}^1$ such that
for $e \in z$, $e \neq 0$ the map $\Delta_e :
V_0 \rightarrow V_1$ is not
an isomorphism.

Let ${^0\mathcal{C}}$ be the following
full subcategory of $\mathcal{C}$.
An object $\mathbf{M}$ of $\mathcal{C}$
belongs to $\Ob ({^0\mathcal{C}})$ if
$\Spec \mathbf{M}$ is a finite 
set. In particular, if 
$\mathbf{M}= ((V_0,V_1),\Delta)
\in \Ob ({^0\mathcal{C}})$
then $\dim_{\mathbb{C}} V_0 =  
\dim_{\mathbb{C}} V_1$, that is
$\dim_{\mathcal{M} (K)} \mathcal{J}
(\mathbf{M}) \in R_+^{\im}$.

One can prove (cf. \cite{Lusztig1992})
that 
${^0\mathcal{C}} = \oplus_{z\in 
\mathbb{CP}^1} 
{^0\mathcal{C}}_z$ (direct coproduct 
of abelian categories), where
${^0\mathcal{C}}_z$ is a full subcategory
of ${^0\mathcal{C}}$ whose objects are
$\mathbf{M} \in \Ob (
{^0 \mathcal{C}})$ such that
$\Spec \mathbf{M}=z$. 

The category ${^0\mathcal{C}}_z$ for
any $z\in \mathbb{CP}^1$ is
equivalent to the category $\mathcal{M} (C_1)$
of nilpotent representations of the Jordan quiver $C_1$.
Let us describe the equivalence functor
$\mathcal{I}_z:
\mathcal{M} (C_1) \rightarrow {^0\mathcal{C}}_z$.
We choose two elements $e$, $e' \in \mathbb{C}^2$
such that $e \neq 0$, $e' \neq 0$, $e \in z$,
$e' \notin z$. Then the functor 
$\mathcal{I}_z$ is given by the following action on 
objects
\begin{gather}\nonumber
\mathcal{I}_z ((V,x)) =
((V,V), \Delta),
\\\nonumber
\text{ where }
\Delta ( v \otimes (ce+c'e') )=cxv +c'v ,
\end{gather}
and the natural action on morphisms.
The equivalence class of $\mathcal{I}_z$
does not depend on the choice of $e$, $e'$.

We denote by $\mathcal{H}_z$
the composition functor
$\mathcal{H}_z = \mathcal{J} \circ
\mathcal{I}_z:
\mathcal{M} (C_1) \rightarrow \mathcal{M} (K)$.
For example, if $z=(1:0)$ the functor 
$\mathcal{H}_{(1:0)}$ is given by
the following action on objects
\begin{equation}\nonumber
\mathcal{H}_{(1:0)} ((V, x)) = 
((V, V), (x, \Id_V ))
\end{equation}
and the natural action on morphisms (here 
we have chosen
$e=(1,0)$, $e'=(0,1)$).

The functor $\mathcal{H}_z$ is a Hall functor 
for any $z \in\mathbb{CP}^1$. 

It follows that $\mathcal{H}_z (\mathbf{J}_n )$ is 
indecomposable in $\mathcal{M} (K)$ of 
dimension $n\delta$, and all indecomposables
in $\Ob (\mathcal{J}({^0 \mathcal{C}}))$ are of 
this form.
One can prove (cf. \cite{Lusztig1992})
that if $\mathbf{M}
\in \Ob (\mathcal{M} (K))$ is
indecomposable then either 
$\mathbf{M}$ is isomorphic to an object, 
obtained from 
a simple one by repeated applications
of reflection functors, or
$\mathbf{M}$ belongs to 
$\Ob (\mathcal{J} ({^0\mathcal{C}}))$.
Therefore any indecomposable object
of $\mathcal{M} (K)$ of dimension
$n\delta$ belongs to
$\Ob (\mathcal{J} ({^{0}\mathcal{C}}))$,
and we can use the map  
\begin{gather}\nonumber
\Spec : 
\mathcal{T}_{n \delta}
\rightarrow 
\mathbb{CP}^1 ,
\\\nonumber
\Spec ([\mathbf{M}])=
\Spec \mathbf{M}
\end{gather}
to identify $\mathcal{T}_{n\delta}$
and $\mathbb{CP}^1$. Abusing notation we
just write
$\mathcal{T}_{n\delta} = \mathbb{CP}^1$.
In this way $\mathcal{T}_{n\delta}$ is
endowed with a structure of an algebraic
variety.

\subsubsection{}

We summarize the above discussion in
the following proposition
which describes the set $\mathcal{T}$ of
isomorphism classes of 
indecomposable objects of $\mathcal{M} (K)$.

\begin{proposition}
One has:
\begin{enumerate}
\item
if $\alpha \notin R_+$ then  
$\mathcal{T}_\alpha = \emptyset$,
\item 
if $\alpha \in R^{\re}_+$
then $\mathcal{T}_\alpha = {pt}$,
the only (up to an isomorphism) indecomposable 
representation
with dimension $\alpha$ being
$\mathbf{U}^0_n$ (resp. $\mathbf{U}^1_n$)
for $\alpha=(n+1,n)$ (resp. $\alpha=(n,n+1)$),
\item  
if $\alpha \in R^{\im}_{+}$
then $\mathcal{T}_\alpha = \mathbb{CP}^1$;
the indecomposable representation
of dimension $\alpha = (n,n)$ corresponding
to $z \in \mathbb{CP}^1$
is $\mathcal{H}_z (\mathbf{J}_n )$.
\end{enumerate}
\end{proposition}

\subsubsection{}\label{SpecMap}

Let $\mathbf{E}_{\alpha}^{ind}=
\{ ( x_a, x_b ) \in
\mathbf{E}_{\alpha} \ | \ 
((\mathbb{C}^{\alpha_0} , \mathbb{C}^{\alpha_1}),
( x_a , x_b )) \text{ is indecomposable } \}$.
Abusing notation we denote by 
$\Spec : \mathbf{E}_{n\delta}^{ind}
\rightarrow \mathbb{CP}^1$ a
regular map obtained by assigning to 
element $(x_a,x_b)\in\mathbf{E}_{n\delta}^{ind}$
the point $z\in \mathbb{CP}^1$ such that
$((\mathbb{C}^n,\mathbb{C}^n),(x_a,x_b))$ is
isomorphic to $\mathcal{H}_z ( \mathbf{J}_n )$.
The map 
$\Spec: \mathbf{E}^{ind}_{n\delta}\to \mathbb{CP}^1$
is a fibration with a constant fiber.

\subsubsection{}

Let us introduce the following
functions on the set 
of isomorphism classes of 
indecomposable objects of $\mathcal{M} (K)$
(which by our usual abuse of notation
are also $\mathbf{G}_{\alpha}$-equivariant functions on
$\mathbf{E}_{\alpha}$):
\begin{equation}\nonumber
E_\alpha = \theta (\mathcal{T}_\alpha ) =
\theta ( \mathbf{E}^{ind}_\alpha )
\in \mathcal{L}_{\alpha} (K),
\end{equation}
where $\theta$ denotes the characteristic 
function of a set.
Note that if $\alpha \in R_+^{\im}$ then $E_{\alpha}$ 
is the constant
function equal to $1$ on $\mathcal{T}_{\alpha}=
\mathbb{CP}^1$.

\subsubsection{}

The following proposition gives some $\ast$-products,
which we use later. Due to Proposition
\ref{Reidtmann} we are only interested in
the restrictions of the $\ast$-products to
the set of indecomposable representations.

\begin{proposition}\label{KProducts}
\begin{align}\label{KProducts01nn1}
E_{(0,1)} \ast E_{(n,n-1)}|_{\mathbf{E}^{ind}_{(n,n)}} &=
E_{(n,n)},
\\\label{KProductsnn101}
E_{(n,n-1)} \ast E_{(0,1)}|_{\mathbf{E}^{ind}_{(n,n)}} &= 
0,
\\\label{KProducts10n1n}
E_{(1,0)} \ast E_{(n-1,n)}|_{\mathbf{E}^{ind}_{(n,n)}} &= 
0,
\\\label{KProductsn1n10}
E_{(n-1,n)} \ast E_{(1,0)} |_{\mathbf{E}^{ind}_{(n,n)}} &=
E_{(n,n)},
\\\label{KProducts01nn}
E_{(0,1)} \ast E_{(n,n)}|_{\mathbf{E}^{ind}_{(n,n+1)}} &=
2E_{(n,n+1)},
\\\label{KProductsnn01}
E_{(n,n)} \ast E_{(0,1)}|_{\mathbf{E}^{ind}_{(n,n+1)}} &= 
0,
\\\label{KProducts10nn}
E_{(1,0)} \ast E_{(n,n)}|_{\mathbf{E}^{ind}_{(n+1,n)}} &= 
0,
\\\label{KProductsnn10}
E_{(n,n)} \ast E_{(1,0)}|_{\mathbf{E}^{ind}_{(n+1,n)}} &= 
2E_{(n+1,n)}.
\end{align}
\end{proposition}
\begin{proof}
\ref{KProductsnn101} and \ref{KProductsnn01}
follow from the fact that 
$\mathbf{U}^1_0$ is projective.

\ref{KProducts10n1n} and \ref{KProducts10nn}
follow from the fact that 
$\mathbf{U}^0_0$ is injective.

Let us prove \ref{KProducts01nn1}.
We need to find the value of the function
$E_{(0,1)} \ast E_{(n,n-1)}$ on the 
isomorphism class of 
the representation $\mathcal{H}_z (\mathbf{J}_n)$. 
We give below the calculation for
$z=(1:0)$.
The case of arbitrary $z \in \mathbb{CP}^1$ is
completely analogous, and, moreover, the
result does not depend on $z$.
The value of
$E_{(0,1)} \ast E_{(n,n-1)}$ on the 
isomorphism class of 
$\mathcal{H}_{(1:0)} (\mathbf{J}_n)$
is equal, by the definition of the $\ast$-product,
to the Euler characteristic of the variety
$N_{\mathbf{U}_0^1, \mathbf{U}_{n-1}^0, 
\mathcal{H}_{(1:0)} (\mathbf{J}_n)}$
of all subrepresentations $\mathbf{V}$ of 
$\mathcal{H}_{(1:0)} (\mathbf{J}_n)$
such that $\mathbf{V}$ is isomorphic to
$\mathbf{U}_0^1$ and
$\mathcal{H}_{(1:0)} (\mathbf{J}_n) / \mathbf{V}$
is isomorphic to
$\mathbf{U}_{n-1}^0$.
Since $\mathcal{H}_{(1:0)} (\mathbf{J}_n)=
((\mathbb{C}^n,\mathbb{C}^n),
(J_n,\Id_{\mathbb{C}^n}))$ it follows 
that $N_{\mathbf{U}_0^1, \mathbf{U}_{n-1}^0; 
\mathcal{H}_{(1:0)} (\mathbf{J}_n)}$
coincides with the set of 
all lines $l\in \mathbb{C}^n$ such that
$((\mathbb{C}^n,\mathbb{C}^n / l ),
(p \circ J_n, p ))$ is indecomposable,
where $p: \mathbb{C}^n \rightarrow
\mathbb{C}^n / l$ is the projection.
The object $((\mathbb{C}^n,\mathbb{C}^n / l ),
(p \circ J_n, p ))$ is indecomposable
if and only if $l \notin \im J_n$. Therefore
$N_{\mathbf{U}_0^1, \mathbf{U}_{n-1}^0; 
\mathcal{H}_{(1:0)} (\mathbf{J}_n)}=\mathbb{C}^{n-1}$, 
and we get \ref{KProducts01nn1}.

Statement \ref{KProductsn1n10} is dual to 
\ref{KProducts01nn1}.

Let us prove \ref{KProducts01nn}.
We need to find the value of the function
$E_{(0,1)} \ast E_{(n,n)}$ on the 
isomorphism class of 
the representation $\mathbf{U}_n^1$. This value 
is equal, by the definition of the $\ast$-product,
to the Euler characteristic of a variety 
$X$ of all subrepresentations $\mathbf{V}$ of 
$\mathbf{U}_n^1$
such that $\mathbf{V}$ is isomorphic to
$\mathbf{U}_0^1$ and
$\mathbf{U}_n^1 / \mathbf{V}$
is isomorphic to
$\mathcal{H}_z (\mathbf{J}_n)$ for some 
$z\in \mathbb{CP}^1$.  
Similar to \ref{SpecMap} one can use the $\Spec$ map
to prove that the variety $X$ is a fibration
over $\mathbb{CP}^1$ with the
fiber over $z$ equal to the variety 
$N_{\mathbf{U}_0^1, \mathcal{H}_z (\mathbf{J}_n);
\mathbf{U}_{n}^1}$
of all subrepresentations $\mathbf{V}$ of
$\mathbf{U}_{n}^1$
such that $\mathbf{V}$ is isomorphic to
$\mathbf{U}_0^1$ and
$\mathbf{U}_{n}^1 / \mathbf{V}$
is isomorphic to
$\mathcal{H}_z (\mathbf{J}_n)$.
Let us consider the case $z=(1:0)$.
Since $\mathbf{U}_n^1 =
((\mathbb{C}^n,\mathbb{C}^{n+1}), (A^t,B^t))$,
where $A$ and $B$ are as in \eqref{AB}, and
$\mathcal{H}_{(1:0)} (\mathbf{J}_n)=
((\mathbb{C}^n,\mathbb{C}^n),
(J_n,\Id_{\mathbb{C}^n}))$, it follows 
that 
$\mathbf{V}=\{ (((0, \dots ,0),(a, \dots ,0)),
(0,0)) \}_{a \in \mathbb{C}} \subset 
\mathbf{U}_n^1$.
Therefore
$N_{\mathbf{U}_0^1, \mathcal{H}_{(1:0)} (\mathbf{J}_n);
\mathbf{U}_{n}^1}$
is a point. The calculation for arbitrary
$z \in \mathbb{CP}^1$ is completely analogous, and,
moreover, the variety 
$N_{\mathbf{U}_0^1, \mathcal{H}_z (\mathbf{J}_n);
\mathbf{U}_{n}^1}$ does not depend
on $z$. Thus
$\chi (X) = \chi (\mathbb{CP}^1) = 2$, which is
the statement \ref{KProducts01nn}.

Statement \ref{KProductsnn10} is dual to 
\ref{KProducts01nn}.
\end{proof}

\subsubsection{}\label{RingelKFirst}

Let us recall that $\mathfrak{n}^{\ast} (K)$ denotes
the Lie algebra generated by $E_{(0,1)}$ and
$E_{(1,0)}$ with respect to the bracket
\eqref{ConvolutionBracket}.
In \ref{DefinitionOfXiAffine} we introduced
a homomorphism  $\Xi_K  :
\mathfrak{n}^{\epsilon} (K)
\rightarrow \mathfrak{n}^{\ast} (K)$.
The following proposition is an analog of Ringel Theorem
(Theorem \ref{FiniteTheorem}) for the quiver $K$.

\begin{proposition}\label{KTheoremFirst}
The  homomorphism  $\Xi_K$
is given by the formulas  
\begin{equation}\nonumber
\begin{split}
\Xi_K 
(\tilde{e}_{\alpha}) &= E_{\alpha}
\text{ for $\alpha\in R^{\re}_+$ },
\\
\Xi_K 
( \alpha_1 (n) ) &= (-1)^{n+1} E_{(n,n)}.
\end{split}
\end{equation}
\end{proposition}
\begin{proof}
We prove the proposition by induction. 
It follows from the definitions that
\begin{gather}\nonumber
\Xi_K (\tilde{e}_{(0,1)}) = E_{(0,1)} ,
\\\nonumber
\Xi_K (\tilde{e}_{(1,0)}) = E_{(1,0)}.
\end{gather}
Suppose that
\begin{gather}\nonumber
\Xi_K (\tilde{e}_{(k,k+1)}) = E_{(k,k+1)},
\\\nonumber
\Xi_K (\tilde{e}_{(k+1,k)}) = E_{(k+1,k)}.
\end{gather}
Then using \eqref{KBracket}
and Proposition \ref{KProducts}
we have
\begin{gather}\nonumber
\Xi_K (\alpha_1 (k+1)) =
\Xi_K ((-1)^k [\tilde{e}_{(0,1)}, 
\tilde{e}_{(k+1,k)}])=
\\\nonumber
=(-1)^k [\Xi_K (\tilde{e}_{(0,1)}) , 
\Xi_K (\tilde{e}_{(k+1,k)})] =
(-1)^k [E_{(0,1)} , E_{(k+1,k)}] =
\\\nonumber
=(-1)^{(k+1)+1} E_{(k+1,k+1)}, 
\end{gather}
which gives $\Xi_K (\alpha_1 (k+1))$, and
\begin{gather}\nonumber
\Xi_K (\tilde{e}_{(k+1,k+2)}) =
\Xi_K (\frac{(-1)^{k+1}}{2} 
[\alpha_1 (k+1), \tilde{e}_{(0,1)}])=
\\\nonumber
=\frac{(-1)^{k+1}}{2} 
[\Xi_K (\alpha_1 (k+1)) , 
\Xi_K (\tilde{e}_{(0,1)})] =
\frac{-1}{2} 
[E_{(k+1,k+1)} , E_{(0,1)}] =
\\\nonumber
= E_{(k+1,k+2)}, 
\end{gather}
\begin{gather}\nonumber
\Xi_K (\tilde{e}_{(k+2,k+1)}) =
\Xi_K (\frac{(-1)^{k}}{2} 
[\alpha_1 (k+1), \tilde{e}_{(1,0)}])=
\\\nonumber
=\frac{(-1)^{k}}{2} 
[\Xi_K (\alpha_1 (k+1)) , 
\Xi_K (\tilde{e}_{(1,0)})] =
\frac{1}{2} 
[E_{(k+1,k+1)} , E_{(1,0)}] =
\\\nonumber
= E_{(k+2,k+1)}, 
\end{gather}
which provides the induction step.

One can avoid some of these calculations 
using an argument similar to 
the proof of Theorem \ref{FiniteTheorem}
to get $\Xi_K (\tilde{e}_{\alpha}) = E_{\alpha}$
for $\alpha\in R^{\re}_+$ (note that
the indecomposable object with dimension 
$\alpha\in R^{\re}_+$ can be obtained by repeated 
applications of reflection functors to a 
simple object). However, one would
still need explicit calculations with the 
$\ast$-product
to get the statement of the theorem in the case
of an imaginary root.
\end{proof}

\subsubsection{}\label{ThetaK}

Let us
introduce the following function $\xi_K: 
R_+ \rightarrow \{ \pm 1 \}$.

\begin{equation}\nonumber
\xi_K (\alpha) = (-1)^{(1+\dim_{\mathbb{C}}
\Hom_{\mathcal{M} (K)} 
(\mathbf{P}, \mathbf{P}))}, 
\end{equation}
where $\mathbf{P} \in \Ob (\mathcal{M} (K))$ 
is indecomposable and 
$\dim_{\mathcal{M} (K)} \mathbf{P} =\alpha$.

Since
$\mathbf{U}_n^0$ and $\mathbf{U}_n^1$
can be obtained by sequences
of reflection functors from simple objects
we have
$\xi_K (\alpha) = 1$
if $\alpha \in R_+^{\re}$.
As any indecomposable object with
dimension $n\delta$ is isomorphic to
$\mathcal{H}_z (\mathbf{J}_n)$ for some 
$z\in \mathbb{CP}^1$, it follows from 
Proposition \ref{HomJJ} that
$\xi_K (n\delta) = (-1)^{n+1}$.
In particular, $\xi_K$ is well-defined.

We put $\Tilde{E}_{\alpha} = \xi_K (\alpha) 
E_{\alpha}$.

\subsubsection{}\label{RingelK}

Here is our final theorem for the Kronecker quiver.

\begin{theorem}\label{KTheorem}
The  homomorphism  
$\Xi_K  :
\mathfrak{n}^{\epsilon} (K)
\rightarrow \mathfrak{n}^{\ast} (K)$
is given by the formulas  
\begin{equation}\nonumber
\begin{split}
\Xi_K 
(\tilde{e}_{\alpha}) &= \Tilde{E}_{\alpha}
\text{ for $\alpha\in R^{\re}_+$ },
\\
\Xi_K 
( \alpha_1 (n) ) &= \Tilde{E}_{(n,n)}.
\end{split}
\end{equation}
\end{theorem}
\begin{proof}
Follows from Proposition \ref{KTheoremFirst}.
\end{proof}

\subsubsection{}
Given an imaginary root $n\delta$,
there are infinitely many
indecomposable representations
with dimension $n\delta$
and only one basic element of the
Lie algebra with degree $n\delta$. 
In a sense the Lie algebra of functions
$\mathfrak{n}^{\ast} (K)$ does not distinguish 
among these representations
(they are "similar", though not isomorphic).
Note also that the Euler characteristic of the set
of indecomposable representations with dimension 
$n\delta$
is equal to $2$ (for 
$\mathcal{T}_{n\delta} = \mathbb{CP}^1$). It is this
$2$ that appears in the structure constants 
\ref{KBracket} !

\subsection{Cyclic quiver}
  
\subsubsection{}

Let us fix $N \in \mathbb{Z}$, $N \geq 2$.
In this section we study the cyclic quiver 
$C_N = (\mathbb{Z}/N\mathbb{Z} , 
\mathbb{Z}/N\mathbb{Z}, \In, \Out )$, where
$\In (k) \equiv k+1 \modN$, 
$\Out (k) \equiv k \modN$.
It is a quiver of affine type, 
with the underlying Dynkin graph of type $A^{(1)}_{N-1}$. 

One can draw $C_7$ as follows
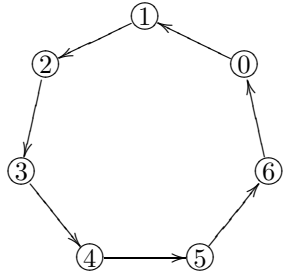
\begin{equation}\nonumber
\def\mypolynode{\ifcase\xypolynode\or
0 \or 1 \or 2 \or 3 \or 4 \or 5 \or 6\fi}
\xypolygon7{~:{/r4pc/:} ~><{} +[o][F-]{\mypolynode}}
\end{equation}

and $C_2$ as follows

\begin{equation}\nonumber
\xygraph{[]
*+[o][F-]{0} :@/^/ [r]*+[o][F-]{1}
*+[o][F-]{1} :@/^/ [l]*+[o][F-]{0}
}
\end{equation}

\subsubsection{}\label{Rotation}

A cyclic quiver has no admissible vertices. 
In particular, there is no Coxeter element. 
However there exists the following "rotation" 
transformation $r$ of the lattice 
$\mathbb{Z}[\mathbb{Z}/N\mathbb{Z}]$.
\begin{gather}\nonumber
r: \mathbb{Z}[\mathbb{Z}/N\mathbb{Z}]
\rightarrow 
\mathbb{Z}[\mathbb{Z}/N\mathbb{Z}] ,
\\\nonumber
r(i) \equiv i+1 \modN .
\end{gather}

\subsubsection{}\label{RootsN3}

It is easy to see that if 
$(a_0, a_1, \dots, a_{N-1} ) \in 
\mathbb{Z}_+ [\mathbb{Z}/N\mathbb{Z}]$
is a positive real root then there exist unique 
$i\in \mathbb{Z}/N\mathbb{Z}$ and 
$l \in \mathbb{Z}_+$,
$l \not\equiv 0 \modN$, 
such that 
$a_k = \# \{ m \in \{ i, \dots , i+l-1\} |\: 
m \equiv k \modN \}$.
We denote the corresponding root by $\alpha_{i,l}$.
In particular the set of simple roots is 
$\{ \alpha_{i,1} \}_{i \in \mathbb{Z}/N\mathbb{Z}}$.
Imaginary roots are integer
multiples of $\delta = (1, 1, \dots , 1 ) \in 
\mathbb{Z} [\mathbb{Z}/N\mathbb{Z}]$.

\subsubsection{}

Next we want to write down explicitly 
the Lie bracket
\eqref{EpsilonCommutatorAffine} in the case
of a cyclic quiver. It appears that expressions
become much simpler if one changes the basis in
$\mathfrak{n}^{\epsilon} (C_N)$ a little.
Namely, let 
\begin{equation}\nonumber
\Tilde{f}_{\alpha_{i,l}}=(-1)^{[l/N]} 
\Tilde{e}_{\alpha_{i,l}},
\end{equation}
where $[l/N]$ denotes the integer part of $l/N$, and
\begin{equation}\nonumber
\Tilde{h} (n) = (-1)^n h (n)
\end{equation}
for $h (n) \in \mathfrak{n}^{\epsilon}_{n\delta}$.
Then
\begin{gather}\label{CNBracket}
[\Tilde{f}_{\alpha_{i,l}}, \Tilde{f}_{\alpha_{j,k}} ] =
\begin{cases}
-\Tilde{f}_{\alpha_{i,k+l}}  
& \text{ if $i+l \equiv j \modN$ } \\
& \text{ and $k+l \not\equiv 0 \modN$ } ,\\  
\Tilde{f}_{\alpha_{j,k+l}}  
& \text{ if $j+k \equiv i \modN$ }\\ 
& \text{ and $k+l \not\equiv 0 \modN$ } ,\\  
(\Tilde{\alpha}_{i,1} + \dots + 
\Tilde{\alpha}_{j-1,1} ) (n)   
& \text{ if $i+l \equiv j \modN$  } \\
& \text{ and $k+l = n N $ } ,
\end{cases}
\\\nonumber
[\Tilde{h} (n), \Tilde{f}_{\alpha_{i,l}} ] =
<h, \alpha_{i,l} >  
\Tilde{f}_{\alpha_{i,l + Nn}} . 
\end{gather}

\subsubsection{}

Let us consider the  case $N=2$. 
The Lie bracket \ref{CNBracket} becomes
\begin{gather}\nonumber
[\Tilde{f}_{\alpha_{i,l}}, \Tilde{f}_{\alpha_{j,k}} ] =
\begin{cases}
\Tilde{\alpha}_{i,1} (n)   
& \text{ if $i \not\equiv j \:(\text{mod}\: 2)$ 
and $k+l =  2n $ } , \\
0 & \text{ otherwise } ,
\end{cases}
\\\nonumber
[\Tilde{h} (n), \Tilde{f}_{\alpha_{i,l}} ] =
<h, \alpha_{i,l} >  
\Tilde{f}_{\alpha_{i,l + 2n}} . 
\end{gather}

The scalar product is given
by \mbox{$<\alpha_{0,1},\alpha_{0,1}>$}$=$
\mbox{$<\alpha_{1,1},\alpha_{1,1}>$}$=$
\mbox{$-<\alpha_{0,1},\alpha_{1,1}>$}$=2$.

Note that the quivers $C_2$ and $K$ have the same 
underlying Dynkin graph $A^{(1)}_1$. Therefore
the algebras $\mathfrak{n}^{\epsilon} (C_2)$ and
$\mathfrak{n}^{\epsilon} (K)$ are 
(non-canonically) isomorphic.
An isomorphism can be given by 
\begin{equation}\nonumber
\eta :\mathfrak{n}^{\epsilon} (C_2) 
\rightarrow \mathfrak{n}^{\epsilon} (K) 
\end{equation}
\begin{equation}\nonumber
\begin{split}
\eta (\Tilde{e}_{\alpha_{0,2n+1}}) &=
\eta ( (-1)^n \Tilde{f}_{\alpha_{0,2n+1}}) = 
(-1)^n \Tilde{e}_{n+1,n} , 
\\
\eta (\Tilde{e}_{\alpha_{1,2n+1}}) &=
\eta ( (-1)^n \Tilde{f}_{\alpha_{1,2n+1}}) = 
(-1)^{n+1} \Tilde{e}_{n,n+1} , 
\\
\eta (\alpha_{1,1} (n)) &=
\eta ((-1)^n \Tilde{\alpha}_{1,1} (n)) = 
\alpha_1 (n) .
\end{split}
\end{equation}

\subsubsection{}

We now turn to the Lie 
algebra $\mathfrak{n}^{\ast} (C_N)$. 
The category  $\mathcal{M} (C_N)$ is
isomorphic to the following category 
$\mathcal{N}_N$. Objects of $\mathcal{N}_N$ are
pairs $(V,x)$, where $V$ is a 
$\mathbb{Z}/N\mathbb{Z}$-graded 
$\mathbb{C}$-linear space and
$x$ is  $\mathbb{C}$-linear 
nilpotent
endomorphism  of $V$, such that 
$\deg (xv) = \deg (v) +1$ for any 
$v\in V$, the set of morphisms from $(V,x)$
to $(W,y)$ is the set of all 
$\mathbb{Z}/N\mathbb{Z}$-graded 
$\mathbb{C}$-linear maps $f$ from $V$ to $W$,
such that $yf=fx$. 

\subsubsection{}

A simple object $\mathbf{P}_{i,1}$ 
($i\in \mathbb{Z}/N\mathbb{Z}$) 
of the category $\mathcal{N}_N$
is a pair 
$\mathbf{P}_{i,1} = (V,x) \in \Ob (\mathcal{N}_N)$,
such that 
\begin{gather}\nonumber
V_i = \mathbb{C} ,
\\\nonumber
V_j = 0 \text{ for $j \neq i$ } ,
\\\nonumber
x=0 .
\end{gather}

The set $\{ [\mathbf{P}_{i,1}] \}_{i\in 
\mathbb{Z}/N\mathbb{Z}}$ is the 
complete set of  
isomorphism classes of
simple objects of $\mathcal{N}_N$.

\subsubsection{}

Given $i \in \mathbb{Z}/N\mathbb{Z}$ and
$l \in \mathbb{Z}_+$, $l \neq 0$,
we denote by $\mathbf{P}_{i,l}$ the following
object  $(V,x)$ of $\mathcal{N}_N$:

\begin{gather}\nonumber
V = \mathbb{C}^l ,
\\\nonumber
\deg (\epsilon_k) \equiv l-k+i \modN ,
\\\nonumber
x=J_l ,
\end{gather}
where 
$\{ \epsilon_1, \dots ,\epsilon_l \}$ 
is the standard basis in $\mathbb{C}^l$ and
$J_l$ is the Jordan block \eqref{Jordan}.

Note that the notation 
$\mathbf{P}_{i,l}$ is consistent with the
notation $\mathbf{P}_{i,1}$ used for
simple objects.

The dimension of $\mathbf{P}_{i,l}$ is given
by 
\begin{equation}\nonumber
\dim_{\mathcal{N}_N} \mathbf{P}_{i,l} = 
\begin{cases}
\alpha_{i,l} & \text{ if $l \not\equiv 0 \modN$ }, 
\\
n\delta & \text{ if $l = nN$ }.
\end{cases}
\end{equation}

\begin{proposition}\label{CIndecomposables}
The set 
$\{ [\mathbf{P}_{i,l}] \}_{\substack{
i \in \mathbb{Z}/N\mathbb{Z} 
\\ l \in \mathbb{Z}_+ \backslash \{ 0 \}}}$
is the complete 
set of isomorphism classes of
indecomposable objects 
of $\mathcal{N}_N$.
In other words, $\mathcal{T}_{\alpha}$ is a one
element set if $\alpha$ is 
a positive real root and 
an $N$ element set if
$\alpha$ is a positive imaginary root.
\end{proposition}

\begin{proof}
An exercise in (graded) linear algebra (see, for
example, \cite[2.19]{Lusztig1992}). 
\end{proof}

\subsubsection{}

We denote by $E_{i,l}
\in \mathcal{L} (C_N)$ the characteristic
function of $[\mathbf{P}_{i,l}]$. 

The following proposition describes  the
restriction of the $\ast$-product to
the set of indecomposable representations.

\begin{proposition}\label{CNProduct}
\begin{equation}\nonumber
(E_{i,l} \ast E_{j,k}) |_{
\mathbf{E}^{ind}_{\alpha_{i,l}
+\alpha_{j,k}}}
= 
\begin{cases}
E_{j,l+k}  &\text{ if $j+k \equiv i \modN$ },\\
0 &\text{ if $j+k \not\equiv i \modN$ }.
\end{cases}
\end{equation}
\end{proposition}
\begin {proof}
Follows from the fact that if 
$\mathbf{X} \in \Ob (\mathcal{N}_N)$
is a subobject of $\mathbf{P}_{m,n}$ then 
$\mathbf{X}$ is isomorphic to $\mathbf{P}_{o,p}$ 
for some
$p \leq n$ and $o = m + n - p \modN$, and
$\mathbf{P}_{m,n} / \mathbf{X}$ is isomorphic to
$\mathbf{P}_{m,n-p}$.
\end{proof}

\subsubsection{}

Now we are ready to prove the analog of 
the Ringel theorem 
for the quiver $C_N$
(cf. \cite{Ringel1993}).

\begin{proposition}\label{CTheoremFirst}
The following formulas describe the map 
$\Xi_{C_N} : \mathfrak{n}^{\epsilon} (C_N)
\rightarrow
\mathfrak{n}^{\ast} (C_N)$:
\begin{equation}\nonumber
\begin{split}
\Xi_{C_N} 
(\Tilde{f}_{\alpha_{i,l}} ) &= E_{i,l} ,\\
\Xi_{C_N} 
(\Tilde{\alpha}_{i,1} (n) ) &= 
E_{i+1,nN} - E_{i,nN} .
\end{split}
\end{equation}
\end{proposition}

\begin{proof}
Let us temporarily denote the map given by the
formulas in the formulation of the proposition by
$\Xi'_{C_N}$. It follows from
\eqref{CNBracket} and Proposition \ref{CNProduct}
that $\Xi'_{C_N}: \mathfrak{n}^{\epsilon} (C_N)
\rightarrow
\mathfrak{n}^{\ast} (C_N)$ is a homomorphism of
Lie algebras, whose value on generators 
$\Tilde{f}_{\alpha_{i,1}}$ of 
$\mathfrak{n}^{\epsilon} (C_N)$
coincides with the value of $\Xi_{C_N}$. Thus
$\Xi'_{C_N} = \Xi_{C_N}$, which proves the proposition.

Another possible proof is an 
induction on length $l$ in $\alpha_{i,l}$ 
similar to the one used in the proof of
Proposition \ref{KTheoremFirst}.
\end{proof}

\subsubsection{}\label{ThetaC}

Let us introduce the following function 
$\xi_{C_N} :
R_+ \rightarrow \{  \pm 1 \}$ (cf. \ref{ThetaK}).

\begin{equation}\nonumber
\xi_{C_N} (\alpha) = (-1)^{(1+\dim_{\mathbb{C}}
\Hom_{\mathcal{N}_N} 
(\mathbf{P}, \mathbf{P}))},
\end{equation} 
where $\mathbf{P} \in \Ob (\mathcal{N}_N)$ 
is indecomposable and 
$\dim_{\mathcal{N}_N} \mathbf{P} =\alpha $.

It follows from a graded analogue of 
Proposition \ref{HomJJ} that
\begin{equation}\nonumber
\begin{split}
\xi_{C_N} (\alpha_{i,l}) &= (-1)^{[l/N]},\\
\xi_{C_N} (n\delta) &= (-1)^{n+1} .
\end{split}
\end{equation}

We put $\Tilde{E}_{i,l} = \xi_{C_N} (\alpha_{i,l}) 
E_{i,l}$ for $l \not\equiv 0 \modN$ and
$\Tilde{E}_{i,nN} = \xi_{C_N} (n\delta) E_{i,nN}$.

\subsubsection{}

Here is our final theorem for a cyclic quiver
(we return to the original generators
$\Tilde{e}_{\alpha}$ and $h(n)$ of
$\mathfrak{n}^{\epsilon} (C_N)$).

\begin{theorem}\label{AffineCTheorem}
The following formulas describe the map 
$\Xi_{C_N} :
\mathfrak{n}^{\epsilon} (C_N)
\rightarrow
\mathfrak{n}^{\ast} (C_N)$:
\begin{equation}\nonumber
\begin{split}
\Xi_{C_N} 
(\Tilde{e}_{\alpha_{i,l}} ) &= \Tilde{E}_{i,l} ,\\
\Xi_{C_N} 
(\alpha_{i,1} (n) ) &= 
\Tilde{E}_{i,nN} - \Tilde{E}_{i+1,nN} .
\end{split}
\end{equation}
\end{theorem}
\begin{proof}
Follows from Proposition \ref{CTheoremFirst}.
\end{proof}

\subsubsection{Remark}
For any real root $\alpha$ there is only one
indecomposable representation of $C_N$ with dimension 
equal to
$\alpha$. However, unlike the case of a finite
type quiver 
the image under the map 
$\Xi_{C_N}$ 
of a basic root vector $\Tilde{e}_{\alpha}$
is
not always the characteristic function
$E_{\alpha}$, but
$\pm E_{\alpha}$ (compare with imaginary
roots of the Kronecker quiver).

\subsubsection{}

We recall (see Proposition \ref{CIndecomposables}) 
that the set of isomorphism classes of indecomposable 
representations with
dimension $n\delta$ is the $N$-element set 
$\mathcal{T}_{n\delta} = \{ [\mathbf{P}_{i,nN}] 
\}_{i \in \mathbb{Z}/N\mathbb{Z}}$. Therefore the set
of functions on $\mathcal{T}_{n\delta}$ is
$N$-dimensional. However
$\dim_\mathbb{C} 
\mathfrak{n}^{\epsilon}_{n\delta} (C_N) = N-1$. 
The following corollary of Theorem 
\ref{AffineCTheorem} describes the image of the 
restriction 
to $\mathfrak{n}^{\epsilon}_{n\delta} (C_N)$
of the map $\Xi_{C_N}$.

\begin{corollary}
$\Xi_{C_N} (\mathfrak{n}^{\epsilon}_{n\delta} (C_N))=
\{ f : \mathcal{T}_{n\delta} 
\rightarrow \mathbb{C} \ | \: 
\sum_{i \in \mathbb{Z}/N\mathbb{Z}}
f ([\mathbf{P}_{i,nN}]) =0 \}$
\end{corollary}
This corollary concludes the discussion of
the cyclic quiver $C_N$. 

\subsection{General affine quiver}

\subsubsection{}

Finally let $Q =(I, \Omega , \In , \Out )$
be a quiver of affine type,
$Q \neq C_N$, $K$. In particular,
$Q$ can be a quiver with cyclic 
(type $A_k^{(1)}$) underlying Dynkin graph, but
with non-cyclic orientation.

In what follows we use classification
of indecomposable objects of $\mathcal{M} (Q)$
as given by V. Dlab and C.M. Ringel 
\cite{DlabRingel}. Let us note that originally
the classification was obtained by 
L.A. Nazarova \cite{Nazarova1973}, and
P. Donovan and M.R. Freislich 
\cite{DonovanFreislich1973}. 
However the description in \cite{DlabRingel}
is more conceptual and suitable for our
purposes. The original proof of Dlab and
Ringel relies on case-by-case consideration.
Since then there appeared proofs 
based on McKay correspondence \cite{Lusztig1992},
and on pure homological algebra
(see, for example, \cite{CrawleyBoevey}).

Our strategy in studying the map
$\Xi_Q$
is to consider different kinds of
roots separately.

\subsubsection{}

Given a root $\alpha$ we denote 
by  $\Supp (\alpha)$ (support of $\alpha$)
the set of all vertices $i\in I$
such that $\alpha_i \neq 0$.

A root $\alpha$ is said to be of \emph{finite type} 
if there exists a proper full subquiver 
$Q' = (I' , \Omega' , \In' \Out' )$
of $Q$ such
that $\Supp (\alpha ) \subset I'$. It follows that
$Q'$ is of finite type and that
$\alpha$ is a real root (see \ref{DynkinFiniteAffine}).

Note that the set $R'_+$ of 
all positive roots of $Q$ with support in $Q'$
coincides with the set of 
positive roots of $Q'$, and
that the restriction of the Euler cocycle 
$\epsilon_Q$ to 
$\mathbb{Z}[I']\times\mathbb{Z}[I']$ coincides with 
$\epsilon_{Q'}$.
Therefore there exists a natural embedding
$\mathfrak{i}^{\epsilon}_{Q'\subset Q} : 
\mathfrak{n}^{\epsilon} (Q') \rightarrow
\mathfrak{n}^{\epsilon} (Q)$.

We recall that $\mathfrak{i}_{Q'\subset Q} : 
\mathfrak{n}^{\ast} (Q') \rightarrow
\mathfrak{n}^{\ast} (Q)$ denotes the Hall
map associated to the embedding 
functor $\mathcal{I}_{Q'\subset Q}:
\mathcal{M} (Q') \rightarrow \mathcal{M} (Q)$
(see \ref{EmbeddingFunctor}).

The following diagram 
of Lie algebra homomorphisms
is commutative
on generators $\{ \Tilde{e}_k \}_{k \in I'}$
of $\mathfrak{n}^{\epsilon} (Q')$,
and, therefore, is commutative.

\begin{equation}\label{AffineFiniteDiagram}
\xymatrix{
\mathfrak{n}^{\epsilon} (Q')
\ar[d]_{\Xi_{Q'}}
\ar[r]^{\mathfrak{i}^{\epsilon}_{Q'\subset Q}} & 
\mathfrak{n}^{\epsilon} (Q)
\ar[d]^{\Xi_Q} \\
\mathfrak{n}^{\ast} (Q') 
\ar[r]^{\mathfrak{i}_{Q'\subset Q}} &
\mathfrak{n}^{\ast} (Q) 
}
\end{equation}

\begin{numproposition}\label{AffineFinite}
Let $\alpha$ be a positive root of finite type.
Then:
\begin{enumerate}
\item
There exists unique (up to an 
isomorphism) indecomposable
object 
$\mathbf{P}_{\alpha} \in 
\Ob (\mathcal{M} (Q))$ with 
$\dim_{\mathcal{M} (Q)}
\mathbf{P}_{\alpha} = \alpha$.
In other words $\mathcal{T}_{\alpha} = pt$. We
denote by $E_{\alpha}\in \mathcal{L} (Q)$ 
the characteristic function
of $\mathcal{T}_{\alpha}$.
\item
$\Xi_Q (\Tilde{e}_{\alpha})
= E_{\alpha}$.
\end{enumerate}
\end{numproposition}
\begin{proof}
The proposition follows from \ref{HallInd},
diagram \eqref{AffineFiniteDiagram}, and Theorem
\ref{FiniteTheorem}.
\end{proof}

\subsubsection{}\label{DefinitionofDefect}

Let $\partial: \mathbb{Z} [I] \rightarrow \mathbb{Z}$
be an additive function given by
\begin{equation}\nonumber
\partial (\alpha) = e (\delta, \alpha),
\end{equation}
where $\delta$ is the first imaginary root
and $e$ is the Euler form \ref{DefinitionOfEpsilon}.

The integer $\partial (\alpha)$ is called 
\emph{defect}
of $\alpha$. 

An element $\alpha \in \mathbb{Z} [I]$
is called \emph{regular} 
(resp. \emph{irregular}) if 
$\partial (\alpha) = 0$ 
(resp. $\partial (\alpha) \neq 0$). 

A root $\alpha \in R$ is called
\emph{regular root} 
(resp. \emph{irregular root})
if it is regular (resp. irregular)
as an element of $\mathbb{Z} [I]$.

We denote the set of regular (resp. irregular)
positive roots by $R_+^{\Reg}$ 
(resp. $R^{\Irr}_+$).

Let us note that each imaginary root is
regular, because 
$e( \delta , \delta ) =
\frac{<\delta , \delta >}{2} = 0$. However there
are regular real roots.

\subsubsection{}

Irregular roots behave similarly
to roots of a quiver of finite type.
The following proposition is proven in
\cite[Sections 1, 2]{DlabRingel}.
\begin{proposition}
Let $Q$ be a quiver of affine type, 
$Q \neq C_N , K$. Then: 
\begin{enumerate}
\item\label{AffineIrregularReflections}
A positive root $\alpha$ is irregular if and only if 
there exists an admissible sequence 
$\{ i_t \}_{t=0}^{k}$
of vertices of $Q$
such that
$\alpha = \sigma_{i_k} 
\sigma_{i_{k - 1}} 
\dots \sigma_{i_{1}} i_{0}$,
and
$\sigma_{i_l} 
\sigma_{i_{l - 1}} 
\dots \sigma_{i_{1}} i_{0} \in R_+$
for any $l \in \{ 1, \dots , k \}$.
\item\label{IrregularIndThenRoot}
If $\mathbf{M} \in \Ob (\mathcal{M} (Q))$
is indecomposable and 
$\dim_{\mathcal{M} (Q)} \mathbf{M}$
is irregular then 
$\dim_{\mathcal{M} (Q)} \mathbf{M}$
is an (irregular) root.
\end{enumerate}
\end{proposition}

\subsubsection{}

Having the admissible sequence of vertices
\ref{AffineIrregularReflections}
one can repeat word-by-word the proof of
Theorem \ref{FiniteTheorem}
to get the following proposition.

\begin{proposition}\label{AffineIrregular}
Let $\alpha \in R^{\Irr}_+$. Then:
\begin{enumerate}
\item\label{AffineIrregularInd}
There exists  unique (up to an isomorphism)
indecomposable object 
$\mathbf{P}_{\alpha} \in \Ob (\mathcal{M} (Q))$
with $\dim_{\mathcal{M} (Q)} 
\mathbf{P}_{\alpha} = \alpha$.
In other words $\mathcal{T}_\alpha = pt$. 
We denote the 
characteristic function 
of $\mathcal{T}_\alpha$ by
$E_{\alpha} \in \mathcal{L} (Q)$.
\item\label{AffineIrregularXi}
$\Xi_Q (\Tilde{e}_{\alpha})=
E_{\alpha} $.
\end{enumerate}
\end{proposition}

\subsubsection{}

Next we turn to regular roots. 

Let us consider a full subcategory
$\mathcal{M}^{\Reg} (Q)$
of $\mathcal{M} (Q)$ with objects being
such $\mathbf{M} \in \Ob (\mathcal{M} (Q))$
that 
\begin{itemize}
\item
$\partial (\dim_{\mathcal{M} (Q)} \mathbf{M}) =0$ 
and
\item
$\partial (\dim_{\mathcal{M} (Q)} 
\mathbf{M}') \leq 0$
for any subobject 
$\mathbf{M}'$ of $\mathbf{M}$ 
in $\mathcal{M} (Q)$.
\end{itemize}

Objects of $\mathcal{M}^{\Reg} (Q)$ are
called \emph{regular} objects
(of $\mathcal{M} (Q)$).

The following is proved in \cite{DlabRingel}.

\begin{proposition}
\begin{enumerate}
\item 
The subcategory $\mathcal{M}^{\Reg} (Q)$
is abelian and \emph{\'epaisse} 
in $\mathcal{M} (Q)$.
\item\label{IndDimReg}
Let $\mathbf{M}\in \Ob (\mathcal{M} (Q))$ be
indecomposable and $\partial 
(\dim_{\mathcal{M} (Q)} \mathbf{M}) =0$.
Then $\mathbf{M} \in 
\Ob (\mathcal{M}^{\Reg} (Q))$.
\end{enumerate}
\end{proposition}

\subsubsection{}

The following proposition gives 
a complete description
of the category $\mathcal{M}^{\Reg}$
due to V. Dlab and C.M. Ringel
\cite{DlabRingel} (see also 
\cite{Lusztig1992}, \cite{CrawleyBoevey}).

\begin{proposition}\label{CProposition}
There exists a set of functors 
$\{ \mathcal{C}_z \}_{z \in \mathbb{CP}^1}$,
parameterized by points of 
$\mathbb{CP}^1$, such that:
\begin{enumerate}
\item\label{CHall}
$\mathcal{C}_z :
\mathcal{M} (C_{N_z})
\rightarrow 
\mathcal{M} (Q)$ is a Hall
functor from the category
of nilpotent representations of 
a cyclic or the Jordan quiver to
the category of nilpotent 
representations of $Q$.
\item\label{RegularCoproduct}
$\mathcal{M}^{\Reg} (Q) = 
\oplus_{z\in \mathbb{CP}^1}
\im \mathcal{C}_z$ 
(coproduct of abelian categories).
\item\label{RegularDim2}
Let $\alpha \in R^{\Reg}_+$. 
Then there exists $z \in \mathbb{CP}^1$
such that $\alpha =\dim\mathcal{C}_z (\alpha')$,
where $\alpha' \in \mathbb{Z}
[\mathbb{Z} / N_z \mathbb{Z}]$.
\item\label{CDelta}
$\dim \mathcal{C}_z (\delta_z) = \delta$
for any $z \in \mathbb{CP}^1$, where
$\delta_z$ (resp. $\delta$) 
is the first imaginary root of 
$C_{N_z}$ (resp. $Q$).
\item\label{CoxeterRotation}
$\dim\mathcal{C}_z \circ r_z = 
c \circ \dim \mathcal{C}_z$ for any 
$z \in \mathbb{CP}^1$, where
$r_z$ is the rotation 
transformation of 
$\mathbb{Z} [\mathbb{Z} / N_z \mathbb{Z}]$
and $c$ is the Coxeter element of $Q$.
\item
$N_z = 1$ for all $z\in \mathbb{CP}^1$ except for
a finite number of points 
(actually no more then $3$). 
We denote the exceptional 
points by $z_1, \dots , z_L$. 
\item\label{NMinusOne}
$\sum_{z\in \mathbb{CP}^1} (N_z - 1) =
\sum_{i=1}^{L} (N_{z_i} - 1) =
|I|-2=\dim_{\mathbb{C}}
\mathfrak{n}^{\epsilon}_{n\delta} (Q) - 1$.
\end{enumerate}
\end{proposition}

We refer the reader to Corollary
\ref{LN} for a description of the
set of integers $L$, 
$\{ N_{z_i} \}_{i=1}^{L}$. We do not
use the exact values of $L$ or
$N_{z_i}$ in this chapter.

An example of explicit 
construction of the functors 
$\mathcal{C}_z$ in $E^{(1)}_6$ case
can be found in \cite[Chapter 11]{AlgebraVIII}.

\subsubsection{}\label{AffineIndFirst}

Let $\alpha_{z,i,l} = 
\dim \mathcal{C}_z (\alpha_{i,l})$
and $\mathbf{P}_{z,i,l} = \mathcal{C}_{z} 
(\mathbf{P}_{i,l})$, where 
$z \in \mathbb{CP}^1$, $i \in 
\mathbb{Z}/N_z\mathbb{Z}$, and
$l \in \mathbb{Z}_+ \backslash \{ 0 \}$ 
(if $N_z=1$ we put $\mathbf{P}_{0,l}=\mathbf{J}_l$;
if $l=nN_z$ we put $\alpha_{z,i,l}=n\delta$). 
Then $\dim_{\mathcal{M} (Q)} \mathbf{P}_{z,i,l}=
\alpha_{z,i,l}$ and
$\mathbf{P}_{z,i,l}$ is indecomposable 
because $\mathcal{C}_{z}$ is a Hall functor.
Moreover it follows from 
\ref{IrregularIndThenRoot},
\ref{AffineIrregularInd},
\ref{IndDimReg},
and \ref{RegularCoproduct}
that 
$\{ [\mathbf{P}_{\alpha}] \}_{\alpha \in R^{\Irr}_+}
\sqcup 
\{ [\mathbf{P}_{z,i,l}] \}_{\substack{
z \in \mathbb{CP}^1 \\
i \in \mathbb{Z}/N_z\mathbb{Z} \\
l \in \mathbb{Z}_+ \backslash \{ 0 \}}}$
is the complete set of isomorphism
classes of indecomposable representations.

\subsubsection{}\label{CCanonical}

The set of functors $\mathcal{C}_z$ satisfying
conditions of Proposition \ref{CProposition}
is defined uniquely up to a
permutation of the set of points of 
$\mathbb{CP}^1$ and
up to automorphisms of categories 
$\mathcal{M} (C_{N_z})$.
More explicitly, one can first identify the set of
isomorphism classes of simple objects of 
$\mathcal{M}^{\Reg} (Q)$.
This set splits into orbits of the Coxeter functor
(product of reflection functors, 
corresponding to the reflections
in the Coxeter element). Each orbit is the set
$\{ [\mathbf{P}_{z,i,1}] \}_{i 
\in \mathbb{Z}/N_z\mathbb{Z}}$
for some $z$. Then
$\Ob (\im \mathcal{C}_z)$ is 
the set of such $\mathbf{M}\in
\Ob (\mathcal{M}^{\Reg} (Q))$ 
that Jordan-H\"{o}lder series 
(in $\mathcal{M}^{\Reg} (Q)$) of
$\mathbf{M}$ contains only objects 
$\mathbf{P}_{z,i,1}$.

\subsubsection{}

The following proposition lists some properties
of the dimension maps $\dim\mathcal{C}_z$.

\begin{proposition}\label{PropertiesOfDimC}
One has
\begin{enumerate}
\item\label{CEpsilonEpsilon0}
$\epsilon_Q 
( \dim \mathcal{C}_z ( \alpha ), 
\dim \mathcal{C}_w ( \beta ) ) = 1$
if $z \neq w$.
\item\label{CEpsilonEpsilon}
$\epsilon_Q 
( \dim \mathcal{C}_z ( \alpha ), 
\dim \mathcal{C}_z ( \beta ) ) = 
\epsilon_{C_{N_z}} ( \alpha  , \beta )$.
\item\label{CFormForm0}
$< \dim \mathcal{C}_z ( \alpha ), 
\dim \mathcal{C}_w ( \beta ) >_Q = 0$
if $z \neq w$.
\item\label{CFormForm}
$< \dim \mathcal{C}_z ( \alpha ), 
\dim \mathcal{C}_z ( \beta ) >_Q = 
<\alpha , \beta>_{C_{N_z}}$.
\item\label{CRegularRoots}
$\{ \alpha_{z,i,l} \}_{\substack{
z \in \mathbb{CP}^1 \\
i \in \mathbb{Z}/N_z\mathbb{Z} \\
l \in \mathbb{Z}_+ \backslash \{ 0 \}
}}$ is
the complete list of regular 
roots, the only equalities being
$\alpha_{z,i,nN_z} = \alpha_{w,j,nN_w}
(=n\delta)$.   
\item\label{DimCAlpha}
Let 
\begin{gather}\nonumber
V=
\mathbb{C} 
[\{ \alpha_{z,i,1} \mod \mathbb{C}\delta 
\}_{\substack{
z \in \mathbb{CP}^1 \\
i \in \mathbb{Z}/N_z\mathbb{Z}
}} ]
= \\\nonumber = 
\mathbb{C} 
[\{ \alpha_{z_k,i,1} \mod \mathbb{C}\delta 
\}_{\substack{
k = 1, \dots , L \\
i \in \mathbb{Z}/N_{z_k} \mathbb{Z}
}} ]
\subset \mathbb{C}[I]/\mathbb{C}\delta,
\end{gather}
where $\mathbb{C}[S]$ denotes the 
$\mathbb{C}$-linear span of a set $S$.
Then 
\begin{equation}\nonumber
\dim_{\mathbb{C}} 
V = |I|-2 = \dim_{\mathbb{C}}
\mathfrak{n}^{\epsilon}_{n\delta} (Q) - 1.
\end{equation}
The following is a  
generating set of linear relations among
elements of the set 
$\{ \alpha_{z_k,i,1} \mod \mathbb{C}\delta 
\}_{\substack{
k = 1, \dots , L \\
i \in \mathbb{Z}/N_{z_k} \mathbb{Z}
}}$.
\begin{equation}\nonumber
\{ \sum_{i \in \mathbb{Z}/N_{z_k}\mathbb{Z}} 
\alpha_{z_k,i,1}
\equiv 0 \: (\mod\: \mathbb{C}\delta )
\}_{k \in \{1, \dots ,L \}}
\end{equation}
\item\label{CFinite}
Let $z=z_k$ (i.e. $N_z \neq 1$). Then 
there exists $i_0 \in \mathbb{Z}/N_z\mathbb{Z}$
such that $\alpha_{z,i_0,1}$ is
of finite type. 
\end{enumerate}
\end{proposition}

\begin{proof}
\ref{CEpsilonEpsilon0} and \ref{CFormForm0} 
follow from \ref{HomExt} and 
\ref{RegularCoproduct}.

\ref{CEpsilonEpsilon} and \ref{CFormForm}
follow from \ref{CHall} and \ref{HallAndEuler}.

\ref{CRegularRoots} follows from the definition 
of root, \ref{RegularCoproduct}, 
\ref{RegularDim2},
\ref{CFormForm0}, and \ref{CFormForm}.

\ref{DimCAlpha} follows from \ref{NMinusOne},
\ref{CFormForm0}, and \ref{CFormForm}.

Let us prove \ref{CFinite}.
Let $p \in I$ be an extending vertex 
(i.e. $\delta_p = 1$).
It follows from the following equality
\begin{equation}\nonumber
\sum_{i\in\mathbb{Z}/N_z\mathbb{Z}}
\alpha_{z,i,1} =
\dim\mathcal{C}_z 
(\sum_{i\in\mathbb{Z}/N_z\mathbb{Z}}
\alpha_{i,1})=
\dim\mathcal{C}_z (\delta_z) =
\delta 
\end{equation}
that either $N_z = 1$ or there exists 
$i_0 \in \mathbb{Z}/N_z\mathbb{Z}$ such
that $(\alpha_{z,i_0,1} )_p=0$, which
proves \ref{CFinite}.
\end{proof}

\subsubsection{}

Now we are ready to give the classification
of indecomposable objects of $\mathcal{M} (Q)$
(cf. \cite{DlabRingel}).

\begin{theorem}\label{AffineInd}
Let $Q$ be a quiver of affine type, 
$Q \neq C_N , K$. 
Let $\mathcal{T}$ 
be the set of isomorphism classes of
indecomposable objects of $\mathcal{M} (Q)$. Then
\begin{enumerate}
\item
$\mathcal{T} = 
\{ [\mathbf{P}_{\alpha}] \}_{\alpha \in R^{\Irr}_+}
\sqcup
\{ [\mathbf{P}_{z,i,l}] \}_{\substack{
z \in \mathbb{CP}^1 \\
i \in \mathbb{Z}/N_z\mathbb{Z} \\
l \in \mathbb{Z}_+ \backslash \{ 0 \}
}}$ ,
\item
$\mathcal{T}_{\alpha} = \emptyset$ 
if $\alpha \notin R_+$ ,
\item
$\mathcal{T}_{\alpha} = pt$ 
if $\alpha \in R^{\re}_+$ ,
\item
$\mathcal{T}_{k\delta} = 
\{ [\mathbf{P}_{z,i,kN_z}] \}_{\substack{
z \in \mathbb{CP}^1 \\
i \in \mathbb{Z}/N_z \mathbb{Z}
}}$ .
\end{enumerate}
\end{theorem}

\begin{proof}
The theorem follows from 
\ref{AffineIndFirst}
and \ref{CRegularRoots}.
\end{proof}

\subsubsection{}

Theorem \ref{AffineInd} allows us 
to introduce the following function $\xi_Q: 
R_+ \rightarrow \{ \pm 1 \}$ (compare with 
\ref{ThetaK},
\ref{ThetaC}).

\begin{equation}\nonumber
\xi_Q (\alpha) = (-1)^{(1+
\dim_{\mathbb{C}}
\Hom_{\mathcal{M} (Q)} 
(\mathbf{P}, \mathbf{P}))}, 
\end{equation}
where 
$\mathbf{P} \in \Ob (\mathcal{M} (Q))$ 
is indecomposable and 
$\dim_{\mathcal{M} (Q)} \mathbf{P} =\alpha $.

The function 
$\xi_Q (\alpha)$ is well-defined
for $\alpha \in R_+^{\re}$ because there is
unique (up to an isomorphism) indecomposable 
object with
dimension $\alpha$.
In the case of an imaginary root one has 
\begin{multline}\nonumber
\xi_Q (n\delta) = (-1)^{(1+
\dim_{\mathbb{C}}\Hom_{\mathcal{M} (Q)} 
(\mathbf{P}_{z,i,nN_z}, \mathbf{P}_{z,i,nN_z}))}=
\\
=(-1)^{(1+\dim_{\mathbb{C}} 
\Hom_{\mathcal{M} (C_{N_z})} 
(\mathbf{P}_{i,nN_z}, \mathbf{P}_{i,nN_z}))}=
(-1)^{n+1}.
\end{multline}
This expression does not depend on $z$ or $i$.

If $\alpha$ is an irregular root
then $\xi_Q (\alpha) = 1$, because
$\mathbf{P}_{\alpha}$ can be obtained from a simple
object by a sequence of reflection functors. 
We put $\Tilde{E}_{\alpha}=
\xi_Q (\alpha) E_{\alpha} = E_{\alpha}$ for
$\alpha\in R_+^{\Irr}$.

\subsubsection{}

We return to the map $\Xi_Q$.

Let $E_{z,j,l}
\in \mathcal{L} (Q)$ be
the characteristic function 
of $[\mathbf{P}_{z,j,l}]
\in \mathcal{T}_{\alpha_{z,j,l}}$ 
and $\Tilde{E}_{z,j,l}=
\xi_Q (\alpha_{z,j,l}) E_{z,j,l}$.

\begin{proposition}\label{CXiFirst}
Let $1 \leq i \leq L$ and 
$j \in \mathbb{Z}/N_{z_i}\mathbb{Z}$. 
Then
\begin{equation}\nonumber
\Xi_Q 
(\Tilde{e}_{\alpha_{z_i,j,1}})=
E_{z_i,j,1} =
\Tilde{E}_{z_i,j,1} .
\end{equation}
\end{proposition}

\begin{proof}
Given $i$ the statement 
is true for one $j$, say $j_0$, 
because of \ref{CFinite} and Proposition 
\ref{AffineFinite}. 
Any other $\alpha_{z_i,j,1}$ can be obtained from
$\alpha_{z_i,j_0,1}$ by repeated applications
of the Coxeter element
(see \ref{CoxeterRotation}). Then one can
use the sequence of reflection functors
corresponding to the sequence of reflections in
the Coxeter element, and
reason similarly to the proof of Theorem
\ref{FiniteTheorem} to get the first
equality in the statement of the proposition
for any $j$. 

The second equality in the statement of the
proposition follows from the fact that
$\xi_Q (\alpha_{z_i,j,1}) =
\xi_{C_{N_{z_i}}} (\alpha_{j,1}) = 1$.
\end{proof}

\subsubsection{}

Next we consider arbitrary
regular positive real roots. 

\begin{proposition}\label{AffineCFrak}
The map $\mathfrak{c}^{\epsilon}_i : 
\mathfrak{n}^{\epsilon} (C_{N_{z_i}})
\rightarrow 
\mathfrak{n}^{\epsilon} (Q)$ given by
\begin{equation}\nonumber
\begin{split}
\mathfrak{c}^{\epsilon}_i 
(\Tilde{e}_{\alpha_{j,l}}) &=
\Tilde{e}_{\alpha_{z_i,j,l}} ,\\
\mathfrak{c}^{\epsilon}_i 
( \alpha_{j,1} (n)) &= 
\alpha_{z_i,j,1} (n)
\end{split}
\end{equation}
is an injective homomorphism of Lie algebras
for any $i \in \{ 1 , \dots , L \}$.
\end{proposition}

\begin{proof}

The proposition follows from 
\ref{CEpsilonEpsilon}, \ref{CFormForm},
and \ref{CDelta}.

\end{proof}

\subsubsection{}

We denote by $\mathfrak{c}_i :
\mathfrak{n}^{\ast} (C_{N_{z_i}}) \rightarrow
\mathfrak{n}^{\ast} (Q)$
the Hall map associated to
the functor $\mathcal{C}_{z_i}:
\mathcal{M} (C_{N_{z_i}}) \rightarrow
\mathcal{M} (Q)$.
In particular, 
$\mathfrak{c}_i (\Tilde{E}_{j,l})=
\Tilde{E}_{z_i,j,l}$.

The following diagram
of Lie algebra homomorphisms

\begin{equation}\nonumber
\xymatrix{
\mathfrak{n}^{\epsilon} (C_{N_{z_i}})
\ar[d]_{\Xi_{C_{N_{z_i}}}}
\ar[r]^{\mathfrak{c}_i^{\epsilon}} & 
\mathfrak{n}^{\epsilon} (Q)
\ar[d]^{\Xi_{Q}} \\
\mathfrak{n}^{\ast} (C_{N_{z_i}}) 
\ar[r]^{\mathfrak{c}_i} &
\mathfrak{n}^{\ast} (Q) 
}
\end{equation}
is commutative on generators 
$\Tilde{e}_{\alpha_{j,1}}$ of
$\mathfrak{n}^{\epsilon} (C_{N_{z_i}})$
(because of Proposition 
\ref{CXiFirst}), and,
therefore, is commutative, and using
Theorem \ref{AffineCTheorem} we
get the following:

\begin{proposition}\label{AffineRegularXi}
Let $1 \leq i \leq L$. Then
\begin{equation}\nonumber
\begin{split}
\Xi_{Q}
(\Tilde{e}_{\alpha_{z_i,j,l}}) =
\Tilde{E}_{z_i,j,l}
\text{ for $l \not\equiv 0 \:
(\text{mod} \: N_{z_i})$ }, \\
\Xi_Q
(\alpha_{z_i,j,1} (n)) = 
\Tilde{E}_{z_i,j,nN_{z_i}}  - 
\Tilde{E}_{z_i,j+1,nN_{z_i}} .
\end{split}
\end{equation}
\end{proposition}

\subsubsection{}

At this point we know values of
$\Xi_Q$ on 
$\mathfrak{n}^{\epsilon}_{\alpha} (Q)$
for any $\alpha \in R^{\re}_+$. As for
an imaginary root space 
$\mathfrak{n}^{\epsilon}_{n\delta} (Q)$
we only know values of 
$\Xi_Q$
on the $\mathbb{C}$-linear span
of 
$\{ \alpha_{z_k,i,1} (n)\}_{\substack{
k = 1, \dots , L \\
i \in \mathbb{Z}/N_{z_k} \mathbb{Z}
}}$, which has codimension $1$ 
in $\mathfrak{n}^{\epsilon}_{n\delta} (Q)$
due to 
Proposition \ref{DimCAlpha}.

Note also that all functions in the
part of the image of $\Xi_Q$
that we know by now vanish
on $[\mathbf{P}_{z,0,n}]$ for 
any $z\neq z_1 , \dots , z_L$ 
(i.e. when $N_z=1$).

\subsubsection{}

The following proposition is due to
V. Dlab and C.M. Ringel \cite{DlabRingel}.

\begin{proposition}\label{KProposition}
There exist
\begin{itemize} 
\item\label{KHall}
a functor 
$\mathcal{K} :
\mathcal{M} (K) 
\rightarrow
\mathcal{M} (Q)$ from the category of
representations of the Kronecker quiver to
the category of representations of $Q$,
\item
a set of elements 
$k_i\in\mathbb{Z}/N_{z_i}\mathbb{Z}$,
one for each $i \in \{ 1 , \dots , L \}$ 
\end{itemize}
such that
\begin{enumerate}
\item
$\mathcal{K}$ is a Hall functor.
\item
$\mathcal{C}_z = \mathcal{K}\circ\mathcal{H}_z$
for any $z\neq z_1 , \dots , z_L$.
\item\label{Intersections} 
$\mathcal{K}\circ\mathcal{H}_{z_i}
(\mathbf{J}_n) = 
\mathbf{P}_{z_i,k_i,nN_{z_i}}$
for any $i \in \{ 1 , \dots , L \}$, 
\item\label{AffineKDim}
Let $\alpha_0 = \dim_{\mathcal{M} (Q)}
\mathcal{K}(\mathbf{U}_0^0)$.
Then $\alpha_0 \in R^{\Irr}_+$ and
$\dim_{\mathcal{M} (Q)}
\mathcal{K}(\mathbf{U}_0^1) =
(\delta-\alpha_0) \in R^{\Irr}_+$.
\end{enumerate}
\end{proposition}

Let us remark that by redefining functors 
$\mathcal{C}_{z_i}$ ("rotating cyclic quivers") 
one can make
$k_i = 0$ for all $i$. 

\subsubsection{}\label{ConstructionOfK}

As the authors could not find a reference for 
Proposition \ref{KProposition} 
except for \cite{DlabRingel}
which relies on case-by-case considerations, let
us give a sketch of a construction of 
the functor $\mathcal{K}$.

We choose an extending vertex $p\in I$, 
which is admissible.
If the underlying Dynkin graph of $Q$ 
is not $A^{(1)}_n$
then any extending vertex has only 
one adjacent edge, and, 
therefore, is admissible. 
If the underlying Dynkin graph
of $Q$ is $A^{(1)}_n$ then one can 
still choose
an admissible extending vertex because $Q$ 
has non-cyclic
orientation.
We assume that $p$ is a sink. 
If $p$ is a source the construction
is analogues or one can
compose the functor $\mathcal{K}$ 
described below with the reflection
functor $\mathcal{S}_p$.

We fix an indecomposable object 
$(\mathbb{C}^{\delta - p} , y )
\in \Ob ( \mathcal{M} (Q) )$
(note that $\delta - p$ is a positive
root of finite type).
Let $q$ ($q_1$, $q_2$ in 
the $A^{(1)}_n$ case) be the vertex 
(vertices) connected
with $p$ by an edge (edges). We denote by $h_0$ 
(resp. $h_1$, $h_2$) the edge (resp. edges)
such that $\In (h_0) = p$ 
(resp. $\In (h_1) = \In (h_2) = p$)
and $\Out (h_0) = q$
(resp. $\Out (h_1) = q_1$, $\Out (h_2) = q_2$).

Instead of $\mathcal{K}$ we construct a functor
$\mathcal{K}' = \mathcal{K} \circ \mathcal{J}: 
\mathcal{C} \rightarrow \mathcal{M} (Q)$
(see \ref{CategoryC} for the definitions of 
$\mathcal{C}$
and $\mathcal{J}$). Then 
$\mathcal{K} = \mathcal{K}' \circ 
\mathcal{J}^{-1}$
(being an isomorphism of categories 
$\mathcal{J}$
has the inverse).

Let $Q$ have the underlying
graph not of $A^{(1)}_n$ type. 
Then it is known that
$( \delta - p )_q = \delta_q = 2$.
The functor $\mathcal{K}' : \mathcal{C} 
\rightarrow 
\mathcal{M} (Q)$ is given by the
following action on objects:
\begin{equation}\nonumber
\mathcal{K}' ((V_0,V_1),\Delta) = (W,z),
\end{equation}
where
\begin{gather}\nonumber
W_i = V_0 \otimes_{\mathbb{C}} 
\mathbb{C}^{( \delta - p )_i} 
\text{ if $i \neq p$ }, 
\\\nonumber
W_p = V_1 , 
\\\nonumber
z_h = \Id_{V_0} \otimes y_h : \:
V_0 \otimes_{\mathbb{C}} 
\mathbb{C}^{( \delta - p )_{\Out (h)}} 
\rightarrow
V_0 \otimes_{\mathbb{C}} 
\mathbb{C}^{( \delta - p )_{\In (h)}}
\text{ if $h \neq h_0$ }, 
\\\nonumber
z_{h_0} = \Delta : \:
V_0 \otimes_{\mathbb{C}} \mathbb{C}^2 ,
\rightarrow V_1 ,
\end{gather}
and the natural action on morphisms.

Let $Q$ have the underlying
graph of $A^{(1)}_n$ type. 
Then 
$( \delta - p )_{q_1} =
( \delta - p )_{q_2} = 1$.
The functor
$\mathcal{K}' : \mathcal{C} \rightarrow 
\mathcal{M} (Q)$ is given by the
following action on objects:
\begin{equation}\nonumber
\mathcal{K}' ((V_0,V_1),\Delta) = (W,z),
\end{equation}
where
\begin{gather}\nonumber
W_i = V_0 \otimes_{\mathbb{C}} 
\mathbb{C}^{( \delta - p )_i}
\text{ if $i \neq p$ }, 
\\\nonumber
W_p = V_1 , 
\\\nonumber
z_h = \Id_{V_0} \otimes y_h : \:
V_0 \otimes_{\mathbb{C}} 
\mathbb{C}^{( \delta - p )_{\Out (h)}} 
\rightarrow
V_0 \otimes_{\mathbb{C}} 
\mathbb{C}^{( \delta - p )_{\In (h)}}
\text{ if $h \neq h_1 , h_2$ }, 
\\\nonumber
z_{h_1} \oplus z_{h_2} = \Delta : \:
V_0 \otimes_{\mathbb{C}} 
(\mathbb{C} \oplus \mathbb{C})
\rightarrow V_1 ,
\end{gather}
and the natural action on morphisms.

One can check that the functor $\mathcal{K}$
defined above satisfies all the conditions in
Proposition \ref{KProposition} (see 
\cite[{\S}9]{CrawleyBoevey} for a similar
argument).

\subsubsection{}

Functor $\mathcal{K} : \mathcal{M} (K) \rightarrow
\mathcal{M} (Q)$ satisfying conditions
of Proposition \ref{KProposition} is not 
unique (cf. \ref{CCanonical}). For example,
in construction \ref{ConstructionOfK} 
one needs to choose an extending vertex, and, moreover,
one can compose the functor $\mathcal{K}$
described in $\ref{ConstructionOfK}$ with
a sequence of reflection functors.

On the other hand a choice of functor $\mathcal{K}$ 
fixes the parameterization of the functors 
$\mathcal{C}_z$ by
points in $\mathbb{CP}^1$. More precisely, 
given a functor
$\mathcal{K} : \mathcal{M} (K) \rightarrow 
\mathcal{M} (Q)$, 
the set of functors $\{ \mathcal{C}_z : 
\mathcal{M} (C_{N_z}) \rightarrow 
\mathcal{M} (Q) \}_{
z \in \mathbb{CP}^1}$, such
that Propositions \ref{CProposition} and 
\ref{KProposition}
hold true is uniquely defined
up to automorphisms of categories 
$\mathcal{M} (C_{N_z})$.

{F}rom now on we fix some particular choice 
of the
functors $\mathcal{K}$ and $\mathcal{C}_z$.

\subsubsection{}\label{KQEpsilon}

Since $\mathcal{K}$ is a Hall functor, 
the following 
proposition follows from
\ref{HallAndEuler}.

\begin{proposition}
A map 
$\mathfrak{k}^{\epsilon} :
\mathfrak{n}^{\epsilon} (K)
\rightarrow
\mathfrak{n}^{\epsilon} (Q)$
given by
\begin{equation}\nonumber
\begin{split}
\mathfrak{k}^{\epsilon} 
(\Tilde{e}_{(n+1,n)}) &= 
\Tilde{e}_{\alpha_0+n\delta} ,\\
\mathfrak{k}^{\epsilon} 
(\Tilde{e}_{(n,n+1)}) &= 
\Tilde{e}_{-\alpha_0+(n+1)\delta} ,\\
\mathfrak{k}^{\epsilon} 
(\alpha'_0 (n) ) &=  \alpha_0 (n)
\end{split}
\end{equation}
is a Lie algebra homomorphism.
Here $\alpha'_0=\dim_{\mathcal{M} (K)} 
\mathbf{U}_0^0 = (1,0)$ and 
$\alpha_0=\dim_{\mathcal{M} (Q)} 
\mathcal{K} (\mathbf{U}_0^0)$.
\end{proposition}

\subsubsection{}

Let $E_{0} (n)$ be the characteristic
function of the set 
$\{ [\mathcal{K}\circ\mathcal{H}_z 
(\mathbf{J}_n)] 
\}_{z \in \mathbb{CP}^1}
\subset \mathcal{T}_{n\delta}$.
In particular, according to 
\ref{Intersections}, 
\begin{equation}\nonumber
E_{0} (n) ([\mathbf{P}_{z_i,j,nN_{z_i}}]) =
\begin{cases}
1 & 
\text{ if $j=k_i$ } ,\\
0 & 
\text{ if $j \neq k_i$  } .
\end{cases}
\end{equation}

We put $\Tilde{E}_0 (n)= 
\xi_Q (n\delta) E_{0} (n)=
(-1)^{n+1} E_{0} (n)$.

\subsubsection{}

Let $\mathfrak{k}:
\mathfrak{n}^{\ast} (K)
\rightarrow
\mathfrak{n}^{\ast} (Q)$
be the Hall map associated with the 
functor $\mathcal K$.

Using \ref{AffineKDim} and \ref{AffineIrregular}
we get the following equalities
\begin{equation}\nonumber
\begin{split}
\Xi_Q 
(\Tilde{e}_{\alpha_0}) &= 
E_{\alpha_0} =
\Tilde{E}_{\alpha_0}    , 
\\
\Xi_Q 
(\Tilde{e}_{\delta - \alpha_0}) &= 
E_{\delta - \alpha_0} =
\Tilde{E}_{\delta - \alpha_0} ,
\end{split}
\end{equation}
which guarantee that the following 
diagram of Lie algebra homomorphisms
is commutative.

\begin{equation}\nonumber
\xymatrix{
\mathfrak{n}^{\epsilon} (K)
\ar[d]^{\Xi_K}
\ar[r]^{\mathfrak{k}^{\epsilon}} & 
\mathfrak{n}^{\epsilon} (Q)
\ar[d]^{\Xi_Q} \\
\mathfrak{n}^{\ast} (K) 
\ar[r]^{\mathfrak{k}} &
\mathfrak{n}^{\ast} (Q) .
}
\end{equation}

Using this diagram and
Theorem \ref{KTheorem} we
get the following proposition.

\begin{proposition}\label{Affine0Xi}
$\Xi_Q
(\alpha_0 (n)) = - \Tilde{E}_{0} (n)$.
\end{proposition}
\subsubsection{}

The last proposition concludes description of the
map $\Xi_Q$ due to the following

\begin{proposition}\label{FullCartan}
The set 
\begin{equation}\nonumber
\{ \alpha_{z_i,j,1} (n) \}_{
\substack{1 \leq i \leq L \\ 
j \in \mathbb{Z}/N_{z_i}\mathbb{Z}}}
\cup \{ \alpha_0 (n) \} \subset
\mathfrak{n}^{\epsilon}_{n\delta} (Q) =
\mathbb{C}[I]/\mathbb{C}\delta
\end{equation}
spans the whole space
$\mathfrak{n}^{\epsilon}_{n\delta} (Q)=
\mathbb{C}[I]/\mathbb{C}\delta$ as a 
$\mathbb{C}$-linear space.
\end{proposition}

\begin{proof}
Note that $\alpha_0$ is linearly independent with
$\{ \alpha_{z_i,j,1} \}_{
\substack{1 \leq i \leq L \\ 
j \in \mathbb{Z}/N_{z_i}\mathbb{Z}}}
\cup \{ \delta \}$
because $\partial (\alpha_0) \neq 0$
whereas $\partial (\alpha_{z_i,j,1})
=\partial{\delta}=0$.
Now the proposition follows from
\ref{DimCAlpha}.
\end{proof}

\subsubsection{}\label{SimplifyNotation}

Let us simplify our notation a little.
We put $N_i=N_{z_i}$, $\Tilde{E}_{i,j} (n) 
= \Tilde{E}_{z_i,j,nN_i}
\in \mathcal{L}_{n\delta} (Q)$, and
$\alpha_{i,j} = \alpha_{z_i, j, 1}
\in \mathbb{Z}[I] \subset \mathbb{C}[I]$.

\subsubsection{}

Here is our main theorem for a general affine
quiver, which follows from Propositions
\ref{AffineIrregular}, \ref{AffineRegularXi},
\ref{CRegularRoots}, \ref{Affine0Xi},
\ref{FullCartan}.

\begin{theorem}\label{AffineTheorem}
Let $Q$ be a quiver of affine type, 
$Q \neq C_N , K$.
Then the following equalities completely
describe the map $\Xi_Q$:
\begin{equation}\nonumber
\begin{split}
\Xi_Q 
(\Tilde{e}_{\alpha}) &= \Tilde{E}_{\alpha} 
\text{ for any $\alpha\in R_+^{\re}$ } , \\
\Xi_Q 
(\alpha_{i,j} (n)) &= 
\Tilde{E}_{i,j} (n) - \Tilde{E}_{i,j+1} (n) ,\\
\Xi_Q (\alpha_{0} (n)) &=  
- \Tilde{E}_{0} (n) .
\end{split}
\end{equation}
\end{theorem}

\subsubsection{}

Theorem \ref{AffineTheorem} is an affine  
analog of the Ringel theorem.
Let us state a few corollaries.

\begin{corollary}
The map $\Xi_Q : \mathfrak{n}^{\epsilon} (Q)
\rightarrow \mathfrak{n}^{\ast} (Q)$ is 
an isomorphism
of Lie algebras.
\end{corollary}
\begin{proof}
The map $\Xi_Q$ is surjective by construction. 
Theorem \ref{AffineTheorem} implies that 
$\Ker \Xi_Q = 0$.
\end{proof}

Since $\Xi_Q$ is an isomorphism given 
by explicit formulas 
(see Theorem \ref{AffineTheorem}) one
can use $\Xi_Q$ together with
\eqref{EpsilonCommutatorAffine}
to obtain explicit expression
for a Lie bracket of any two elements
of $\mathfrak{n}^{\ast} (Q)$
(cf. Corollary \ref{FiniteCorollary}).

\subsubsection{}

The next corollary of Theorem
\ref{AffineTheorem} is a kind of inverse
to Proposition \ref{ReidtmannFinal}. 

\begin{corollary}
The support of the space of functions 
$\mathfrak{n}^{\ast} (Q)$
is equal to the set of all isomorphism classes
of indecomposable 
objects of $\mathcal{M} (Q)$.
In other words, for any indecomposable object
$\mathbf{M}$
there exists an element $f$ of the Lie algebra
$\mathfrak{n}^{\ast} (Q)$
such that $f([\mathbf{M}]) \neq 0$.
\end{corollary}

\subsubsection{}

The most interesting feature of the
category $\mathcal{M} (Q)$ for a
general affine quiver $Q$ is the structure
of the set $\mathcal{T}_{n\delta}$ 
of indecomposable objects
with dimension equal to an imaginary
root (it does not depend on the root in
question) and the restriction of the image
of the map $\Xi_Q$
to this set. We refer the reader to 
\ref{AffineInd} and \ref{AffineTheorem} for
the corresponding results.

Let us consider a map

\begin{equation}\nonumber
\begin{split}
\mu : \mathcal{T}_{n\delta} 
\rightarrow \mathbb{CP}^1 , 
\\
\mu ( [\mathbf{P}_{z,i,nN_z}]) = z   . 
\end{split}
\end{equation}

Note that $\mu^{-1} (z)$ consists of
$N_z$ points. In particular,
the map $\mu$ is
injective on $\mathbb{CP}^1 \backslash 
\{ z_1 , \dots , z_L \}$.

We denote by $\mu_{\ast}$ the 
following map from 
the space of $\mathbb{C}$-valued functions on
$\mathcal{T}_{n\delta}$ to the space
of $\mathbb{C}$-valued functions on
$\mathbb{CP}^1$:

\begin{equation}\nonumber
\mu_* (\phi) (z) = 
\sum_{[\mathbf{P}] \in \mu^{-1}(z)} 
\phi ([\mathbf{P}]) .
\end{equation}

\subsubsection{}\label{NAstIm}

The following is a corollary of 
Theorem \ref{AffineTheorem}.

\begin{corollary}
\begin{equation}\nonumber
\mathfrak{n}^{\ast}_{n\delta} (Q) =
\im (\Xi_Q 
|_{\mathfrak{n}^{\epsilon}_{n\delta} (Q)})
=\{ \phi: \mathcal{T}_{n\delta} 
\rightarrow \mathbb{C} \: | \:
\mu_{\ast} (\phi) 
\text{ is a constant function} \}.
\end{equation}
\end{corollary}

\subsubsection{}

Instead of introducing the functions 
$\Tilde{E}_{\alpha}$
we could change the definition of the Lie algebra
$\mathfrak{n}^{\epsilon}(Q)$. Namely let
$\mathfrak{n}^{\epsilon'}(Q)$ be the Lie algebra 
defined
in the same way as $\mathfrak{n}^{\epsilon}(Q)$ 
(see
\ref{AffineNEpsilon}) but using the following 
cocycle 
\begin{equation}\label{Coboundary}
\epsilon'_Q (\alpha, \beta)=
\epsilon_Q (\alpha, \beta) 
\xi_Q (\alpha +\beta)
(\xi_Q (\alpha))^{-1} 
(\xi_Q (\beta))^{-1}
\end{equation}
instead of the Euler cocycle $\epsilon_Q$.
Here we use arbitrary extension of $\xi_Q$ from 
$R_+$ to $\mathbb{Z}[I]$.
To avoid confusion we denote the generators of
$\mathfrak{n}^{\epsilon'}(Q)$ by
$\Tilde{f}_{\alpha}$ and $\Tilde{h} (n)$
instead of 
$\Tilde{e}_{\alpha}$ and $h (n)$.
Then it follows from Theorem \ref{AffineTheorem} 
that the map
$\Xi'_Q : \mathfrak{n}^{\epsilon'} (Q)\rightarrow 
\mathfrak{n}^{\ast} (Q)$ given by
\begin{equation}\nonumber
\begin{split}
\Xi'_Q 
(\Tilde{f}_{\alpha}) &= E_{\alpha} 
\text{ for any $\alpha\in R_+^{\re}$ } , \\
\Xi'_Q 
(\Tilde{\alpha}_{i,j} (n)) &= 
E_{i,j} (n) - E_{i,j+1} (n) ,\\
\Xi'_Q (\Tilde{\alpha}_{0} (n)) &=  
-E_{0} (n) 
\end{split}
\end{equation}
is an isomorphism of Lie algebras.

Let us remark that equation
\eqref{Coboundary} means that
the cocycles $\epsilon_Q$ and
$\epsilon'_Q$ differ by a 
coboundary, that is belong to
the same cohomology class
in $H^2 (\mathbb{Z}[I], 
\mathbb{Z}/2\mathbb{Z})$.

\subsection{Lie algebras over $\mathbb{Z}$}

\subsubsection{}

Throughout this section 
$Q =(I , \Omega , \In , \Out )$ 
is a quiver of affine type,
$Q \neq C_N, K$.

\subsubsection{}

Let $\alpha \in R_+^{\re}$. We put
$\aleph^{\epsilon}_{\alpha} (Q)=
\mathbb{Z}\Tilde{e}_{\alpha} 
\subset \mathfrak{n}^{\epsilon}_{\alpha} (Q)$.
We denote by
$\aleph^{\epsilon}_{n\delta} (Q)
\subset \mathfrak{n}^{\epsilon}_{n\delta} (Q)$
the lattice additively generated by
$\{ i \mod \mathbb{C}\delta \}_{i \in I}$. 
The lattice $\aleph^{\epsilon}_{n\delta} (Q)$
is isomorphic to
the root lattice $\mathbb{Z}[I']$ of the 
quiver $Q'=(I' , \Omega' , \In' \Out' )$
obtained by removing from $Q$ an extending
vertex and the adjacent edges.

We put $\aleph^{\epsilon} (Q)= 
\oplus_{\alpha \in R_+}
\aleph^{\epsilon}_{\alpha} (Q)$. The
lattice $\aleph^{\epsilon} (Q)$ is a Lie
algebra over $\mathbb{Z}$ with respect 
to the Lie bracket \eqref{EpsilonCommutatorAffine}.

\begin{proposition}\label{AlephEpsilonGenerators}
The Lie algebra $\aleph^{\epsilon} (Q)$
is generated by $\{ \Tilde{e}_i \}_{i \in I}$.
\end{proposition}

\begin{proof}
Induction on $n ( \alpha )$, where 
$n ( \alpha ) = \sum_{i \in I} \alpha_i$
for a root $\alpha$. It is crucial
that $Q$ is simply laced.
\end{proof}

\subsubsection{}

Let $\alpha \in R_+^{\re}$. We denote by
$\aleph^{\ast}_{\alpha} (Q)$ 
the set of all $\mathbb{Z}$-valued
functions on $\mathcal{T}_{\alpha} = pt$. We put
(cf. Corollary \ref{NAstIm})
\begin{equation}\nonumber
\aleph^{\ast}_{n\delta} (Q) =
\{ \phi: \mathcal{T}_{n\delta} 
\rightarrow \mathbb{Z} \: | \:
\mu_{\ast} (\phi) 
\text{ is a constant function} \}.
\end{equation}
Let $\aleph^{\ast} (Q)= 
\oplus_{\alpha \in R_+}
\aleph^{\ast}_{\alpha} (Q)$.
Our goal is to prove that 
$\Xi_Q (\aleph^{\epsilon} (Q)) =
\aleph^{\ast} (Q)$.

\subsubsection{}

Let $\Hat{\aleph}^{\ast} (Q)= 
\oplus_{\alpha \in R_+}
\Hat{\aleph}^{\ast}_{\alpha} (Q)$,
where
$\Hat{\aleph}^{\ast}_{\alpha} (Q)$
is the set of all constructible
$\mathbf{G}_{\alpha}$-invariant
functions $f :
\mathbf{E}^{nil}_{\alpha} 
\rightarrow 
\mathbb{Z}$ such that $f (x) =0$
if $(\mathbb{C}^{\alpha}, x)$ is decomposable.

It follows from the definition of
the $\ast$-product and  
Proposition \ref{Reidtmann}
that $\Hat{\aleph}^{\ast} (Q)$
is a Lie algebra over $\mathbb{Z}$
with respect to the bracket
\eqref{ConvolutionBracket}.

Using definition of the
map $\Xi_Q$ and Proposition 
\ref{AlephEpsilonGenerators}
one gets
$\Xi_Q (\aleph^{\epsilon} (Q)) \subset
\Hat{\aleph}^{\ast} (Q)$.
On the other hand 
it follows from definitions that 
$\Xi_Q (\aleph^{\epsilon} (Q)) \subset
\mathfrak{n}^{\ast}(Q)$
and from Corollary \ref{NAstIm} 
that $\aleph^{\ast}(Q) =
\mathfrak{n}^{\ast}(Q) 
\cap \Hat{\aleph}^{\ast} (Q)$.
Thus we get an inclusion
\begin{equation}\label{Aleph1}
\Xi_Q (\aleph^{\epsilon} (Q)) \subset
\aleph^{\ast} (Q) . 
\end{equation}

It follows from Theorem \ref{AffineTheorem}
that 
\begin{gather}\nonumber
\Tilde{E}_{i,j} (n) - \Tilde{E}_{i,j+1} (n) 
\in \Xi_Q (\aleph^{\epsilon} (Q)) ,
\\\nonumber
\Tilde{E}_{0} (n) 
\in \Xi_Q (\aleph^{\epsilon} (Q)),
\\\nonumber
\Tilde{E}_{\alpha}
\in \Xi_Q (\aleph^{\epsilon} (Q))
\text{ for $\alpha \in R_+^{\re}$ }.
\end{gather}
Since $\{ \Tilde{E}_{i,j} (n) 
- \Tilde{E}_{i,j+1} (n) \}_{\substack{
1 \leq i \leq L \\ 
j \in \mathbb{Z}/N_i\mathbb{Z}
}}
\cup
\{ E_{0} (n) \}$ is a set of additive generators
of $\aleph^{\ast}_{n\delta} (Q)$
and $\Tilde{E}_{\alpha}$ is 
an additive generator of
$\aleph^{\ast}_{\alpha} (Q)$
for $\alpha \in R^{\re}_+$,
we conclude that
\begin{equation}\label{Aleph2}
\aleph^{\ast} (Q) \subset
\Xi_Q (\aleph^{\epsilon} (Q)) .
\end{equation}

Combining \eqref{Aleph1} and \eqref{Aleph2}
we get the following

\begin{proposition}\label{XiAleph}
$\Xi_Q (\aleph^{\epsilon} (Q))
= \aleph^{\ast} (Q)$.
\end{proposition}

It follows, in particular, that
$\aleph^{\ast} (Q)$ is
closed with respect to the
Lie bracket
\eqref{ConvolutionBracket} and is
generated as a Lie algebra over $\mathbb{Z}$
by $\{ E_i \}_{i \in I}$.

\section{Fine structure of 
affine root systems}\label{FineStructure}

\setcounter{subsection}{1}
\setcounter{subsubsection}{0}

In this section we deliberate on
relations of our results
to the structure theory of affine 
Lie algebras. 

\subsubsection{}

Let $Q=(I,\Omega , \In , \Out )$
be a quiver of affine type,
$Q \neq C_N, K$.

The set $\{ \alpha_{i,j} \in \mathbb{Z}[I] \}_{
\substack{1 \leq i \leq L \\ j \in 
\mathbb{Z}/N_i\mathbb{Z}}}$ 
(notation as in \ref{SimplifyNotation})
plays a prominent role in our consideration
of an imaginary root space 
$\mathfrak{n}^{\ast}_{n\delta} (Q)$. 
This set
of elements of the root lattice 
$\mathbb{Z}[I]$ of $Q$
is canonically associated to the quiver.
We call $\alpha_{i,j}$
a \emph{cyclic root}.

Tables of cyclic roots for particular 
quivers of affine type can be found
in \cite[Section 6]{DlabRingel}.

Let us give another description of the set
of cyclic roots.
It follows from \ref{CoxeterRotation} and 
\ref{CRegularRoots} that regular roots have
finite orbits under the action of the Coxeter
element $c$. The converse is also true,
that is any root having finite $c$-orbit is
regular (see \cite[Section 1]{DlabRingel}).

We call a finite $c$-orbit lowest if an element
of this orbit cannot be represented as a sum
of two regular positive roots (this condition
does not depend on the choice of the element
of the orbit).
Let $\mathcal{O}_i=
\{ \alpha_{i,j} \}_{
j \in \mathbb{Z}/N_i\mathbb{Z}}$.
It follows from \ref{CoxeterRotation},
\ref{CRegularRoots} and \ref{DimCAlpha}
that 
$\{ \mathcal{O}_i \}_{
1 \leq i \leq L }$
is the complete set of the lowest $c$-orbits.

Let $R_i \subset \mathbb{Z}[I]$ 
be the additive span of 
elements of an orbit $\mathcal{O}_i$.
It follows from 
Proposition \ref{PropertiesOfDimC} that
$<R_i,R_j>_Q = 0$ if $i\neq j$,
$\oplus_{i=1}^{L} R_i\otimes_{\mathbb{Z}} 
\mathbb{C}$ has 
codimension $1$ in $\mathbb{C}[I]$, and
each $R_i$ is isomorphic
to the root lattice of type $A_{N_i-1}^{(1)}$.
In other words there are $L$ root lattices
of type $A^{(1)}_n$ "hidden" inside
the root lattice of $Q$, and to uncover them one
can use finite orbits of the Coxeter element.

\subsubsection{}
It follows from Proposition
\ref{XiAleph} that the set
$\{ \alpha_{i,j} \mod \mathbb{C}\delta
= \Xi_Q^{-1} (\Tilde{E}_{i,j} (n) - 
\Tilde{E}_{i,j+1} (n)) \}_{\substack{
1 \leq i \leq L \\ 
j \in \mathbb{Z}/N_i\mathbb{Z}
}} \cup
\{ \alpha_{0} \mod \mathbb{C}\delta
= - \Xi_Q^{-1} (\Tilde{E}_{0} (n) )\}$ is a set
of additive generators of 
$\aleph^{\epsilon}_{n\delta} (Q)$.
Since $\delta_p = 1$ for an extending vertex
$p \in I$ it follows that 
if $\alpha$, $\beta \in \mathbb{Z}[I]
\subset \mathbb{C}[I]$
and $\alpha \equiv \beta \:(\text{mod} \: 
\mathbb{C}\delta )$ then
$\alpha \equiv \beta \:(\text{mod} \: 
\mathbb{Z}\delta)$. Combining with 
\ref{DimCAlpha}, \ref{CDelta} and Proposition
\ref{FullCartan} we get the
following

\begin{proposition}\label{ZIPresentation}
The abelian group $\mathbb{Z}[I]$ has the
following presentation:
\begin{equation}\nonumber
< \{ \alpha_{i,j} \}_{\substack{
1 \leq i \leq L \\ 
j \in \mathbb{Z}/N_i\mathbb{Z}
}} \cup
\{ \alpha_{0} \} \cup \{ \delta \} 
\: | \:
\{ \sum_{j \in \mathbb{Z}/N_i\mathbb{Z}} 
\alpha_{i,j} = \delta 
\}_{1 \leq i \leq L } >,
\end{equation}
that is
$\{ \alpha_{i,j} \}_{\substack{
1 \leq i \leq L \\ 
j \in \mathbb{Z}/N_i\mathbb{Z}
}} \cup
\{ \alpha_{0} \} \cup \{ \delta \}$ is a
set of additive generators of $\mathbb{Z}[I]$
and the set
$\{ \sum_{j \in \mathbb{Z}/N_i\mathbb{Z}} 
\alpha_{i,j} = \delta 
\}_{1 \leq i \leq L}$
is a generating set of relations among
the generators.
\end{proposition}

\subsubsection{}

From now till the end of the chapter
we assume that 
$\delta - \alpha_0$ is a simple root,
$\delta - \alpha_0 = p \in I$.
We also assume that $p$ is an extending  vertex.
Then it follows from the definition of 
$\alpha_0$ that $p$ is a sink.

Let $Q'=(I', \Omega' , \In' , \Out')$ be the
quiver obtained from $Q$ by 
removing the vertex $p$
and adjacent edges. 
It is a quiver of finite type. 
The lattice
$\aleph^{\epsilon}_{n\delta} (Q)$ is isomorphic to
$\mathbb{Z}[I']$.

Assume $Q$ is not of $A^{(1)}_n$ type. 
It follows from \ref{CDelta} that
$(\sum_{j \in \mathbb{Z}/N_i\mathbb{Z}}
\alpha_{i,j})_p = \delta_p = 1$. Hence 
for each $i \in \{ 1 , \dots , L \}$
there exists $n_i \in \mathbb{Z}/N_i\mathbb{Z}$
such that
\begin{gather}\nonumber
(\alpha_{i,n_i})_p = 1,
\\\nonumber
(\alpha_{i,l})_p = 0
\text{ if $l \neq n_i$ }. 
\end{gather}

Note that $p$ being a sink 
we may assume that 
the reflection $\sigma_p$ 
is the first reflection in the Coxeter 
element $c$. Therefore
\begin{equation}\nonumber
0=(\alpha_{i,n_i + 1})_p =
(c \alpha_{i,n_i})_p =
(\alpha_{i,n_i})_q
-(\alpha_{i,n_i})_p =
(\alpha_{i,n_i})_q - 1,
\end{equation}
where $q$ is the vertex connected with
$p$ by an edge. We conclude that
$(\alpha_{i,n_i})_q = 1$. Note that
$(\sum_{j \in \mathbb{Z}/N_i\mathbb{Z}}
\alpha_{i,j})_q = \delta_q = 2$. 
Hence there exists $m_i \in \mathbb{Z}/N_i\mathbb{Z}$,
$m_i \neq n_i$ such that
\begin{gather}\nonumber
(\alpha_{i,m_i})_q = (\alpha_{i,n_i})_q = 1 
\\\nonumber
(\alpha_{i,l})_q = 0
\text{ if $l \neq n_i$, $m_i$ }. 
\end{gather}
Now we again use the fact that $\sigma_p$
is the first reflection in the Coxeter element
$c$:
\begin{equation}\nonumber
(\alpha_{i,m_i + 1})_p =
(c \alpha_{i,m_i})_p =
(\alpha_{i,m_i})_q
-(\alpha_{i,m_i})_p = 1.
\end{equation}
Therefore $n_i = m_i + 1$. 

It follows from the above considerations and from
Proposition \ref{ZIPresentation}
that the abelian 
group $\mathbb{Z}[I']$ is 
freely generated by the set
$\{ \alpha_{i,j} \}_{\substack{
1 \leq i \leq L \\ 
j \in \mathbb{Z}/N_i\mathbb{Z}, 
j \neq n_i 
}} \cup
\{ - \alpha_{0} \}$
and that
\begin{equation}\label{EPAlpha}
\begin{split}
e_Q (-\alpha_0, \alpha_{i,m_i}) = 
e_Q (p, \alpha_{i,m_i}) &= 0,
\\
e_Q (\alpha_{i,m_i}, -\alpha_0)  = 
e_Q (\alpha_{i,m_i}, p) &= -1,
\\
e_Q (-\alpha_0, \alpha_{i,j})  = 
e_Q (p, \alpha_{i,j}) &= 0 
\text{ if $j \neq n_i$, $m_i$ },
\\
e_Q (\alpha_{i,j} , -\alpha_0)  = 
e_Q (\alpha_{i,j} , p ) &= 0 
\text{ if $j \neq n_i$, $m_i$ }.
\end{split}
\end{equation}

In the $A^{(1)}_n$ case one also
has a set of pairs 
$\{ ( n_i , m_i ) \}_{i = 1}^{L}$
such that the above statement is true. 
The proof is similar, except that
one should consider two vertices $q_1$ and $q_2$
connected to $p$ instead of one vertex $q$.
{F}rom now on we do not distinguish 
between $A^{(1)}_n$ and non-$A^{(1)}_n$
cases.

Let $\Hat{Q} = (
\Hat{I} , \Hat{\Omega}, \Hat{\In} , \Hat{\Out})$ 
be a quiver given as follows.
The set of vertices $\Hat{I}$ is equal to 
$\{ \Hat{\alpha}_{i,j} \}_{\substack{
1 \leq i \leq L \\ 
j \in \mathbb{Z}/N_i\mathbb{Z}, 
j \neq n_i 
}} \cup
\{ \spadesuit \}$, 
the set of edges 
is equal to the set of all ordered pairs
$(h,t)$, where either $h = \Hat{\alpha}_{i,j}$
and $t=\Hat{\alpha}_{i,j-1}$, or
$h=\spadesuit$ and $t=\Hat{\alpha}_{i,m_i }$,
and the maps $\Hat{\In}$ and $\Hat{\Out}$ are
given by $\Hat{\In}((h,t))=h$, $\Hat{\Out}((h,t))=t$.
One can draw the quiver $\Hat{Q}$
as follows.

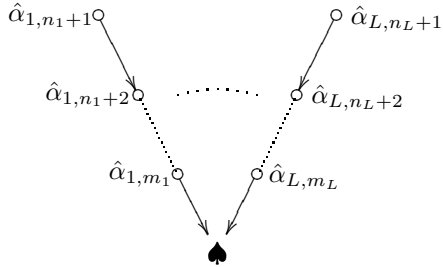
\begin{equation}\nonumber
\xygraph{
!{<0pt,0pt>;<15pt,0pt>:}
[] \spadesuit
-@{<-{}} [uul]
*\cir<2pt>{}
!{\save -<14pt,0pt>*\txt{$\Hat{\alpha}_{1,m_1}$} 
\restore}
-@{{}.{}} [uul]
*\cir<2pt>{}
!{\save -<19pt,0pt>*\txt{$\Hat{\alpha}_{1,n_1+2}$} 
\restore}
-@{<-{}} [uul]
*\cir<2pt>{}
!{\save -<18pt,0pt>*\txt{$\Hat{\alpha}_{1,n_1+1}$} 
\restore}
[ddddddrrr]
\spadesuit
-@{<-{}} [uur]
*\cir<2pt>{}
!{\save +<18pt,-1pt>*\txt{$\Hat{\alpha}_{L,m_L}$} 
\restore}
-@{{}.{}} [uur]
*\cir<2pt>{}
!{\save +<23pt,-1pt>*\txt{$\Hat{\alpha}_{L,n_L+2}$} 
\restore}
-@{<-{}} [uur]
*\cir<2pt>{}
!{\save +<23pt,-1pt>*\txt{$\Hat{\alpha}_{L,n_L+1}$} 
\restore}
[ddddddlll]
[uuuul]
-@/^2pt/@{{}.{}} [rr]
}
\end{equation}

We summarize the results of this subsection in
the following proposition.

\begin{proposition}\label{Nu}
The map 
$\nu : \mathbb{Z} [\Hat{I}]
\rightarrow \mathbb{Z} [I']$
given by 
$\nu (\Hat{\alpha}_{i,j}) = \alpha_{i,j}$,
$\nu (\spadesuit) = - \alpha_{0}$
is an isomorphism of lattices
and
$e_{Q'} (\nu (\alpha ), \nu (\beta ))
= e_{\Hat{Q}} (\alpha, \beta)$ for
any $\alpha$, $\beta \in 
\mathbb{Z}[\Hat{I}]$. 
\end{proposition}

\subsubsection{}

It follows from Proposition \ref{Nu} that 
$<\nu (\alpha ), \nu (\beta )>_{Q'}
= <\alpha, \beta >_{\Hat{Q}}$ for
any $\alpha$, $\beta \in 
\mathbb{Z}[\Hat{I}]$.
In other words, $\nu$ is an isometry of lattices. 

We recall (see
\cite{Bourbaki}) that the set of
simple roots in a 
root system of finite type is unique up to
an action of the semi-direct product
of the automorphism group of the
Dynkin graph and the Weyl group.
In our situation we state it as follows
(abusing notation we denote by the
same letter $\pi$ a morphism of Dynkin
graphs 
$\pi: (I' , \Omega', \{ \In' , \Out' \})
\rightarrow
(\Hat{I} , \Hat{\Omega}, 
\{ \Hat{\In} , \Hat{\Out} \})$
and the induced morphism of
the root lattices
$\pi: \mathbb{Z}[I']
\rightarrow
\mathbb{Z}[\Hat{I}]$).
\begin{proposition}
There exist a morphism of Dynkin graphs 
$\pi: (I' , \Omega', \{ \In' , \Out' \})
\rightarrow
(\Hat{I} , \Hat{\Omega}, 
\{ \Hat{\In} , \Hat{\Out} \})$ and
an element $w$ of the Weyl group 
of 
$(\Hat{I} , \Hat{\Omega}, 
\{ \Hat{\In} , \Hat{\Out} \})$
such that 
$\pi \circ \nu = w$, or, more
explicitly,
\begin{equation}\nonumber
\pi (\alpha_{i,j}) = w (\Hat{\alpha}_{i,j}),
\end{equation}
for any $1 \leq i \leq L$, 
$j \in \mathbb{Z}/N_i\mathbb{Z}, 
j \neq n_i$
and 
\begin{equation}\nonumber
\pi (\alpha_{0}) = - w (\spadesuit ) .
\end{equation}
\end{proposition}

A similar result for such an orientation
of $Q$ that each vertex is
admissible has been proven by
R. Steinberg \cite{Steinberg}.

In particular, the underlying Dynkin graphs of
$\Hat{Q}$ and $Q'$ coincide,
and we get the following corollary
which describes the set of numbers $L$, $\{ N_i 
\}_{i=1}^{L}$
for a quiver with a given underlying Dynkin graph
(cf. \cite{DlabRingel}).

\begin{corollary}\label{LN}
One has:
\begin{enumerate}
\item
If $Q$ is of $A^{(1)}_n$ type
then either $L=1$ and 
$N_1=n$, or $L=2$ and
$N_1 + N_2 = n+1$.
\item
If $Q$ is of $D^{(1)}_n$ type
then $L=3$, 
$N_1=N_2=2$, and $N_3=n-2$
(up to a permutation of $\{ N_i \}$).
\item
If $Q$ is of $E^{(1)}_n$ type $(6 \leq n \leq 8)$ 
then $L=3$, 
$N_1=2$, $N_2=3$, and $N_3=n-3$
(up to a permutation of $\{ N_i \}$).
\end{enumerate}
\end{corollary}

Corollary \ref{LN} does not give the numbers $L$, 
$\{ N_i \}_{i=1}^L$ in the
$A^{(1)}_n$ case. Let us describe these numbers 
for the sake of 
completeness (see \cite{DlabRingel} for the proof).
We denote by $J_1$ (resp $J_2$) the number of 
clockwise
(resp. counterclockwise) oriented arrows in $Q$.
Note that $J_1, J_2\neq 0$ since $Q$ is not cyclic.
If $J_1=1$ (resp. $J_2=1$), then $L=1$ and $N_1=J_2$
(resp. $N_1=J_1$). If $J_1\neq 1$ and $J_2\neq 1$, 
then   
$L=2$ and $N_i=J_i$ up to the transposition of 
$N_1$ and $N_2$.

\subsubsection{}

Finally we would like to discuss relation of
the quiver construction of the Lie algebra 
$\mathfrak{n}$ and
the representation theory of finite groups.
In this section we assume that $Q$ is a quiver of 
$A_{2n+1}^{(1)}$, $D_{n}^{(1)}$, or
$E_{n}^{(1)}$ type with such an orientation
that each vertex is admissible.
Lusztig \cite{Lusztig1992}
has reinterpreted representations of such quivers
using McKay correspondence \cite{McKay}. 

Let us recall that McKay correspondence 
establishes an isomorphism of lattices 
$\eta: K\mathbb{C}\Gamma \rightarrow
\mathbb{Z}[I]$, where $K\mathbb{C}\Gamma$
is the Grothendieck group of the category
of $\mathbb{C}$-linear finite dimensional 
representations of a finite
subgroup $\Gamma$ of $SL (2,\mathbb{C})$ 
and $\mathbb{Z}[I]$ is the root lattice
of an affine Dynkin graph 
$(I, E, \Ends )$.
In particular, $\eta^{-1} (i)$ is the
class of a simple
$\mathbb{C}\Gamma$-module, which we denote
by $\rho_i$.
The following formula describes
the pull-back of the bilinear form
$<,>$ under $\eta^*$.
 
\begin{equation}\nonumber
<\eta (X),\eta (Y)> =
\dim_{\mathbb{C}} \Hom_{\mathbb{C}\Gamma} 
(X \otimes \mathbb{C}^2 , Y) -
\dim_{\mathbb{C}} \Hom_{\mathbb{C}\Gamma} 
(X \otimes \rho , Y),
\end{equation}
where $X, Y$ are $\mathbb{C}\Gamma$-modules,
$\mathbb{C}^2$ is equipped with
the trivial $\Gamma$-action, and
$\rho$ is the natural $2$-dimensional
representation of $\Gamma \subset
SL (2,\mathbb{C})$.

The isomorphism of lattices yields the 
following bijection
between the set of affine Dynkin graphs and the
set of finite subgroups of $SL (2, \mathbb{C})$. 

$A^{(1)}_{n} \leftrightarrow 
\mathbb{Z}/(n+1)\mathbb{Z}$ -
cyclic group of order $n+1$,

$D^{(1)}_{n} \leftrightarrow 
\mathbb{D}_{n-2}$ - binary
dihedral group,

$E^{(1)}_{6} \leftrightarrow 
\mathbb{T}$ - binary
tetrahedral group,

$E^{(1)}_{7} \leftrightarrow 
\mathbb{O}$ - binary
octahedral group,

$E^{(1)}_{8} \leftrightarrow 
\mathbb{I}$ - binary
icosahedral group.

Note that affine Dynkin graphs of type 
$A^{(1)}_{2n}$ do not have an orientation
such that each vertex is admissible and
are excluded from our considerations.
The rest of the graphs admit precisely two
such orientations.

For an affine quiver $Q$ with the special 
orientation Lusztig was able to re-obtain the 
classification of indecomposable
representations of $Q$ entirely in terms of 
the representation theory of $\Gamma$.
Below we discuss relation of
maximal cyclic subgroups of $\Gamma$
and cyclic roots of $Q$ following
\cite{Lusztig1992}.

Let $F$ be the set of lines $l$ in 
$P ( \rho )$ whose isotropy group
$\Gamma_{l}$
(which is a cyclic group)
has order greater than $2$.
The set $F$ is finite. Let us
choose representatives 
$\{ l_i \}_{i = 1}^{L'}$
for $\Gamma$-orbits in $F$.
We assume that if $\Gamma_{l'}=\Gamma_{l_i}$
for some $i$, but
$l'$ and $l_i$ are not in the same $\Gamma$-orbit
then $l' = l_j$ for some $j$.
Let $\Gamma_i = \Gamma_{l_i}$. It follows from 
our assumptions on the quiver $Q$ that
the order of $\Gamma_i$ is even.
We put $N'_i= \frac{| \Gamma_i |}{2}$.
Note that $l_i$ is a one-dimensional
$\mathbb{C}\Gamma_i$-module, and 
that $l_i^{\otimes (2N'_i)}$ is a trivial 
$\mathbb{C}\Gamma_i$-module. Let $\varkappa_{i,j}=
l_i^{\otimes(2j)}\oplus l_i^{\otimes(2j+1)}$
for $j \in \mathbb{Z}/N'_i\mathbb{Z}$. 

One can identify the number of exceptional points
of $\mathbb{CP}^1$
and their multiplicities appearing in the
classification of indecomposable 
representations of $Q$ in terms of the group
$\Gamma$ and its maximal cyclic subgroups
(see \cite[Section 1]{Lusztig1992}).
We recall that $\{ \alpha_{i,j} \in \mathbb{Z}[I] 
\}_{\substack{1 \leq i \leq L \\ j \in 
\mathbb{Z}/N_i\mathbb{Z}}}$
is the set of cyclic roots.

\begin{proposition}\label{LLNNaa}
One has
\begin{enumerate}
\item
$L=L'$.
\item
$N_i=N'_i$ up to a permutation of indexes 
$i \in \{ 1, \dots , L \}$.
\item
$\alpha_{i,j} = \eta ( \Ind_{\Gamma_i}^{\Gamma}
\varkappa_{i,j})$ up to a permutation of indexes
$i \in \{ 1, \dots , L \}$ and a 
cyclic permutation of indexes
$j \in \mathbb{Z}/N_i\mathbb{Z}$.
\end{enumerate}
\end{proposition}

In particular, combining Propositions 
\ref{LN} and \ref{LLNNaa} 
one can identify the lengths of the branches
of the Dynkin graph with the orders of the 
maximal cyclic subgroups of $\Gamma$.

Using Frobenius Reciprocity we get
the following formula for the cyclic
roots in terms of the representation theory
of $\Gamma$ and its cyclic subgroups.

\begin{equation}\nonumber
(\alpha_{i,j})_k =
\dim_{\mathbb{C}}\Hom_{\mathbb{C}\Gamma}
(\Ind_{\Gamma_i}^{\Gamma}\varkappa_{i,j},
\rho_k)=
\dim_{\mathbb{C}}\Hom_{\mathbb{C}\Gamma_i}
(\varkappa_{i,j},\rho_k).
\end{equation}

We would like to remark that both 
our construction of the Lie algebra 
$\mathfrak{n}$  and the proof of
Theorem \ref{AffineTheorem} can be reformulated 
in the language of the
representation theory of $\Gamma$. In particular,
instead of the reflection functors which involve
all possible quiver orientations one might
use ``the square root'' of the Coxeter functor 
and its inverse employed by Lusztig 
\cite{Lusztig1992}.

\end{document}